\documentclass[12pt,a4paper]{article}
\usepackage[utf8]{inputenc}
\usepackage{amssymb,geometry,nicefrac,enumerate,mathtools,comment,bbm,mathrsfs,color,verbatim,graphicx,bm,dsfont,xparse,textcomp}
\usepackage{amsmath}
\allowdisplaybreaks
\usepackage{listings}
\usepackage{amsfonts}
\usepackage{cite}
\usepackage[english]{babel}
\usepackage{subfig}
\usepackage{datetime}
\usepackage[inline]{enumitem}
\usepackage{subcaption}
\usepackage{footnote}
\usepackage[most]{tcolorbox}
\usepackage{colortbl}
\usepackage[table]{xcolor}

\usepackage{etoolbox}
\usepackage{amsthm}
\usepackage{thmtools}
\usepackage[hypertexnames=false]{hyperref}
\hypersetup{
	breaklinks=true,
    colorlinks,
    linkcolor={blue!90!black},  
    citecolor={green},
    urlcolor={green!60!black}   
}
\usepackage[capitalise,nameinlink,sort]{cleveref}
\crefname{figure}{Figure}{Figures}
\Crefname{figure}{Figure}{Figures}
\crefname{enumi}{item}{items}
\crefname{equation}{}{}
\crefname{subsection}{Subsection}{Subsections}
\Crefname{subsection}{Subsection}{Subsections}
\theoremstyle{plain}
\newtheorem{theorem}{Theorem}[section]
\newtheorem{lemma}[theorem]{Lemma}
\newtheorem{remark}[theorem]{Remark}

\newtheorem{definition}[theorem]{Definition}
\newtheorem{corollary}[theorem]{Corollary}
\newtheorem{proposition}[theorem]{Proposition}
\newtheorem{setting}[theorem]{Setting}

\numberwithin{equation}{section}

\DeclareFontEncoding{LS1}{}{}
\DeclareFontSubstitution{LS1}{stix}{m}{n}
\DeclareMathAlphabet{\mathscr}{LS1}{stixscr}{m}{n}

\usepackage{mleftright}
\usepackage{mathtools}

\DeclarePairedDelimiter{\rbr}{(}{)}
\newcommand{\pa}[1]{\left({#1}\right)}
\newcommand{\pb}[1]{\left[{#1}\right]}

\newcommand{\R}{\mathbb{R}}
\newcommand{\N}{\mathbb{N}}

\newcommand{\smallsum}{\textstyle\sum}

\newcommand{\norm}[1]{\lVert #1 \rVert}

\newcommand{\normmm}[1]{{\left\vert\kern-0.25ex\left\vert\kern-0.25ex\left\vert #1
    \right\vert\kern-0.25ex\right\vert\kern-0.25ex\right\vert}}

\newcommand{\abs}[1]{\lvert #1\rvert}
\newcommand{\Abs}[1]{\big| #1 \big|}

\newcommand{\babs}[1]{\bigl\lvert #1\bigr\rvert}

\newcommand{\scalprod}[2]{\langle #1, #2 \rangle}

\newcommand{\eps}{\varepsilon}

\newcommand{\mf}[1]{\mathfrak{#1}}

\newcommand{\cB}{\mathcal{B}}

\newcommand{\cF}{\mathcal{F}}

\newcommand{\cM}{\mathcal{M}}
\newcommand{\cN}{\mathcal{N}}

\newcommand{\cR}{\mathcal{R}}

\newcommand{\fA}{\mathfrak{A}}
\newcommand{\fB}{\mathfrak{B}}

\newcommand{\fH}{\mathfrak{H}}

\newcommand{\fN}{\mathfrak{N}}

\newcommand{\fP}{\mathfrak{P}}

\newcommand{\fV}{\mathfrak{V}}

\newcommand{\fa}{\mathfrak{a}}
\newcommand{\fb}{\mathfrak{b}}
\newcommand{\fc}{\mathfrak{c}}

\newcommand{\ff}{\mathfrak{f}}
\newcommand{\fg}{\mathfrak{g}}

\newcommand{\fj}{\mathfrak{j}}

\newcommand{\fn}{\mathfrak{n}}

\newcommand{\fq}{\mathfrak{q}}

\newcommand{\fv}{\mathfrak{v}}
\newcommand{\fw}{\mathfrak{w}}
\newcommand{\fx}{\mathfrak{x}}
\newcommand{\fy}{\mathfrak{y}}

\newcommand{\bfM}{\mathbf{M}}

\newcommand{\scrP}{\mathscr{P}}

\newcommand{\scrU}{\mathscr{U}}

\newcommand{\scra}{\mathscr{a}}
\newcommand{\scrb}{\mathscr{b}}

\newcommand{\scrd}{\mathscr{d}}

\newcommand{\scrf}{\mathscr{f}}

\newcommand{\scrq}{\mathscr{q}}
\newcommand{\scrr}{\mathscr{r}}

\newcommand{\scrv}{\mathscr{v}}

\newcommand{\scrx}{\mathscr{x}}

\NewDocumentEnvironment{cproof}{m}
{\phantomsection\begin{proof}[Proof of \cref{#1}]}%
{\noindent The proof of \cref{#1} is thus complete.
\end{proof}}

\NewDocumentEnvironment{cproof2}{m}
{\phantomsection\begin{proof}[Proof of \cref{#1}]}%
{\noindent This completes the proof of \cref{#1}.
\end{proof}}

\ExplSyntaxOn

\NewDocumentCommand{\enum}{ O{;} m o }
 {
  \my_enum:nnn { #1 } { #2 } { #3 }
 }

\seq_new:N \l__my_enum_seq
\tl_new:N \l__my_enum_item_tl
\prop_const_from_keyval:Nn \l__verbs
{
show = shows ,
imply = implies ,
demonstrate = demonstrates ,
prove = proves ,
establish = establishes ,
ensure = ensures ,
assure = assures
}
\int_new:N \l__number_of_args

\cs_new_protected:Nn \my_enum:nnn
 {
  \seq_set_split:Nnn \l__my_enum_seq { #1 } { #2 }
  \seq_remove_all:Nn \l__my_enum_seq {}
  \int_set_eq:NN \l__number_of_args { \seq_count:N \l__my_enum_seq }
  \seq_use:Nnnn \l__my_enum_seq { ~and~ } { ,~ } { ,~and~ }
  \IfNoValueTF{#3}{}{
    \space
    \int_compare:nNnTF{ \l__number_of_args } < {2}{ \prop_item:Nn \l__verbs {#3} }{ #3 }
  }
 }

\seq_new:N \g_cflist_loaded
\seq_new:N \g_cflist_pending

\NewDocumentCommand{\cfadd}{ m }
{
  \seq_if_in:NnF \g_cflist_loaded { #1 } {
    \seq_if_in:NnF \g_cflist_pending { #1 } {
      \seq_gput_right:Nn \g_cflist_pending { #1 }
    }
  }
}

\NewDocumentCommand{\cfconsiderloaded}{ m }{
  \seq_gput_right:Nn \g_cflist_loaded {#1}
}

\NewDocumentCommand{\cfremove}{ m }
{
  \seq_gremove_all:Nn \g_cflist_pending { #1 }
}

\NewDocumentCommand{\cfload}{ o }
{
  \seq_if_empty:NTF \g_cflist_pending {\unskip} {
    (cf.\ \cref{\seq_use:Nn \g_cflist_pending {,}})\IfValueTF{#1}{#1~}{\unskip}
    \seq_gconcat:NNN \g_cflist_loaded \g_cflist_loaded \g_cflist_pending
    \seq_gclear:N \g_cflist_pending
  }
}

\NewDocumentCommand{\cfclear} {} {
  \seq_gclear:N \g_cflist_loaded
  \seq_gclear:N \g_cflist_pending
}

\NewDocumentCommand{\cfout}{ o }
{
  \seq_if_empty:NTF \g_cflist_pending {\unskip} {
    (cf.\ \cref{\seq_use:Nn \g_cflist_pending {,}})\IfValueTF{#1}{#1~}{\unskip}
    \seq_gclear:N \g_cflist_pending
  }
}

\NewDocumentCommand{\ifnocf} { m } {
  \seq_if_empty:NT \g_cflist_pending { #1 }
}

\ExplSyntaxOff

\ExplSyntaxOn

\bool_new:N \g_noteobserve

\NewDocumentCommand{\nobs}{}{
  \bool_if:nTF { \g_noteobserve } {
    \bool_gset_false:N \g_noteobserve
    note~
  } {
    \bool_gset_true:N \g_noteobserve
    observe~
  }
}

\NewDocumentCommand{\Nobs}{}{
  \bool_if:nTF { \g_noteobserve } {
    \bool_gset_false:N \g_noteobserve
    Note~
  } {
    \bool_gset_true:N \g_noteobserve
    Observe~
  }
}

\ExplSyntaxOff

\ExplSyntaxOn

\bool_new:N \g_hencetherefore

\NewDocumentCommand{\hence}{}{
  \bool_if:nTF { \g_hencetherefore } {
    \bool_gset_false:N \g_hencetherefore
    hence~
  } {
    \bool_gset_true:N \g_hencetherefore
    therefore~
  }
}

\NewDocumentCommand{\Hence}{}{
  \bool_if:nTF { \g_hencetherefore } {
    \bool_gset_false:N \g_hencetherefore
    Hence,~we~obtain~
  } {
    \bool_gset_true:N \g_hencetherefore
    Therefore,~we~obtain~
  }
}

\ExplSyntaxOff

\NewDocumentEnvironment {athm} {m m} {%
\begin{#1}\label{#2}\global\def\loc{#2}%
}{%
\end{#1}%
}

\NewDocumentEnvironment {adef} {m} {%
\begin{definition}\label{#1}\global\def\loc{#1}%
}{%
\end{definition}%
}

\NewDocumentEnvironment{aproof} {} {%
\begin{proof}[Proof~of~\cref{\loc}]%
}{%
\end{proof}%
}

\ExplSyntaxOff

\newcommand{\idMatrix}{\mathbb{I}}

\newcommand{\andq}{\text{and} \quad}

\newcommand{\andqq}{\text{and} \qquad}
\newcommand{\qandq}{\quad \text{and} \quad}
\newcommand{\qqandqq}{\qquad \text{and} \qquad}

\newcommand{\functionSNN}{\mathcal{N}_{\infty}}
\newcommand{\riskR}{\mathcal{R}_{\infty}}
\newcommand{\riskRR}{\mathcal{R}_{\infty}}

\newcommand{\width}{\mathfrak{h}}
\newcommand{\ones}[1]{\scrv_{#1}}
\newcommand{\with}{\curvearrowleft}

\definecolor{darkgreen}{rgb}{0.0, 0.5, 0.0}
\definecolor{cyan(process)}{rgb}{0.0, 0.72, 0.92}
\definecolor{cinnamon}{rgb}{0.82, 0.41, 0.12}
\definecolor{darkred}{rgb}{0.55, 0.0, 0.0} 
\definecolor{persianblue}{rgb}{0.11, 0.22, 0.73}
\definecolor{resnamecolor}{rgb}{0.0, 0.5, 0.0}
\newcommand{\resname}[1]{\textcolor{resnamecolor}{#1}}
\newcommand{\Exists}{\exists\,}
\newcommand{\Forall}{\forall\,}
\newcommand{\rk}{\operatorname{rank}}
\renewcommand{\d}{\mathrm{d}} 
\renewcommand{\emptyset}{\varnothing}
\newcommand{\Hs}{\operatorname{Hess}}
\newcommand{\spectrum}[1]{\cfadd{def:spectrum}\sigma(#1)}
\DeclarePairedDelimiter{\spro}{\langle}{\rangle}
\newcommand{\submanifold}{\cfadd{def:submanifold}\text{submanifold }}
\newcommand{\PP}{\cfadd{def:Rd_subset_unique_projection}\mathscr{P}}
\newcommand{\pp}{\cfadd{def:Rd_subset_unique_projection}\mathscr{p}}
\newcommand{\bfPP}{\cfadd{def:max_open_set_including_M}\mathbf{P}}
\newcommand{\dimension}{\mathfrak{d}}
\newcommand{\tangentspace}[2]{\cfadd{def:tangent_space}\mathcal{T}_{#2}^{#1}}
\newcommand{\tubular}[2]{\mathbb{B}_{#1}^{#2}}
\newcommand{\gradientG}{\mathcal{G}}
\newcommand{\locmin}{\cfadd{def:loc_min}local minimum }
\newcommand{\locmax}{\cfadd{def:loc_max}local maximum }
\newcommand{\locextrema}{\cfadd{def:loc_extrema}local extremum }
\newcommand{\saddle}{\cfadd{def:saddle}saddle }
\newcommand{\Rect}{\mathcal{A}}

\newcommand{\indicator}[1]{\mathbbm{1}_{\smash{#1}}}

\newcommand{\realapprox}[2]{\mathcal{N} ^{#1}_ {#2}}
\definecolor{codegreen}{rgb}{0,0.6,0}
\definecolor{codegray}{rgb}{0.5,0.5,0.5}
\definecolor{codepurple}{rgb}{0.58,0,0.82}
\definecolor{backcolour}{rgb}{0.95,0.95,0.92}
\lstdefinestyle{mystyle}{
    backgroundcolor=\color{backcolour},   
    commentstyle=\color{codegreen},
    keywordstyle=\color{magenta},
    numberstyle=\tiny\color{codegray},
    stringstyle=\color{codepurple},
    basicstyle=\ttfamily\footnotesize,
    breakatwhitespace=false,         
    breaklines=true,                 
    captionpos=b,                    
    keepspaces=true,                 
    numbers=left,                    
    numbersep=5pt,                  
    showspaces=false,                
    showstringspaces=false,
    showtabs=false,                  
    tabsize=2
}
\makeatletter
\ExplSyntaxOn
\seq_new:N \g_abbrs
\prop_new:N \g_abbr_counts
\tl_new:N \l_abbr_count_tl

\ExplSyntaxOn

\NewDocumentCommand{\abbr}{m m O{#1} m m O{#4} m}{
	\expandafter\newcommand\csname#3\endcsname[1][]{
		\seq_if_in:NnTF \g_abbrs {#1} {
			\prop_get:NnN \g_abbr_counts {#1} \l_abbr_count_tl
			\prop_gput:Nnx \g_abbr_counts {#1} {\int_eval:n {\l_abbr_count_tl + 1}}
			\hyperref[#1]{#7}
		} {
			\seq_gput_left:Nn \g_abbrs {#1}
			\prop_gput:Nnn \g_abbr_counts {#1} {1}
			\expandafter\gdef\csname#1@def\endcsname{#2}
			\phantomsection\label{#1}
			\str_if_eq:nnTF{##1}{}{\emph{#2}}{##1}~(\hyperref[#1]{#7})
		}
	}
	\expandafter\newcommand\csname#6\endcsname[1][]{
		\seq_if_in:NnTF \g_abbrs {#1} {
			\prop_get:NnN \g_abbr_counts {#1} \l_abbr_count_tl
			\prop_gput:Nnx \g_abbr_counts {#1} {\int_eval:n {\l_abbr_count_tl + 1}}
			\hyperref[#1]{#4}
		} {
			\expandafter\gdef\csname#1@def\endcsname{#5}
			\seq_gput_left:Nn \g_abbrs {#1}
			\prop_gput:Nnn \g_abbr_counts {#1} {1}
			\phantomsection\label{#1}
			\str_if_eq:nnTF{##1}{}{\emph{#5}}{##1}~(\hyperref[#1]{#4})
		}
	}
}

\ExplSyntaxOff

\makeatother
\abbr{OP}{optimization problem}{OPs}{optimization problems}{OP}
\abbr{ANN}{artificial neural network}{ANNs}{artificial neural networks}{ANN}
\abbr{DNN}{deep artificial neural network}{DNNs}{deep artificial neural networks}{DNN}
\abbr{GD}{Gradient descent}{GDs}{gradient descent}{GD}
\abbr{GF}{gradient flow}{GFs}{gradient flow}{GF}
\abbr{SGD}{stochastic gradient descent}{SGDs}{stochastic gradient descent}{SGD}
\abbr{DL}{Deep learning}{DLs}{Deep learning}{DL}
\abbr{AI}{artificial intelligence}{AIs}{artificial intelligence}{AI}
\abbr{PDE}{partial differential equation}{PDEs}{partial differential equations}{PDE}
\abbr{ReLUs}{rectified linear unit}{ReLU}{rectified linear unit}
\makeatother

\begin{document}

\title{On the existence of infinitely many realization \\
functions of non-global local minima for \\
artificial neural networks with ReLU activation
}

\author{Shokhrukh Ibragimov$^{1}$, Arnulf Jentzen$^{2,3}$, \\ Timo Kr\"oger$^{4}$, and Adrian Riekert$^{5}$\bigskip\\
\small{$^1$ Applied Mathematics: Institute for Analysis and Numerics,}\vspace{-0.1cm}\\
\small{University of M\"unster, Germany; e-mail: \texttt{sibragim}\textcircled{\texttt{a}}\texttt{uni-muenster.de}}\smallskip\\
\small{$^2$ School of Data Science and School of Artificial Intelligence,} \vspace{-0.1cm}\\
\small{The Chinese University of Hong Kong, Shenzhen (CUHK-Shenzhen),} \vspace{-0.1cm}\\
\small{China; e-mail: \texttt{ajentzen}\textcircled{\texttt{a}}\texttt{cuhk.edu.cn}}\smallskip\\
\small{$^3$ Applied Mathematics: Institute for Analysis and Numerics,}\vspace{-0.1cm}\\
\small{University of M\"unster, Germany; e-mail: \texttt{ajentzen}\textcircled{\texttt{a}}\texttt{uni-muenster.de}}\smallskip\\
\small{$^4$ Applied Mathematics: Institute for Analysis and Numerics,}\vspace{-0.1cm}\\
\small{University of M\"unster, Germany; e-mail: \texttt{timo.kroeger}\textcircled{\texttt{a}}\texttt{uni-muenster.de}}\smallskip\\
\small{$^5$ Applied Mathematics: Institute for Analysis and Numerics,}\vspace{-0.1cm}\\
\small{University of M\"unster, Germany; e-mail: \texttt{ariekert}\textcircled{\texttt{a}}\texttt{uni-muenster.de}}}

\date{\today}

\maketitle

\vspace{-0.5cm}

\begin{abstract}
Gradient descent (GD) type optimization schemes are standard tools for training fully connected feedforward artificial neural networks (ANNs) with rectified linear unit (ReLU) activation and can be considered as temporal discretizations of solutions of gradient flow (GF) differential equations. It is known that, for shallow ReLU ANNs, the risk of every bounded GF trajectory converges under suitable assumptions to the risk of a critical point, a zero of the corresponding gradient function. This motivates studying necessary and sufficient conditions in GF and GD convergence analysis for critical points and their dependence on problem parameters, such as the target function.

In the first main result, we prove for shallow ReLU ANNs that for every dimension $d \in \N$, every $\scra \in \R$, $\scrb \in (\scra,\infty)$, and every positive $\delta \in (0,\infty)$, there exists a Lipschitz continuous target function $\ff \colon [\scra, \allowbreak \scrb]^d \allowbreak \to \R$ such that, for every hidden-layer width $\width \in \N \cap (1, \infty)$, there are uncountably many different realization functions of non-global local minimum points whose risks exceed the infimum of the risk function by more than $\delta$.

In the second main result, we show that, for shallow ReLU ANNs with only one neuron on each layer and a continuous piecewise polynomial target function, there exist at most finitely many different realization functions of critical points.
\end{abstract}

\newpage
\tableofcontents

\section{Introduction}

\GD\ type optimization schemes are the standard instruments to train fully connected feedforward \ANNs\ with \ReLU\ activation. Although there are a huge number of numerical simulations which indicate that \GD\ type optimization schemes effectually train \ANNs\ with \ReLU\ activation, until today there is no mathematical theory which rigorously explains the success of \GD\ type optimization schemes in the training of such \ANNs\ (however, cf., e.g., \cite{Ge2015, LeeJordanRecht2016, chizat2018global, zou2018stochastic, jacot2020neural,du2019gradienta} and the references mentioned therein for several promising mathematical analysis approaches for \GD\ type optimization schemes).

\GD\ type optimization schemes can be considered as temporal discretizations of solutions of \GF\ differential equations and most of the key challenges in the mathematical convergence analysis of \GD\ type optimization schemes seem to already be present in the analysis of \GF\ differential equations. In Eberle et al.~\cite[Theorem~1.2]{eberle2021existence} (cf.\ Bolte \& Pauwels~\cite[Theorem~4]{bolte2020mathematical}, Davis et al.~\cite[Corollary~5.11]{Davis2020}, Jentzen \& Riekert~\cite[Item~$($iv$)$ in Theorem~1.1]{Adrian2021GradientFlows}, and Jentzen \& Riekert~\cite[Theorem~1.3]{jentzen2021existence}) it has been proved, under suitable piecewise polynomial assumptions on the target function and the density of the input data, that the risk of every non-divergent \GF\ trajectory converges in the training of \ANNs\ with one hidden layer and \ReLU\ activation to the risk of a critical point, by which we mean a zero point of the corresponding gradient function, and, taking this into account, it is one of the key research issues in the mathematical convergence analysis of \GF\ trajectories and \GD\ type optimization schemes, respectively, to study sufficient and necessary conditions for critical points of the risk function and, thereby, to obtain an understanding about the appearance of critical points in dependence of the problem parameters such as the target function (cf., e.g., Cheridito et al.~\cite{Florian2021LandscapeAnalysis}).

In the training of \ANNs\ with one hidden layer and \ReLU\ activation the risk function can have saddle points, global minimum points, and non-global local minimum points (cf., e.g., Cheridito et al.~\cite[Theorem~2.4]{Florian2021LandscapeAnalysis}). To establish convergence of the risk of a bounded \GF\ trajectory to the risk of a global minimum point, we thus need to exclude convergence to the risk of any non-global critical point. Besides non-global local minimum points and saddle points, this includes possible zeros of the generalized gradient at which the risk function fails to be differentiable. In the case of saddle points, the articles \cite{Wang2019AvoidSaddlePoints, Ge2015, LeePanageasRecht2019, LeeJordanRecht2016, PanageasPiliouras2017} suggest and study a promising approach which might be successful to verify that the risk of an appropriate \GF\ trajectory does not converge to the risk of a saddle point. From this point of view it seems particularly important to analyze the risk function in terms of its non-global local minimum points in order to better understand the success of \GD\ type optimization schemes in the training of \ANNs.

The scientific literature has dealt with these non-global local minimum points in a variety of ways. There are several examples for finitely many training data and architectures of \ANNs\ for which the considered risk function has non-global local minimum points (see, e.g., \'{S}wirszcz et al.~\cite{swirszcz2017local}). Moreover, non-global local minimum points could be found in the risk landscape of \ANNs\ with one hidden layer and \ReLU\ activation in special student-teacher setups with the probability distribution of the input data given by the normal distribution (see Safran \& Shamir~\cite{safran2018spurious}). In other cases, where the target function has a very simple form, the critical points of the risk landscape are fully characterized and thus all local minimum points are known (see Cheridito et al.~\cite[Corollary~2.15]{CHERIDITO2022101646}, Cheridito et al.~\cite{Florian2021LandscapeAnalysis}, and Jentzen \& Riekert~\cite[Corollary~2.11]{Adrian2021StochasticGDconvergence}). Additionally, in the case of \ANNs\ with linear activation and finitely many training data it was shown that all local minimum points of the risk function corresponding to the squared error loss are global minimum points (cf.\ Kawaguchi~\cite{kawaguchi2016deep} and Laurent \& von Brecht~\cite{pmlr-v80-laurent18a}). For a connection between the critical points of the risk function and the critical points of the risk function with regard to a larger network width we refer to the articles by Zhang et al.~\cite{zhang2021embeddingb,zhang2021embeddinga}.

Further progress in this regard has been made in the so-called overparameterized regime. In this regime it was demonstrated for different situations that the set of all \ANNs\ with risk equal to 0 forms a high-dimensional submanifold of the parameter space (see Cooper~\cite{cooper2018loss}). Analogous results were also shown in the non-overparameterized regime (cf.\ Dereich \& Kassing~\cite{dereich2021minimal}, Fehrman et al.~\cite{FehrmanGessJentzen2020}, and Jentzen \& Riekert~\cite{Adrian2021ConvergencePiecewise}). In addition, for \ANNs\ with one hidden layer and quadratic activation with finitely many training data, using different assumptions, it was shown that all local minimum points in the risk landscape corresponding to the squared error loss are global minimum points (see Du \& Lee~\cite{du2018power} and Soltanolkotabi et al.~\cite{soltanolkotabi2018theoretical}). In the case of mild overparameterization in special student-teacher setups for \ANNs\ with one hidden layer and quadratic activation with random training data there are mathematical analyses for the probability of the occurrence of non-global local minimum points (see Mannelli et al.~\cite{mannelli2020optimization}). In particular, the influence of the number of the training data, the input dimension, and the number of the hidden neurons of the teacher \ANN\ was examined more closely.

We also want to mention approaches to visualize the structure of the risk landscape, which is particularly interesting in a local environment of critical points. The article by Li et al.~\cite{li2018visualizing} presents different ways to get visual access to the high-dimensional risk landscape, discusses disadvantages of these approaches, and suggests an alternative with the so-called filter normalization. Without this filter normalization, numerical experiments suggest that arbitrary two-dimensional patterns can be found in the risk landscape of wide and deep \ANNs\ using common training data sets such as FashionMNIST and CIFAR10 (see Skorokhodov \& Burtsev~\cite{skorokhodov2019loss}). There are different attempts to mathematically explain this phenomenon and to prove that such patterns can be found around approximate global minimum points in special situations, for example, using the universal approximation theorem (see Czarnecki et al.~\cite{czarnecki2020deep}).

\begin{figure}[!htb]
\centering
\subfloat{\includegraphics[width=8.2cm]{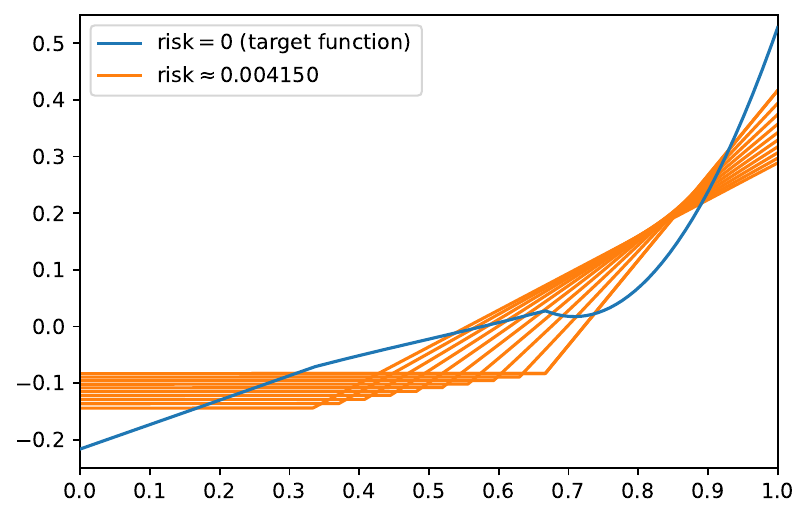}} \,
\subfloat{\includegraphics[width=8.2cm]{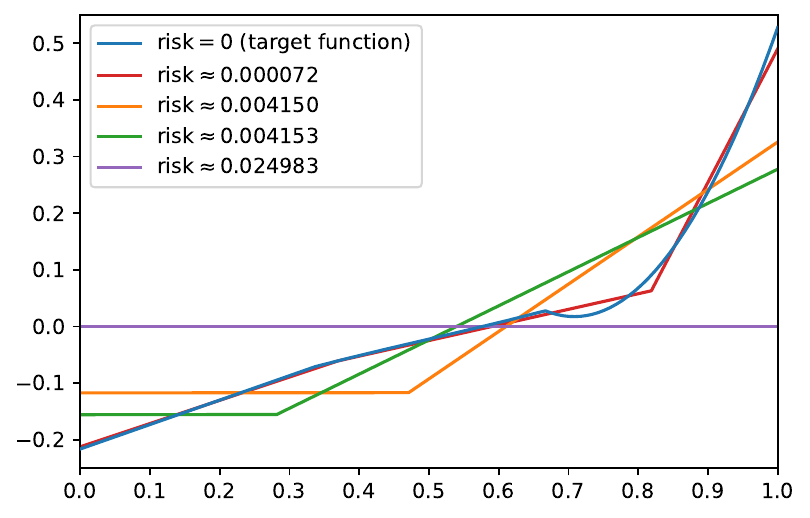}}
\caption{Numerical simulations associated to \cref{cor:uncountably_many_realizations_of_non_global_local_minimas_with_fixed_delta} in the case where $d = 1$, $\width = 4$, $\dimension = 13$, $\scra = 0$, $\scrb = 1$, $\alpha = \nicefrac{1}{3}$, and $\beta = \nicefrac{2}{3}$ in \cref{cor:uncountably_many_realizations_of_non_global_local_minimas_with_fixed_delta}: On the left picture we approximately plot the target function $\ff \colon [0, 1] \to \R$ (cf.\ \cref{0_eq:f}) and 10 different realization functions of non-global local minimum points (cf.\ \cref{0_eq:M}) of $ \riskR \colon \R^{13} \to \R $. On the right picture we plot the result from a simulation with the following setting. We randomly initialize 50 \ANNs\ with the Xavier initialization (we initialize the weights to be normally distributed with mean $0$ and variance $\nicefrac{2}{5}$ and we initialize the biases to $0$), then we approximately train these \ANNs\ with the \GD\ optimization method using a learning rate of $\nicefrac{1}{20}$ until the maximum norm of the generalized gradient function (cf.\ \cref{eq:loss:gradient}) evaluated at the current position of the \GD\ process is strictly less than $10^{-4}$, and, thereafter, we gradually plot the realization functions of the resulting \ANNs\ whereby a realization function is not drawn if a realization function with an $L^2$-distance strictly less than $10^{-2}$ has already been drawn. This simulation resulted in four different realization functions (red, orange, green, and purple) on the right picture whereby one realization function (orange) can also approximately be found in the left picture. We also refer to \cref{list:python} for the {\sc Python} source code used to create \cref{fig:local:minimum:points}.}
\label{fig:local:minimum:points}
\end{figure}

In view of these scientific findings, we are in this article particularly interested in the study of non-global local minimum points of the risk functions. In the main results of this work we establish two basic results regarding the appearance of critical points in the training of \ANNs\ with one hidden layer and \ReLU\ activation. Specifically, in the first main result of this work, see \cref{theorem:negative} below, we prove in the training of shallow \ReLU\ \ANNs\ that for every arbitrary dimension $d \in \N$, every $\scra \in \R$, $\scrb \in (\scra,\infty)$, and every arbitrarily large positive $\delta \in (0,\infty)$ we have that there exists a Lipschitz continuous target function $\ff \colon [\scra, \scrb]^d \to \R$ such that for every number $\width \in \N \cap (1, \infty)$ of neurons on the hidden layer we have that the risk function has uncountably many different realization functions of non-global local minimum points whose risks are strictly larger than the sum of the infimum of the risk function and the arbitrarily large positive real number $\delta$ (see also \cref{fig:local:minimum:points} in \cref{sec:negative_result} below for a graphical illustration related to the statement of \cref{theorem:negative}). \cref{theorem:negative} thus suggests even in the situation where the target function is Lipschitz continuous that the training problem might be very challenging due to the appearance of infinitely many different realization functions of non-global local minimum points. To the best of our knowledge, \cref{theorem:negative} is the first result in the scientific literature which rigorously proves in the training of fully connected \ANNs\ with \ReLU\ activation that there exists a target function such that the risk function has infinitely many different realization functions of non-global local minimum points. We now present the precise statement of \cref{theorem:negative}.
\begin{samepage}
\begin{theorem}\label{theorem:negative}
Let $d \in \N$, $\delta, \scra \in \R$, $\scrb \in (\scra, \infty)$, for every $\width \in \N$ let $\dimension_{\width} \in \N$ satisfy $\dimension_{\width} = \allowbreak d \width + \allowbreak 2 \width \allowbreak + 1$, for every $\width \in \N$, $\theta = \allowbreak (\theta_1, \allowbreak \dots, \allowbreak \theta_{\dimension_{\width}}) \in \R^{\dimension_{\width}}$ let $\cN_{\width}^{\theta} \in C([\scra, \scrb]^d, \R)$ satisfy for all $x = \allowbreak (x_1, \allowbreak \dots, \allowbreak x_d) \allowbreak \in \allowbreak [\scra, \allowbreak \scrb]^d$ that
\begin{equation}
\textstyle \cN_{\width}^{\theta}(x) = \theta_{\dimension_{\width}} + \sum_{i = 1}^{\width} \theta_{d \width + \width + i} \max\bigl\{\theta_{d \width + i} + \sum_{j = 1}^d \theta_{(i - 1)d + j} x_j, 0\bigr\},
\end{equation}
and for every $\width \in \N$, $f \in C([\scra, \scrb]^d, \R)$ let $\cR_{f, \width} \colon \R^{\dimension_{\width}} \to \R$ satisfy for all $\theta \in \R^{\dimension_{\width}}$ that $\cR_{f, \width}(\theta) = \allowbreak \int_{[\scra, \scrb]^d} \abs{\cN_{\width}^{\theta}(y) - f(y)}^2 \, \d y$. Then there exists a Lipschitz continuous $\ff \colon \allowbreak [\scra, \allowbreak \scrb]^d \allowbreak \to \R$ such that for all $\width \in \N \cap (1, \infty)$ it holds that
\begin{multline}
\textstyle \bigl\{v \in C([\scra, \scrb]^d, \R) \colon \bigl[\Exists \theta \in \{\vartheta \in \R^{\dimension_{\width}} \colon v = \cN_{\width}^{\vartheta}\}, \eps \in (0, \infty) \colon \\
\textstyle \cR_{\ff, \width}(\theta) = \inf_{\vartheta \in [-\eps, \eps]^{\dimension_{\width}}} \cR_{\ff, \width}(\theta + \vartheta) > \delta + \inf_{\vartheta \in \R^{\dimension_{\width}}} \cR_{\ff, \width}(\vartheta)\bigr]\bigr\}
\end{multline}
is an uncountable set.
\end{theorem}
\end{samepage}

\cref{theorem:negative} is an immediate consequence of \cref{cor:uncountably_many_realizations_of_non_global_local_minimas_with_varying_delta} in \cref{subsec:existence_of_inf_many_realizations}. In the second main result of this work, see \cref{theorem:positive} below, we provide in a special situation sufficient conditions to ensure that there are at most finitely many different realization functions of non-global local minimum points. Specifically, in \cref{item1:theorem:positive} in \cref{theorem:positive} below we show in the training of shallow \ReLU\ \ANNs\ in the special situation where there is only one neuron on each layer (corresponding to the case $d = 1 = \width$ in \cref{theorem:negative} above) and where the target function is continuous and piecewise polynomial that there exist at most finitely many different realization functions of critical points. This enables us to conclude in \cref{item2:theorem:positive} in \cref{theorem:positive} that (in contrast to the situation of \cref{theorem:negative} above) there exist at most finitely many different realization functions of (non-global) local minimum points. In addition, \cref{item1:theorem:positive} in \cref{theorem:positive} together with \cite[Item~$($v$)$ in Theorem~1.1]{Adrian2021GradientFlows} and \cite[Theorem~1.2]{eberle2021existence} allows us to conclude in \cref{item3:theorem:positive} in \cref{theorem:positive} that in the training of such \ANNs\ we have that the risk of every non-divergent \GF\ trajectory converges to the risk of a global minimum point provided that the initial risk is sufficiently small. To describe a \GF\ trajectory, we need to specify an appropriate generalized gradient function in \cref{theorem:positive} as the risk function is not differentiable in the case of \ANNs\ with \ReLU\ activation (due to the fact that the \ReLU\ activation function $\R \ni x \mapsto \max\{x, 0\} \in \R$ fails to be differentiable at the origin). As in \cite{Adrian2021StochasticGDconvergence} (cf., e.g., also Cheridito et al.~\cite{CHERIDITO2022101646}) we accomplish this by means of an approximation procedure in which the \ReLU\ activation function $\R \ni x \mapsto \max\{x, 0\} \in \R$ is approximated through appropriate continuously differentiable functions whose derivatives converge pointwise to the left-derivative of the \ReLU\ activation function; see \cref{eqn:theorem:positive:rect} in \cref{theorem:positive}. We now present the precise statement of \cref{theorem:positive}.

\begin{samepage}
\begin{theorem}\label{theorem:positive}
Let $n \in \N$, $\fx_0, \fx_1, \dots, \fx_n, \scra, \scrb \in \R$ satisfy $\scra \allowbreak = \fx_0 \allowbreak < \fx_1 \allowbreak < \dots \allowbreak < \fx_n \allowbreak = \scrb$, let $f \allowbreak \in C([\scra, \allowbreak \scrb], \allowbreak \R)$ satisfy for all $j \in \{1, 2, \dots, n\}$ that $f|_{[\fx_{j - 1}, \fx_j]}$ is a polynomial, let $\Rect_r \allowbreak \colon \allowbreak \R \allowbreak \to \allowbreak \R$, $r \in \N \cup \{\infty\}$, satisfy for all $x \in \R$ that $(\cup_{r \in \N} \{\Rect_r\}) \subseteq C^1(\R, \R)$, $\Rect_{\infty}(x) = \max\{x, 0\}$, $\sup_{r \in \N} \sup_{y \in [-\abs{x}, \abs{x}]} \abs{(\Rect_r)'(y)} < \infty$, and
\begin{equation}\label{eqn:theorem:positive:rect}
\textstyle \limsup_{r \to \infty} (\abs{\Rect_r(x) - \Rect_{\infty}(x)} + \abs{(\Rect_r)'(x) - \mathbbm{1}_{(0, \infty)}(x)}) = 0,
\end{equation}
for every $r \in \N \cup \{\infty\}$ let $\cN_r = (\cN_r^{\theta})_{\theta \in \R^4} \colon \R^4 \to C([\scra, \scrb], \R)$ and $\cR_r \colon \R^{4} \to \R$ satisfy for all $\theta \allowbreak = (\theta_1, \allowbreak \dots, \allowbreak \theta_4) \allowbreak \in \R^{4}$, $x \allowbreak \in [\scra, \allowbreak \scrb]$ that $\cN_r^{\theta}(x) \allowbreak = \theta_4 \allowbreak + \theta_3 [\Rect_r(\theta_2 \allowbreak + \theta_1 x)]$ and $\cR_r(\theta) \allowbreak = \int_{\scra}^{\scrb} \abs{\cN_r^{\theta}(y) - f(y)}^2 \, \d y$, and let $\gradientG \colon \R^{4} \to \R^{4}$ satisfy for all $\theta \in \{\vartheta \allowbreak \in \R^{4} \colon \allowbreak ((\nabla \cR_r)(\vartheta))_{r\in \N} \text{ is convergent}\}$ that $\gradientG(\theta) \allowbreak = \lim_{r \to \infty} (\nabla \cR_r)(\theta)$. Then
\begin{enumerate}[label = (\roman*)]
\item
\label{item1:theorem:positive} it holds that $\{v \in C([\scra, \scrb], \R) \colon (\Exists \theta \in \gradientG^{-1}(\{0\}) \colon v = \cN_{\infty}^{\theta})\}$ is a finite set,
  
\item
\label{item2:theorem:positive} it holds that
\begin{multline}
\textstyle \bigl\{ v \in C([\scra, \scrb], \R) \colon \bigl[ \Exists \theta \in \{ \vartheta \in \R^4 \colon v = \cN_{\infty}^{\vartheta} \}, \varepsilon \in (0,\infty) \colon \\
\textstyle \riskRR (\theta) = \inf\nolimits_{ \vartheta \in [-\varepsilon,\varepsilon]^4 } \riskRR (\theta + \vartheta) \bigr] \bigr\}
\end{multline}
is a finite set, and

\item
\label{item3:theorem:positive} there exists $\eps \in (0, \infty)$ which satisfies\footnote{\Nobs that for every $\dimension \in \N$, $x = (x_1, \dots, x_{\dimension})$, $y = (y_1, \dots, y_{\dimension}) \in \R^{\dimension}$ it holds that $\spro{x, y} = \sum_{j = 1}^{\dimension} x_j y_j$ (standard scalar product) and $\norm{x} = (\spro{x, x})^{\nicefrac{1}{2}}$ (standard norm).} for all $\Theta = (\Theta_t)_{t \in [0, \infty)} \allowbreak \colon [0, \allowbreak \infty) \to \R^4$ with $\riskRR(\Theta_0) \allowbreak \le \allowbreak \eps \allowbreak + \inf_{\vartheta \in \R^4} \riskRR(\vartheta)$, $\Forall t \allowbreak \in [0, \allowbreak \infty) \allowbreak \colon \allowbreak \Theta_t \allowbreak = \Theta_0 - \int_0^t \gradientG(\Theta_s) \, \d s$, and $\liminf_{t \to \infty} \norm{\Theta_{t}} < \infty $ that
\begin{equation}
\textstyle \limsup_{t \to \infty} \riskRR(\Theta_t) = \inf_{\vartheta \in \R^4} \riskRR(\vartheta).
\end{equation}
\end{enumerate}
\end{theorem}
\end{samepage}

\cref{theorem:positive} is an immediate consequence of \cref{cor:convergence_to_global_minima} in \cref{subsec:main_positive_results}. The remainder of this article is organized in the following way.

In \cref{sec:values_of_gen_grad_at_loc_min} we prove in \cref{prop:loss:diff:vc} and \cref{prop:loss:diff:wb} a few basic differentiability properties for the risk function and we establish in \cref{prop:local:minima:gradient} that every local minimum point of the risk function is a critical point (a zero of the generalized gradient function). In \cref{sec:differential_geometric_prelimiaries} we recall some basic concepts and elementary results from differential geometry and we collect in \cref{prop:class:crit:points} some necessary and sufficient conditions for local extremum and saddle points. In \cref{sec:negative_result} we employ \cref{prop:class:crit:points} from \cref{sec:differential_geometric_prelimiaries} to establish in \cref{cor:uncountably_many_realizations_of_non_global_local_minimas_with_varying_delta} that there exist uncountably many Lipschitz continuous target functions such that the associated risk functions have infinitely many realization functions of non-global local minimum points. \cref{theorem:negative} above is a direct consequence of \cref{cor:uncountably_many_realizations_of_non_global_local_minimas_with_varying_delta}. We also refer to \cref{fig:local:minimum:points} in \cref{sec:negative_result} for a graphical illustration related to the statement of \cref{cor:uncountably_many_realizations_of_non_global_local_minimas_with_varying_delta}. Finally, in \cref{sec:on_finitely_many_realization_functions} we prove in \cref{cor:finite_number_of_realizations} in the special situation where the target function is continuous and piecewise polynomial and where the considered shallow \ReLU\ \ANNs\ have one neuron on each layer that there exist only finitely many different realization functions of all critical points of the risk function (of all zeros of the generalized gradient function). \cref{theorem:positive} above can then be shown by combining \cref{prop:local:minima:gradient}, \cref{cor:finite_number_of_realizations}, \cite[Item~$($v$)$ in Theorem~1.1]{Adrian2021GradientFlows}, and \cite[Theorem~1.2]{eberle2021existence}. This is precisely the subject of \cref{cor:convergence_to_global_minima} in \cref{sec:on_finitely_many_realization_functions}.

\section{Values of the generalized gradient function at local minimum points}
\label{sec:values_of_gen_grad_at_loc_min}

In this section we establish in \cref{prop:local:minima:gradient} in \cref{subsec:values_of_general_grad_at_loc_min} below that every local minimum point $\theta \in \R^{\dimension} = \R^{d \width + 2 \width + 1}$ of the risk function $\riskRR \colon \R^{\dimension} \to \R$ is a critical point in the sense that it is a zero of the generalized gradient function $\gradientG \colon \R^{\dimension} \to \R^{\dimension}$. Our proof of \cref{prop:local:minima:gradient} uses the essentially well-known representation result for the generalized gradient function $ \gradientG \colon \R^{ \dimension } \to \R^{ \dimension } $ in \cref{prop:limit:lr} and the elementary relationships between the generalized gradient function $ \gradientG \colon \R^{ \dimension } \to \R^{ \dimension } $ and the first-order partial derivatives of the risk function $ \riskRR \colon \R^{ \dimension } \to \R $ in \cref{prop:loss:diff:vc} and \cref{prop:loss:diff:wb} in \cref{subsec:differentiability_prop_for_general_grad} below. The proof of \cref{prop:limit:lr} can be derived analogously to the proof of \cite[Proposition~2.2]{Adrian2021GradientFlows}. \cref{prop:loss:diff:vc} and \cref{prop:loss:diff:wb} are slight generalizations of \cite[Lemma~2.6]{CHERIDITO2022101646} and \cite[Lemma~2.7]{CHERIDITO2022101646}, respectively. Our proofs of \cref{prop:loss:diff:vc} and \cref{prop:loss:diff:wb} make use of \cref{cor:interchange} in \cref{subsec:diff_properties_for_parameter_integrals} below, which is a direct corollary of the elementary differentiability result in \cref{lem:interchange} in \cref{subsec:diff_properties_for_parameter_integrals}. \cref{lem:interchange} and \cref{cor:interchange} are slight generalizations of \cite[Lemma~4]{Adrian2021ConvergencePiecewise}. Only for completeness we also include in this section detailed proofs for \cref{lem:interchange} and \cref{cor:interchange}.

In \cref{setting:snn0} in \cref{subsec:ANNs_multidimensional_input_hidden_layer} below we describe our mathematical setup to introduce the target function $f \colon [\scra, \scrb]^d \to \R$, the unnormalized probability distribution of the input data $ \mu \colon \cB([\scra, \allowbreak \scrb]^d) \allowbreak \to [0, \allowbreak \infty]$, the realization functions $\realapprox{\theta}{\infty} \colon \R^d \to \R$, $\theta \in \R^{\dimension}$, the risk function $\riskRR \colon \R^{\dimension} \to \R$, and the generalized gradient function $\gradientG \colon \R^{\dimension} \to \R^{\dimension}$. For the convenience of the reader we recall the notions of a local minimum point, of a local maximum point, of a local extremum point, and of a saddle point in \cref{def:loc_min,def:loc_max,def:loc_extrema,def:saddle} in \cref{subsec:local_extrema_and_saddle_points} below.


\subsection{Artificial neural networks (ANNs) with multidimensional input and hidden layer}
\label{subsec:ANNs_multidimensional_input_hidden_layer}

\cfclear
\begin{setting}\label{setting:snn0}
Let $d, \width, \dimension \in \N$, $\scra \in \R$, $\scrb \in (\scra, \infty)$ satisfy $\dimension = d \width + 2\width + 1$, let $e_1, \allowbreak e_2, \allowbreak \dots, \allowbreak e_d \allowbreak \in \R^d$ satisfy $e_1 \allowbreak = (1, \allowbreak 0, \allowbreak 0, \allowbreak \dots, \allowbreak 0)$, $e_2 \allowbreak = (0, \allowbreak 1, \allowbreak 0, \allowbreak \dots, \allowbreak 0)$, $\dots$, $e_d \allowbreak = (0, \allowbreak 0, \allowbreak \dots, \allowbreak 0, \allowbreak 1) \allowbreak \in \R^d$, let $f \colon [\scra, \allowbreak \scrb]^d \allowbreak \to \R$ be measurable, let $\Rect_r \colon \R \to \R$, $r \in \N \cup \{\infty\}$, satisfy for all $x \in \R$ that $(\cup_{r \in \N} \{\Rect_r \}) \subseteq C^1 (\R, \R)$, $\Rect_\infty(x) = \max \{x, 0\}$, $\sup_{r \in \N} \sup_{y \in [- \abs{x}, \abs{x}]} \abs{(\Rect_r)'(y)} < \infty$, and
\begin{equation}\label{setting:snn0:assumption:rect}
\limsup\nolimits_{r \to \infty} (\abs{\Rect_r(x) - \Rect_\infty(x)} + \abs{(\Rect_r)'(x) - \indicator{(0, \infty)}(x)}) = 0,
\end{equation}
let $\mu \colon \cB ([\scra, \scrb]^d) \to [0, \infty]$ be a finite measure, assume $\int_{[\scra, \scrb]^d}\abs{f(x)}^2 \,\mu(\d x) \allowbreak < \infty$, let $\lambda \colon \cB([\scra, \allowbreak \scrb]^d) \allowbreak \to [0, \allowbreak \infty]$ be the Lebesgue-Borel measure on $[\scra, \allowbreak \scrb]^d$, for every $\theta = \allowbreak (\theta_1, \allowbreak \dots, \allowbreak \theta_{\dimension}) \allowbreak \in \R^{\dimension}$, $i \in \allowbreak \{1, \allowbreak 2, \allowbreak \dots, \allowbreak \width\}$ let $\fw_i^{\theta} \allowbreak = (\fw_{i, 1}^{\theta}, \allowbreak \dots, \allowbreak \fw_{i, d}^{\theta}) \allowbreak \in \R^d$, $\fb_i^{\theta}, \fv_i^{\theta}, \fc^{\theta} \in \R$ satisfy $\wedge_{j = 1}^d [\fw_{i, j}^{\theta} \allowbreak = \theta_{(i - 1)d + j}]$, $\fb_i^{\theta} \allowbreak = \theta_{d\width + i}$, $\fv_i^{\theta} \allowbreak = \theta_{d\width + \width + i}$, and $\fc^{\theta} \allowbreak = \theta_{\dimension}$ and let $I_i^\theta \subseteq \R^d$ satisfy $I_i^\theta = \allowbreak \{x \allowbreak \in [\scra, \scrb]^d \colon \allowbreak \fb_i^{\theta} + \allowbreak \spro{\fw_i^{\theta}, x} \allowbreak > 0\}$, for every $r \in \allowbreak \N \allowbreak \cup \{\infty\}$ let $\realapprox{}{r} \allowbreak = (\realapprox{\theta}{r})_{\theta \in \R^{\dimension}} \colon \R^{\dimension} \to C(\R^d , \R) $ and $\cR_r \colon \R^{\dimension} \to \R$ satisfy for all $\theta \allowbreak \in \R^{\dimension}$, $x \allowbreak \in \R^d$ that 
\begin{equation}\label{setting:snn0:eq:realization_risk}
\textstyle \realapprox{\theta}{r}(x) = \fc^{\theta} + \smallsum_{i = 1}^\width \fv_i^{\theta} \bigl[\Rect_r (\fb_i^{\theta} + \spro{\fw_i^{\theta}, x})\bigr] \qandq \cR_r(\theta) = \int_{[\scra, \scrb]^d} \abs{\realapprox{\theta}{r}(y) - f(y)}^2 \, \mu (\d y),
\end{equation}
and let $\gradientG = (\gradientG_1, \dots, \gradientG_{\dimension}) \colon \R^{\dimension} \to \R^{\dimension}$ satisfy for all $\theta \in  \{\varphi \in \R^{\dimension} \colon ((\nabla \cR_r )(\varphi))_{r \in \N} \allowbreak \ \text{is convergent}\}$ that $\gradientG(\theta) = \lim_{r \to \infty} (\nabla \cR_r)(\theta)$.
\end{setting}

\subsection{Differentiability properties for parameter-dependent Lebesgue integrals}
\label{subsec:diff_properties_for_parameter_integrals}


\cfclear
\begin{lemma}\label{lem:interchange}
Let $d \in \N$, $\scra \in \R$, $\scrb \in (\scra, \infty)$, let $\phi \colon \R \times [\scra, \scrb]^d \to \R$ be measurable, let $\mu \colon \allowbreak \cB ([\scra, \allowbreak \scrb]^d) \allowbreak \to [0, \allowbreak \infty]$ be a measure, assume for all $x \in \R$ that $\int_{[\scra, \scrb]^d} \abs{\phi(x,s)} \, \mu(\d s) \allowbreak < \infty$, let $ \Phi \colon \R \to \R $ satisfy for all $x \in \R$ that
\begin{equation}\label{eq:lem:interchange:Phi}
\textstyle \Phi(x) = \int_{[\scra, \scrb]^d} \phi(x, s) \, \mu (\d s),
\end{equation}
let $E \in \cB ([\scra, \scrb]^d)$ satisfy $\mu ([\scra, \scrb]^d \backslash E) = 0$, let $c \colon E \to \R$ be measurable, let $x \in \R$, $\delta \in (0, \infty)$ satisfy for all $h \in (-\delta, \delta)$, $s \in E$ that $\abs{\phi(x + h, s) - \phi(x, s)} \le \abs{h} \abs{c(s)}$, assume $\int_E \abs{c(s)} \, \mu (\d s) \allowbreak < \infty$, and assume for all $s \in E$ that $\R \ni v \mapsto \phi(v, s) \in \R$ is differentiable at $x$. Then
\begin{enumerate}  [label = (\roman*)]
\item it holds that $\Phi$ is differentiable at $x$ and

\item it holds that
\begin{equation}
\textstyle \Phi'(x) = \int_E \rbr[\bigr]{\tfrac{\partial}{\partial x} \phi }(x, s) \, \mu( \d s).
\end{equation}
\end{enumerate}
\end{lemma}
\begin{cproof}{lem:interchange}
\Nobs that \cref{eq:lem:interchange:Phi} and the assumption that $\mu ([\scra, \scrb]^d \backslash E) = 0$ demonstrate for all $h \in \R \backslash \{0\}$ that
\begin{equation}\label{lem:interchange:eq1}
\begin{split}
\textstyle h^{-1} [\Phi(x + h) - \Phi(x)] & = \textstyle \int_{[\scra, \scrb]^d} h^{-1} [\phi(x + h, s) - \phi(x, s)] \, \mu(\d s) \\
& \textstyle = \int_E h^{-1}[\phi(x + h, s) - \phi(x, s)] \, \mu(\d s).
\end{split}
\end{equation}
In the next step we \nobs that the assumption that for all $s \in E$ it holds that $\R \ni v \mapsto \phi(v, s) \in \R$ is differentiable at $x$ shows that for all $s \in E$ it holds that
\begin{equation}\label{lem:interchange:eq2}
\limsup\nolimits_{\R\backslash\{0\} \ni h \to 0} \babs{h^{-1} [\phi(x + h, s) - \phi(x, s)] - \bigl(\tfrac{\partial}{\partial x} \phi\bigr)(x, s)} = 0.
\end{equation}
Furthermore, we \nobs that the assumption that for all $h \in (-\delta, \delta)$, $s \in E$ it holds that $\abs{\phi(x + h, s) - \phi(x, s)} \le \abs{h} \abs{c(s)}$ proves that for all $h \in (-\delta, \delta) \backslash \{0\}$, $s \in E$ it holds that
\begin{equation}
\abs{ h^{-1} [\phi(x + h, s) - \phi(x, s) ]} \le \abs{ c(s) }.
\end{equation}
Combining \cref{lem:interchange:eq1}, \cref{lem:interchange:eq2}, the assumption that $\int_E \abs{c(s)} \, \mu(\d s) < \infty$, and Lebesgue's dominated convergence theorem \hence assures that
\begin{equation}
\begin{split}
& \lim\nolimits_{\R\backslash\{0\} \ni h \to 0} \bigl( h^{-1} [ \Phi (x + h) - \Phi (x) ] \bigr) \\
& = \textstyle \int_E \bigl[ \lim\nolimits_{\R\backslash\{0\} \ni h \to 0} \bigl( h^{-1} [ \phi (x + h, s) - \phi (x, s) ] \bigr) \bigr] \, \mu(\d s) = \textstyle \int_E \rbr[\bigr]{\tfrac{\partial}{\partial x} \phi}(x, s )\, \mu( \d s ).
\end{split}
\end{equation}
\end{cproof}


\cfclear
\begin{corollary}\label{cor:interchange}
Let $d, n \in \N$, $\scra \in \R$, $\scrb \in (\scra, \infty)$, let $\phi \colon \R^n \times [\scra, \scrb]^d \to \R$ be measurable, let $\mu \colon \cB([\scra, \scrb]^d) \to [0, \infty] $ be a measure, assume for all $x \in \R^n$ that $\int_{[\scra, \scrb]^d} \abs{\phi(x,s)} \, \mu(\d s) \allowbreak < \infty$, let $\Phi \colon \R^n \to \R$ satisfy for all $x \in \R^n$ that
\begin{equation}
\textstyle \Phi(x) = \int_{[\scra, \scrb]^d} \phi(x, s) \, \mu (\d s),
\end{equation}
let $E \in \cB([\scra, \scrb]^d)$ satisfy $\mu([\scra, \scrb]^d \backslash E) = 0$, let $c \colon E \to \R$ be measurable, let $x_1, \allowbreak x_2, \allowbreak \dots, \allowbreak x_n \allowbreak \in \R$, $j \allowbreak \in \allowbreak \{1, \allowbreak 2, \dots, \allowbreak n\}$, $\delta \in (0, \infty)$ satisfy for all $h \in (- \delta, \delta)$, $s \in E$ that
\begin{equation}
\abs{\phi(x_1, \dots, x_{j - 1}, x_j + h, x_{j + 1}, \dots, x_n, s) - \phi(x_1, \dots, x_n, s)} \le \abs{h} \abs{c(s)},
\end{equation}
assume $\int_E \abs{c(s)} \, \mu(\d s) < \infty$, and assume for all $ s \in E $ that $\R \ni v \mapsto \allowbreak \phi(x_1, \allowbreak \dots, \allowbreak x_{j-1}, \allowbreak v, \allowbreak x_{j+1}, \allowbreak \dots, \allowbreak x_n, \allowbreak s) \allowbreak \in \R$ is differentiable at $x_j$. Then
\begin{enumerate} [label = (\roman*)]
\item
\label{cor:interchange:item1} it holds that $\R \ni v \mapsto \Phi(x_1, \dots, x_{j-1}, v , x_{j+1}, \dots, x_n) \in \R$ is differentiable at $x_j$ and

\item
\label{cor:interchange:item2} it holds that
\begin{equation}
\textstyle \rbr[\bigr]{\tfrac{\partial}{\partial x_j} \Phi}(x_1, \dots, x_n) = \int_E \rbr[\bigr]{\tfrac{\partial}{\partial x_j} \phi}(x_1, \dots, x_n, s) \, \mu(\d s).
\end{equation}
\end{enumerate}
\end{corollary}
\begin{cproof}{cor:interchange}
\Nobs that \cref{lem:interchange} establishes \cref{cor:interchange:item1,cor:interchange:item2}.
\end{cproof}

\subsection{Differentiability properties for the risk function}
\label{subsec:differentiability_prop_for_general_grad}


\cfclear
\begin{proposition}[\resname{Explicit representations for the generalized gradient}]\label{prop:limit:lr}
Assume \cref{setting:snn0} and let $\theta \in \R^{\dimension}$. Then it holds for all $i \in \{1, \allowbreak 2, \allowbreak \dots, \allowbreak \width\}$, $j \in \{1, \allowbreak 2, \allowbreak \dots, \allowbreak d\} $ that
\begin{equation}\label{eq:loss:gradient}
\begin{split}
\textstyle \gradientG_{(i - 1) d + j}(\theta) & = \textstyle 2 \fv_i^{\theta} \int_{I_i^\theta} \scalprod{e_j}{x} (\realapprox{\theta}{\infty}(x) - f(x)) \, \mu (\d x), \\
\gradientG_{d \width + i}(\theta) & = \textstyle 2 \fv_i^{\theta} \int_{I_i^\theta} (\realapprox{\theta}{\infty}(x) - f(x)) \, \mu (\d x), \\
\gradientG_{d \width + \width + i}(\theta) & = \textstyle 2 \int_{[\scra, \scrb]^d} \bigl[\Rect_\infty(\fb_i^{\theta} + \spro{\fw_i^{\theta}, x})\bigr] (\realapprox{\theta}{\infty}(x) - f(x)) \, \mu (\d x), \\
\text{and} \qquad \gradientG_{\dimension}(\theta) & = \textstyle 2 \int_{[\scra, \scrb]^d} (\realapprox{\theta}{\infty}(x) - f(x)) \, \mu (\d x)\ifnocf.
\end{split}
\end{equation}
\cfout[.]
\end{proposition}
\begin{cproof}{prop:limit:lr}
\Nobs that \cref{setting:snn0:assumption:rect} and \cref{setting:snn0:eq:realization_risk} establish \cref{eq:loss:gradient} (cf., e.g., \cite[Proposition~2.2]{Adrian2021GradientFlows} and \cite[Items~$($v$)$ and $($vi$)$ in Theorem~2.9]{hutzenthaler2021convergence}).
\end{cproof}


\cfclear
\begin{lemma}[\resname{Differentiability w.r.t.\ outer weights and an outer bias}]\label{prop:loss:diff:vc}
Assume \cref{setting:snn0} and let $\theta = (\theta_1, \dots, \theta_{\dimension}) \in \R^{\dimension}$. Then
\begin{enumerate} [label = (\roman*)]
\item\label{prop:loss:diff:vc:item1} it holds for all $i \in \N \cap (d \width + \width, \dimension]$ that $\R \ni v \mapsto \riskRR(\theta_1, \dots, \theta_{i-1}, v, \theta_{i+1}, \dots, \theta_{\dimension} ) \in \R$ is differentiable at $\theta_i$ and

\item\label{prop:loss:diff:vc:item2} it holds for all $i \in \N \cap (d \width + \width, \dimension]$ that $(\frac{\partial}{\partial \theta_i} \riskRR)(\theta) = \gradientG_i(\theta)$.
\end{enumerate}
\end{lemma}
\begin{cproof}{prop:loss:diff:vc}
Throughout this proof let $ \phi \colon \R^{\dimension} \times [\scra, \scrb]^d \to \R $ satisfy for all $ \vartheta \in \R^{\dimension} $, $x \in \allowbreak [\scra, \allowbreak \scrb]^d $ that $ \phi( \vartheta, x ) = \abs{ \realapprox{\vartheta }{\infty} (x) - f(x) }^2$. \Nobs that the fact that $\Rect_\infty$ is Lipschitz continuous establishes that
\begin{equation}\label{eq:prop:loss:diff:vc:realization}
\R^{\dimension} \times [\scra, \scrb]^d \ni (\vartheta, x) \mapsto  \realapprox{\vartheta }{\infty} (x) \in \R
\end{equation}
is locally Lipschitz continuous. The fact that for all $ v, w \in \R $ it holds that $ v^2 - w^2 = (v-w)(v+w) $ \hence ensures that for all $ v, w \in \{ \vartheta \in \R^{\dimension} \colon \norm{ \theta - \vartheta } \le 1 \} $, $ x \in [\scra, \scrb]^d $ it holds that
\begin{equation}\label{eq:prop:loss:diff:vc:abs:phi}
\begin{split}
& \abs{ \phi( v, x ) - \phi( w, x ) } = \babs{ \bigl[ \realapprox{v}{\infty} (x) - f(x) \bigr]^2 - \bigl[ \realapprox{w}{\infty} (x) - f(x) \bigr]^2 } \\
& = \babs{ \bigl[ \realapprox{v}{\infty} (x) - \realapprox{w}{\infty} (x) \bigr] \bigl[ \realapprox{v}{\infty} (x) + \realapprox{w}{\infty} (x) - 2 f(x) \bigr] } \\
& \le 2 \abs{ \realapprox{v}{\infty} (x) - \realapprox{w}{\infty} (x) } \bigl[ \abs{ f(x) } + \sup\nolimits_{ \vartheta \in \R^{\dimension}, \norm{ \theta - \vartheta } \le 1 } \sup\nolimits_{ y \in [\scra, \scrb]^d } \abs{ \realapprox{\vartheta}{\infty} (y) } \bigr] < \infty \ifnocf.
\end{split}
\end{equation}
\cfload[.] Furthermore, \nobs that Hölder's inequality and the fact that for all $ \vartheta \in \R^{\dimension} $ it holds that
\begin{equation}\label{eq:prop:loss:diff:vc:risk:finite}
\textstyle \riskRR (\vartheta) = \int_{ [\scra, \scrb]^d } \abs{ \realapprox{\vartheta}{\infty} (x) - f(x) }^2 \, \mu( \d x ) = \int_{ [\scra, \scrb]^d } \phi( \vartheta, x ) \, \mu( \d x ) < \infty
\end{equation}
assure that for all $ \vartheta \in \R^{\dimension} $ it holds that
\begin{equation}\label{eq:prop:loss:diff:vc:integral:f}
\begin{split}
\textstyle \int_{ [\scra, \scrb]^d } \abs{ f(x) } \, \mu( \d x ) & \textstyle \le \bigl[ \mu( [\scra, \scrb]^d ) \bigr]^{\nicefrac{1}{2}} \bigl[\int_{ [\scra, \scrb]^d } \abs{ f(x) }^2 \, \mu(\d x)\bigr]^{\nicefrac{1}{2}} \\
& \le \textstyle \bigl[ \mu( [\scra, \scrb]^d ) \bigr]^{\nicefrac{1}{2}} \bigl[\int_{[\scra, \scrb]^d} \abs{ f(x) - \realapprox{\vartheta}{\infty} (x) }^2 \, \mu(\d x) \bigr]^{\nicefrac{1}{2}} \\
& \quad \textstyle + \bigl[ \mu( [\scra, \scrb]^d ) \bigr]^{\nicefrac{1}{2}} \bigl[\int_{ [\scra, \scrb]^d } \abs{ \realapprox{\vartheta}{\infty} (x) }^2 \, \mu( \d x )\bigr]^{\nicefrac{1}{2}} \\
& \le \textstyle \bigl[ \mu( [\scra, \scrb]^d ) \bigr]^{\nicefrac{1}{2}} \bigl[ \riskRR (\vartheta) \bigr]^{\nicefrac{1}{2}} + \bigl[ \mu( [\scra, \scrb]^d ) \bigr] \bigl[ \sup\nolimits_{ x \in [\scra, \scrb]^d } \abs{ \realapprox{\vartheta}{\infty} (x) } \bigr] < \infty.
\end{split}
\end{equation}
In addition, \nobs that the chain rule and the fact that for all $ i \in \N \cap (d \width + \width , \dimension ] $, $ x \in [\scra, \scrb]^d $ it holds that $ \R \ni v \mapsto \realapprox{ ( \theta_1, \dots, \theta_{i-1}, v, \theta_{i+1}, \dots, \theta_{\dimension} ) }{\infty} (x) \in \R $ is differentiable at $ \theta_i $ imply that for all $i \in \allowbreak \N \allowbreak \cap (d \width + \width, \dimension ] $, $ x \in [\scra, \scrb]^d $ it holds that $ \R \ni v \mapsto \phi( \theta_1, \dots, \theta_{i-1}, v, \theta_{i+1}, \dots, \theta_{\dimension}, x ) \in \R $ is differentiable at $ \theta_i $. Combining this, \cref{eq:prop:loss:diff:vc:realization}, \cref{eq:prop:loss:diff:vc:abs:phi}, \cref{eq:prop:loss:diff:vc:risk:finite}, \cref{eq:prop:loss:diff:vc:integral:f}, and \cref{cor:interchange} demonstrates that for all $ i \in \N \cap (d \width + \width , \dimension ] $ it holds that $ \R \ni v \mapsto \riskRR ( \theta_1, \dots, \theta_{i-1}, v, \theta_{i+1}, \dots, \theta_{\dimension} ) \in \R $ is differentiable at $ \theta_i $ and
\begin{equation}\label{eq:prop:loss:diff:vc:risk}
\bigl( \tfrac{\partial}{\partial \theta_i } \riskRR \bigr) ( \theta ) = \textstyle \int_{ [\scra, \scrb]^d } \bigl( \tfrac{\partial}{\partial \theta_i} \phi \bigr) ( \theta_1, \dots, \theta_{\dimension}, x ) \, \mu (\d x).
\end{equation}
Next \nobs that the chain rule shows for all $ j \in \{ 1, 2, \dots, \width \} $, $x \in [\scra, \scrb]^d$ that
\begin{equation}
\bigl(\tfrac{\partial}{\partial \theta_{d \width + \width + j}} \phi \bigr) ( \theta_1, \dots, \theta_{\dimension}, x ) = 2 ( \realapprox{\theta }{\infty} (x) - f(x) ) \bigl[ \Rect_\infty \bigl(\fb_j^{\theta} + \spro{\fw_j^{\theta}, x}\bigr) \bigr]
\end{equation}
and
\begin{equation}
\bigl( \tfrac{\partial}{\partial \theta_{\dimension}} \phi \bigr) ( \theta_1, \dots, \theta_{\dimension}, x ) = 2 ( \realapprox{\theta }{\infty} (x) - f(x) ).
\end{equation}
This and \cref{eq:prop:loss:diff:vc:risk} prove that for all $ j \in \{ 1, 2, \dots, \width \} $ it holds that
\begin{equation}
\bigl( \tfrac{\partial}{\partial \theta_{\width d + \width + j}} \riskRR \bigr) ( \theta ) = \textstyle 2 \int_{ [\scra, \scrb]^d } \bigl[\Rect_\infty \bigl(\fb_j^{\theta} + \spro{\fw_j^{\theta}, x}\bigr) \bigr] ( \realapprox{\theta }{\infty} (x) - f(x) ) \, \mu(\d x)
\end{equation}
and
\begin{equation}
\textstyle ( \frac{\partial}{\partial \theta_{\dimension} } \riskRR )(\theta ) = 2 \int_{ [\scra, \scrb]^d } ( \realapprox{\theta }{\infty} (x) - f(x) )  \, \mu(\d x).
\end{equation}
\cref{prop:limit:lr} \hence establishes for all $ i \in \N \cap (d \width + \width , \dimension ] $ that $ ( \frac{\partial}{\partial \theta_i } \riskRR )(\theta ) = \gradientG_i ( \theta)$.
\end{cproof}

\cfclear
\begin{lemma}[\resname{Differentiability w.r.t.\ inner weights and inner biases}]\label{prop:loss:diff:wb}
Assume \cref{setting:snn0}, assume $\mu \ll \lambda$, and let $\theta = (\theta_1, \dots, \theta_{\dimension}) \in \R^{ \dimension }$. Then
\begin{enumerate} [label = (\roman*)]
\item \label{prop:loss:diff:wb:item1} it holds for all $ j \in \{1, 2, \dots, \width \}$, $i \in \N \cap ((j - 1) d, j d] \cup \{d \width + j\}$ with $\abs{\fb_j^{\theta}} + \norm{\fw_j^{\theta}} >0 $ that $ \R \ni v \mapsto \riskRR ( \theta_1, \dots, \theta_{i-1}, v, \theta_{i+1}, \dots, \theta_{\dimension}) \in \R $ is differentiable at $ \theta_i $ and

\item \label{prop:loss:diff:wb:item2} it holds for all $ j \in \{1, 2, \dots, \width \} $, $ i \in \N \cap ( (j-1)d, jd ] \cup \{d \width + j\} $
with $\abs{\fb_j^{\theta}} + \norm{\fw_j^{\theta}} > 0$ that $( \frac{\partial}{\partial \theta_i } \riskRR )(\theta ) = \gradientG_{i}(\theta)$.
\end{enumerate}
\end{lemma}
\begin{cproof}{prop:loss:diff:wb}
Throughout this proof let $ \phi \colon \R^{\dimension} \times [\scra, \scrb]^d \to \R $ satisfy for all $ \vartheta \in \R^{\dimension} $, $ x \in [\scra, \scrb]^d $ that $ \phi( \vartheta, x ) = \abs{ \realapprox{\vartheta }{\infty} (x) - f(x) }^2 $ and for every $j \in \{ 1, 2, \dots, \width \}$ let $E_j \subseteq \R^d$ satisfy
\begin{equation}\label{prop:loss:diff:wb:def:E}
E_j = \bigl\{ x \in [\scra, \scrb]^d \colon \fb_j^{\theta} + \spro{\fw_j^{\theta}, x} \neq 0 \bigr\}.
\end{equation}
\Nobs that the integral transformation theorem ensures for all $\vartheta \in \R^d \backslash \{ 0 \}$, $c \in \R$ that
\begin{equation}\label{prop:loss:diff:wb:integral:hyperplane}
\begin{split}
\textstyle \int_{ \R^d } \indicator{ \{ x \in \R^d \colon c + \scalprod{ \vartheta }{ x } = 0 \} } (y) \, \d y
&= \textstyle \int_{ \R^d } \indicator{ \{ 0 \} } ( c + \scalprod{ \vartheta }{ y } ) \, \d y = \int_{ \R^d } \indicator{ \{ 0 \} } ( \scalprod{ \vartheta }{ y + c \norm{ \vartheta }^{-2} \vartheta } ) \, d y \\
&= \textstyle \int_{ \R^d } \indicator{ \{ 0 \} } ( \scalprod{ \vartheta }{ y } ) \, d y = \int_{ \R^d } \indicator{ \{ x \in \R^d \colon \scalprod{ \vartheta }{ x } = 0 \} } (y) \, \d y \ifnocf.
\end{split}
\end{equation}
\cfload[.]Moreover, \nobs that the rank-nullity theorem demonstrates for all $\vartheta \in \R^d \backslash \{ 0 \}$ that
\begin{equation}
\dim_{ \R } ( \{ x \in \R^d \colon \scalprod{ \vartheta }{ x } = 0 \} ) = d - \dim_{ \R } ( \{ y \in \R \colon [ \Exists x \in \R^d \colon y = \scalprod{ \vartheta }{ x } ] \} ) = d - 1.
\end{equation}
\Hence for all $ \vartheta \in \R^d \backslash \{ 0 \} $ that $ \int_{ \R^d } \indicator{ \{ x \in \R^d \colon \scalprod{ \vartheta }{ x } = 0 \} } (y) \, \d y = 0 $. Combining this with \cref{prop:loss:diff:wb:integral:hyperplane} and the fact that for all $ \vartheta \in \R^d $, $ c \in \R $ it holds that $ \{ x \in [\scra, \scrb]^d \colon c + \scalprod{ \vartheta }{ x } = 0 \} \subseteq \{ x \in \R^d \colon c + \scalprod{ \vartheta }{ x } = 0 \} $ shows that for all $ \vartheta \in \R^d \backslash \{ 0 \} $, $ c \in \R $ it holds that $ \lambda( \{ x \in [\scra, \scrb]^d \colon c + \scalprod{ \vartheta }{ x } = 0 \} ) \le \int_{ \R^d } \indicator{ \{ x \in \R^d \colon c + \scalprod{ \vartheta }{ x } = 0 \} } ( y ) \, \d y = 0 $. This and \cref{prop:loss:diff:wb:def:E} prove that for all $ j \in \{ 1, 2, \dots, \width \} $ with $\abs{\fb_j^{\theta}} + \norm{\fw_j^{\theta}} > 0 $ it holds that
\begin{equation}\label{eq:prop:loss:diff:wb:E}
\lambda ( [\scra, \scrb]^d \backslash E_j ) = 0.
\end{equation}
Next \nobs that the fact that $ \Rect_\infty $ is Lipschitz continuous implies that
\begin{equation}\label{eq:prop:loss:diff:wb:realization}
\R^{ \dimension } \times [\scra, \scrb]^d \ni (\vartheta, x) \mapsto \realapprox{\vartheta}{\infty} ( x )  \in \R
\end{equation}
is locally Lipschitz continuous. The fact that for all $ v, w \in \R $ it holds that $ v^2 - w^2 = (v-w)(v+w) $ \hence shows that for all $v, w \in \{ \vartheta \in \R^{\dimension} \colon \norm{ \theta - \vartheta } \le 1 \} $, $ x \in [\scra, \scrb]^d$ it holds that
\begin{equation}\label{eq:prop:loss:diff:wb:abs:phi}
\begin{split}
& \abs{ \phi( v, x ) - \phi( w, x ) } = \babs{ \bigl[ \realapprox{v}{\infty} (x) - f(x) \bigr]^2 - \bigl[ \realapprox{w}{\infty} (x) - f(x) \bigr]^2 } \\
& = \babs{ \bigl[ \realapprox{v}{\infty} (x) - \realapprox{w}{\infty} (x) \bigr] \bigl[ \realapprox{v}{\infty} (x) + \realapprox{w}{\infty} (x) - 2 f(x) \bigr] } \\
& \le 2 \abs{ \realapprox{v}{\infty} (x) - \realapprox{w}{\infty} (x) } \bigl[ \abs{ f(x) } + \sup\nolimits_{ \vartheta \in \R^{\dimension}, \norm{ \theta - \vartheta } \le 1 } \sup\nolimits_{ y \in [\scra, \scrb]^d } \abs{ \realapprox{\vartheta}{\infty} (y) } \bigr] < \infty.
\end{split}
\end{equation}
Moreover, \nobs that Hölder's inequality and the fact that for all $ \vartheta \in \R^{\dimension} $ it holds that
\begin{equation}\label{eq:prop:loss:diff:wb:risk:finite}
\textstyle \riskRR (\vartheta) = \int_{ [\scra, \scrb]^d } \abs{ \realapprox{\vartheta}{\infty} (x) - f(x) }^2 \, \mu( \d x ) = \int_{ [\scra, \scrb]^d } \phi( \vartheta, x ) \, \mu( \d x ) < \infty
\end{equation}
prove that for all $ \vartheta \in \R^{\dimension} $ it holds that
\begin{equation}\label{eq:prop:loss:diff:wb:integral:f}
\begin{split}
\textstyle \int_{ [\scra, \scrb]^d } \abs{ f(x) } \, \mu( \d x ) & \textstyle \le \textstyle \bigl[ \mu( [\scra, \scrb]^d ) \bigr]^{\nicefrac{1}{2}}  \bigl[\int_{ [\scra, \scrb]^d } \abs{ f(x) }^2 \, \mu(\d x)\bigr]^{\nicefrac{1}{2}} \\
& \le \textstyle \bigl[ \mu( [\scra, \scrb]^d ) \bigr]^{\nicefrac{1}{2}} \bigl[\int_{[\scra, \scrb]^d} \abs{ f(x) - \realapprox{\vartheta}{\infty} (x) }^2 \, \mu(\d x)\bigr]^{\nicefrac{1}{2}} \\
& \quad \textstyle + \bigl[ \mu( [\scra, \scrb]^d ) \bigr]^{\nicefrac{1}{2}} \bigl[\int_{ [\scra, \scrb]^d } \abs{ \realapprox{\vartheta}{\infty} (x) }^2 \, \mu(\d x)\bigr]^{\nicefrac{1}{2}} \\
& \le \textstyle \bigl[ \mu( [\scra, \scrb]^d ) \bigr]^{\nicefrac{1}{2}} \bigl[ \riskRR (\vartheta) \bigr]^{\nicefrac{1}{2}} + \bigl[ \mu( [\scra, \scrb]^d ) \bigr] \bigl[ \sup\nolimits_{ x \in [\scra, \scrb]^d } \abs{ \realapprox{\vartheta}{\infty} (x) } \bigr] < \infty.
\end{split}
\end{equation}
In addition, \nobs that the chain rule and the fact that for all $ x \in \R \backslash \{ 0 \} $ it holds that $ \Rect_\infty $ is differentiable at $ x $ and $ (\Rect_\infty )' ( x ) = \indicator{(0, \infty)} ( x ) $ demonstrate that for all $ j \in \{1, 2, \dots, \width \}$, $i \in \N \cap ((j-1)d, jd] \cup \{d \width + j\}$, $x \in E_j$ it holds that $ \R \ni v \mapsto \realapprox{(\theta_1, \dots, \theta_{i-1}, v, \theta_{i+1}, \dots, \theta_{ \dimension })}{\infty} ( x ) \in \R $ is differentiable at $ \theta_i $. This ensures that for all $j \in \{1, 2, \dots, \width\}$, $i \in \N \cap ((j-1)d, jd] \cup \{d \width + j\}$, $x \in E_j$ it holds that
\begin{equation}
\R \ni v \mapsto \phi(\theta_1, \dots, \theta_{i-1}, v, \theta_{i+1}, \dots, \theta_{ \dimension }, x) \in \R
\end{equation}
is differentiable at $\theta_i$. \cref{cor:interchange}, \cref{eq:prop:loss:diff:wb:E}, \cref{eq:prop:loss:diff:wb:realization}, \cref{eq:prop:loss:diff:wb:abs:phi}, \cref{eq:prop:loss:diff:wb:risk:finite}, \cref{eq:prop:loss:diff:wb:integral:f}, and the assumption that $\mu \ll \lambda$ \hence establish for all $j \in \{1, 2, \dots, \width\}$, $i \in \N \cap \allowbreak ((j-1)d, jd] \allowbreak \cup \{d \width + j\}$ with $\abs{\fb_j^{\theta}} + \norm{\fw_j^{\theta}} > 0$ that $\R \ni v \mapsto \riskRR ( \theta_1, \dots, \theta_{i-1}, v, \allowbreak \theta_{i+1}, \dots, \theta_{\dimension}) \in \R$ is differentiable at $ \theta_i $ and
\begin{equation}\label{eq:prop:loss:diff:wb:risk}
\bigl( \tfrac{\partial}{\partial \theta_i } \riskRR \bigr) ( \theta ) = \textstyle \int_{ E_j } \bigl( \tfrac{\partial}{\partial \theta_i} \phi \bigr) ( \theta_1, \dots, \theta_{\dimension}, x ) \, \mu (\d x).
\end{equation}
Next \nobs that the chain rule shows for all $j \in \{1, 2, \dots, \width\}$, $k \in \{1, 2, \dots, d\}$, $x \allowbreak \in E_j$ that
\begin{equation}
\begin{split}
& \bigl(\tfrac{\partial}{\partial \theta_{(j - 1) d + k}} \phi \bigr) ( \theta_1, \dots, \theta_{\dimension}, x ) \\
& = 2 (\realapprox{\theta }{\infty}(x) - f(x)) \fv_j^{\theta} \spro{e_k, x} \indicator{(0, \infty)} \bigl(\fb_j^{\theta} + \spro{\fw_j^{\theta}, x}\bigr) = 2 (\realapprox{\theta }{\infty}(x) - f(x)) \fv_j^{\theta} \spro{e_k, x} \indicator{I_j^\theta}(x)
\end{split}
\end{equation}
and
\begin{equation}
\begin{split}
& \bigl(\tfrac{\partial}{\partial \theta_{d \width + j}} \phi\bigr) (\theta_1, \dots, \theta_{\dimension}, x) \\
& = 2 (\realapprox{\theta}{\infty}(x) - f(x)) \fv_j^{\theta} \indicator{(0, \infty)} \bigl(\fb_j^{\theta} + \spro{\fw_j^{\theta}, x}\bigr) = 2 (\realapprox{\theta }{\infty}(x) - f(x)) \fv_j^{\theta} \indicator{I_j^\theta}(x).
\end{split}
\end{equation}
This and \cref{eq:prop:loss:diff:wb:risk} prove for all $j \in \{1, 2, \dots, \width\}$, $k \in \{1, 2, \dots, d\}$ with $\abs{\fb_j^{\theta}} + \norm{\fw_j^{\theta}} > 0$ that
\begin{equation}
\begin{split}
\bigl(\tfrac{\partial}{\partial \theta_{(j-1)d+k}} \riskRR\bigr) (\theta) & = \textstyle 2 \fv_j^{\theta} \int_{ E_j } \scalprod{ e_k }{ x } (\realapprox{\theta }{\infty} (x) - f(x) ) \indicator{I_j^\theta} ( x ) \, \mu (\d x) \\
& = \textstyle 2 \fv_j^{\theta} \int_{ I_j^\theta } \scalprod{ e_k }{ x } ( \realapprox{\theta }{\infty} (x) - f(x) ) \, \mu (\d x)
\end{split}
\end{equation}
and
\begin{equation}
\begin{split}
\bigl(\tfrac{\partial}{\partial \theta_{d \width + j}} \riskRR\bigr) (\theta) & = \textstyle 2 \fv_j^{\theta} \int_{E_j} (\realapprox{\theta}{\infty}(x) - f(x)) \indicator{I_j^\theta}(x) \, \mu(\d x) = \textstyle 2 \fv_j^{\theta} \int_{I_j^\theta} (\realapprox{\theta }{\infty}(x) - f(x)) \, \mu(\d x).
\end{split}
\end{equation}
\cref{prop:limit:lr} \hence establishes for all $j \in \{1, 2, \dots, \width\}$, $i \in \N \cap ((j-1)d, jd] \cup \{d \width + j\}$ with $\abs{\fb_j^{\theta}} + \norm{\fw_j^{\theta}} > 0$ that $(\frac{\partial}{\partial \theta_i} \riskRR)(\theta) = \gradientG_i(\theta)$.
\end{cproof}

\subsection{Local extrema and saddle points}
\label{subsec:local_extrema_and_saddle_points}

\cfclear
\begin{definition}[\resname{Local minima}]\label{def:loc_min}
Let $\dimension \in \N$, let $U \subseteq \R^{\dimension}$ be a set, let $f \colon U \to \R$ be a function, and let $x \allowbreak \in U $. Then we say that $x$ is a local minimum point of $f$ if and only if there exists $\varepsilon \in (0, \infty)$ such that $f(x) = \inf_{y \in \{v \in U \colon \norm{v - x} \le \varepsilon\}} f(y)$ \cfout.
\end{definition}

\cfclear
\begin{definition}[\resname{Local maxima}]\label{def:loc_max}
Let $\dimension \in \N$, let $U \subseteq \R^{\dimension}$ be a set, let $f \colon U \to \R$ be a function, and let $x \in U$. Then we say that $ x $ is a local maximum point of $ f $ if and only if there exists $\varepsilon \in (0, \infty)$ such that $f(x) = \sup_{y \in \{v \in U \colon \norm{v - x} \le \eps\}} f(y)$ \cfout.
\end{definition}

\cfclear
\begin{definition}[\resname{Local extrema}]\label{def:loc_extrema}
Let $\dimension \in \N$, let $U \subseteq \R^{\dimension}$ be a set, let $f \colon U \to \R$ be a function, and let $x \allowbreak \in U$. Then we say that $x$ is a local extremum point of $f$ if and only if $x \in \allowbreak \{ y \in \allowbreak \R^{\dimension} \colon \allowbreak (y \text{ is a \locmin point of } f)\} \cup \allowbreak \{y \in \allowbreak \R^{\dimension} \colon \allowbreak (y \text{ is a \locmax point of } f)\}$ \cfout.
\end{definition}

\cfclear
\begin{definition}[\resname{Saddle points}]\label{def:saddle}
Let $\dimension \in \N$, let $U \subseteq \R^{\dimension}$ be open, let $f \colon U \to \R$ be a function, let $x \in U$, and assume that $f$ is differentiable at $x$. Then we say that $x$ is a saddle point of $f$ if and only if we have that
\begin{enumerate}[label = (\roman*)]
\item it holds that $x$ is not a \locextrema point of $f$ and
\item it holds that $(\nabla f)(x) = 0$ \ifnocf.
\end{enumerate}
\cfout[.]
\end{definition}

\subsection{Values of the generalized gradient function at local minimum points}
\label{subsec:values_of_general_grad_at_loc_min}

\cfclear
\begin{proposition}\label{prop:local:minima:gradient}
Assume \cref{setting:snn0}, assume $\mu \ll \lambda$, and let $\theta = (\theta_1, \dots, \theta_{\dimension}) \in \R^\dimension$ be a \locmin point of $ \riskRR $ \cfload. Then $\gradientG(\theta) = 0$.
\end{proposition}
\begin{cproof}{prop:local:minima:gradient}
\Nobs that \cref{prop:loss:diff:vc} ensures for all $i \in \N \cap (d \width + \width, \dimension]$ that $\R \ni v \allowbreak \mapsto \allowbreak \riskRR (\theta_1, \dots, \theta_{i-1}, v, \theta_{i+1}, \dots, \theta_{\dimension}) \in \R$ is differentiable at $\theta_i $ and $ ( \frac{\partial}{\partial \theta_i } \riskRR )(\theta ) = \gradientG_i (\theta)$. The assumption that $\theta$ is a \locmin point of $\riskRR$ \hence shows for all $i \in \N \cap (d \width + \width, \dimension]$ that
\begin{equation}\label{eq:prop:local:minima:gradient:vc}
\gradientG_i ( \theta) = ( \tfrac{\partial}{\partial \theta_i } \riskRR )(\theta ) = 0 \ifnocf.
\end{equation}
\cfload[.]Moreover, \nobs that \cref{prop:loss:diff:wb} ensures for all $j \in \{1, 2, \dots, \width\}$, $i \in \N \cap \allowbreak ((j-1)d, \allowbreak jd] \allowbreak \cup \allowbreak \{d \width \allowbreak + j\}$ with $\abs{\fb_j^{\theta}} + \norm{\fw_j^{\theta}} > 0$ that $ \R \ni v \mapsto \riskRR (\theta_1, \dots, \theta_{i-1}, v, \allowbreak \theta_{i+1}, \dots, \theta_{\dimension} ) \in \R $ is differentiable at $ \theta_i $ and $(\frac{\partial}{\partial \theta_i} \riskRR )(\theta ) \allowbreak = \gradientG_i(\theta)$. This and the assumption that $\theta$ is a \locmin point of $\riskRR$ imply for all $j \in \{1, 2, \dots, \width\}$, $i \in \N \cap ((j-1)d, jd] \cup \{d \width + j\}$ with $\abs{\fb_j^{\theta}} + \norm{\fw_j^{\theta}} > 0$ that
\begin{equation}\label{eq:prop:local:minima:gradient:wb}
\gradientG_i(\theta) = (\tfrac{\partial}{\partial \theta_i} \riskRR)(\theta) = 0.
\end{equation}
In addition, \nobs that \cref{prop:limit:lr} and the fact that for all $j \in \{1, 2, \dots, \width\}$ with $\abs{\fb_j^{\theta}} \allowbreak + \allowbreak \norm{\fw_j^{\theta}} = 0$ it holds that $I_j^\theta = \emptyset$ demonstrate for all $j \in \{1, 2, \dots, \width\}$, $i \allowbreak \in \allowbreak \N \allowbreak \cap \allowbreak ((j-1)d, \allowbreak jd] \allowbreak \cup \{d \width \allowbreak + j\}$ with $\abs{\fb_j^{\theta}} \allowbreak + \allowbreak \norm{\fw_j^{\theta}} = 0$ that $\gradientG_i(\theta) = 0$. This and \cref{eq:prop:local:minima:gradient:wb} assure for all $j \allowbreak \in \allowbreak \{1, \allowbreak 2, \allowbreak \dots, \allowbreak \width\}$, $i \in \allowbreak \N \allowbreak \cap ((j-1)d, \allowbreak jd] \allowbreak \cup \{d \width \allowbreak + j\}$ that $\gradientG_i(\theta) = 0$. \Hence for all $i \in \N \cap (0, d \width + \width]$ that $\gradientG_i(\theta) \allowbreak = 0$. Combining this with \cref{eq:prop:local:minima:gradient:vc} establishes for all $i \in \{1, 2, \dots, \dimension\}$ that $\gradientG_i(\theta) = 0$.
\end{cproof}

\section{Differential geometric preliminaries}
\label{sec:differential_geometric_prelimiaries}

This section is devoted to establishing some essentially well-known necessary and sufficient conditions for local extremum and saddle points in \cref{prop:class:crit:points} in \cref{subsec:necessary_and_sufficient_conditions} below. Our proof of \cref{prop:class:crit:points} uses the well-known rank bound for the Hessian matrix in \cref{lem:rank:bound} in \cref{subsec:tangent_spaces_associated_to_submanifolds} below, whose proof can be found, e.g., in \cite[Chapter~1, Section~2]{GuilleminPollack2010DiffGeo}, the essentially well-known sufficient condition for a local minimum point in \cref{lem:eigenvalues:locmin:point}, and the well-known necessary condition for a local minimum point in \cref{lem:eigenvalues:nonnegative} in \cref{subsec:necessary_and_sufficient_conditions} below.
In the proof of \cref{lem:eigenvalues:locmin:point} we employ the well-known Taylor-type estimate from \cref{lem:mult:taylor} in \cref{subsec:necessary_and_sufficient_conditions}.

For the convenience of the reader we also recall in this section the notion of the spectrum of a matrix as well as some basic differential geometric concepts such as the notions of an immersion, of a submanifold, of the unique projection on a nonempty set, and of the tangent space; see \cref{def:immersion,def:submanifold,def:Rd_subset_unique_projection,def:max_open_set_including_M,def:tangent_space,def:spectrum}.

\subsection{Immersions}

\cfclear
\begin{definition}[\resname{Immersion}]\label{def:immersion}
Let $\dimension, k \in \N$, $n \in \N \cup \{\infty\}$ and let $U \subseteq \R^k$ be open. Then we say that $\varphi$ is a $C^n$-immersion from $U$ to $\R^{\dimension}$ if and only if we have that
\begin{enumerate}[label = (\roman*)]
\item it holds that $\varphi \in C^n(U, \R^{\dimension})$ and 

\item it holds for all $x \in U$ that $\rk(\varphi'(x)) = k$.
\end{enumerate}
\end{definition}

\subsection{Submanifolds of Euclidean spaces}

\cfclear
\begin{definition}[\resname{Submanifold}]\label{def:submanifold}
Let $\dimension, k \in \N$, $n \in \N \cup \{\infty\}$ \cfload. Then we say that $\cM$ is a $k$-dimensional $C^n$-submanifold of $\R^{\dimension}$ if and only if it holds for all $x \in \cM$ that there exist $U \in \{V \subseteq \R^k \colon V \text{ is open}\}$, $\eps \in (0, \infty)$, $\varphi \in C (U, \R^{\dimension})$ such that
\begin{enumerate}[label = (\roman*)] \cfadd{def:immersion}
\item it holds that $\cM \subseteq \R^{\dimension}$,

\item it holds that $\varphi$ is a $C^n$-immersion from $U$ to $\R^{\dimension}$,

\item it holds that $\varphi(U) = \cM \cap \{y \in \R^{\dimension} \colon \norm{x - y} < \eps\}$, and

\item it holds that $U \ni y \mapsto \varphi(y) \in \varphi(U)$ is a homeomorphism\ifnocf.
\end{enumerate}
\cfout[.]
\end{definition}

\subsection{Nonlinear projections}

\cfclear
\begin{definition}\label{def:Rd_subset_unique_projection}
\cfconsiderloaded{def:Rd_subset_unique_projection} Let $\dimension \in \N$ and let $\cM \subseteq \R^{\dimension}$ satisfy $\cM \neq \emptyset$. Then we denote by $\PP_{\cM} \subseteq \R^{\dimension}$ the set given by
\begin{equation}
\PP_{\cM} = \{x \in \R^{\dimension} \colon (\exists_{1} \, y \in \cM \colon \norm{x - y} = \inf\nolimits_{z \in \cM} \norm{x - z})\}
\end{equation}
and we denote by $\pp_{\cM} \colon \PP_{\cM} \to \R^{\dimension}$ the function which satisfies for all $x \in \PP_{\cM}$ that $\pp_{\cM}(x) \in \cM$ and
\begin{equation}
\norm{x - \pp_{\cM}(x)} = \inf\nolimits_{y \in \cM} \norm{x - y}\ifnocf.
\end{equation}
\cfload[.]
\end{definition}

\cfclear
\begin{definition}\label{def:max_open_set_including_M}
\cfconsiderloaded{def:max_open_set_including_M} Let $\dimension, k \in \N$, let $\cM \subseteq \R^{\dimension}$ be a $k$-dimensional $C^2$-\submanifold of $\R^{\dimension}$, and assume $\cM \neq \emptyset$ \cfload. Then we denote by $\bfPP_{\cM} \subseteq \R^{\dimension}$ the set given by
\begin{equation}
\textstyle \bfPP_{\cM} = \bigcup_{\substack{U \subseteq \R^{\dimension} \; \text{is open}, \; U \subseteq \PP_{\cM}, \\ \text{and} \; \pp_{\cM}|_{U} \in C^1(U, \R^{\dimension})}} U\ifnocf.
\end{equation}
\cfload[.]
\end{definition}

\subsection{Tangent spaces associated to submanifolds of Euclidean spaces}
\label{subsec:tangent_spaces_associated_to_submanifolds}

\cfclear
\begin{definition}[\resname{Tangent space}]\label{def:tangent_space}
\cfconsiderloaded{def:tangent_space} Let $\dimension \in \N$, let $\cM \subseteq \R^{\dimension}$ be a set, and let $x \in \cM$ \cfload. Then we denote by $\tangentspace{x}{\cM}$ the set given by
\begin{equation}
\tangentspace{x}{\cM} = \{v \in \R^{\dimension} \colon (\Exists \gamma \in C^1(\R, \R^{\dimension}) \colon ([\gamma(\R) \subseteq \cM] \wedge [\gamma(0) = x] \wedge [\gamma'(0) = v]))\}\ifnocf.
\end{equation}
\cfout[.]
\end{definition}

\cfclear
\begin{lemma}\label{lemma:tangent_space}
Let $\dimension, k \in \N$, let $\cM$ be a $k$-dimensional $C^1$-\submanifold of $\R^{\dimension}$, and let $x \in \cM$ \cfload. Then it holds that $\tangentspace{x}{\cM}$ is a $k$-dimensional vector subspace of $\R^{\dimension}$\cfout.
\end{lemma}
\begin{cproof}{lemma:tangent_space}
\Nobs that, e.g., \cite[Chapter~1, Section~2]{GuilleminPollack2010DiffGeo} ensures that $\tangentspace{x}{\cM}$ is a $k$-di\-men\-si\-onal vector subspace of $\R^{\dimension}$\cfload.
\end{cproof}

\cfclear
\begin{lemma} \label{lem:rank:bound}
Let $\dimension, k \in \N$, let $U \subseteq \R^{\dimension}$ be open, let $f \in C^2(U, \R)$, let $\cM \subseteq U$ satisfy $\cM \subseteq \{ x \in U \colon (\nabla f) (x) \allowbreak = 0\}$, assume that $\cM$ is a $k$-dimensional $C^1$-\submanifold of $\R^{\dimension}$, and let $x \in \cM$ \cfload. Then
\begin{enumerate} [label = (\roman*)]
\item \label{lem:rank:bound:item1} it holds for all $v \in \tangentspace{x}{\cM}$ that $((\Hs f)(x)) v = 0$ and

\item \label{lem:rank:bound:item2} it holds that  $\rk((\Hs f)(x)) \leq \dimension - k$\ifnocf.
\end{enumerate}
\cfout[.]
\end{lemma}
\begin{cproof}{lem:rank:bound}
Throughout this proof let $v \in \tangentspace{x}{\cM}$ and let $\gamma \in C^1(\R, \allowbreak \cM)$ satisfy $\gamma(0) = x$ and $\gamma'(0) = v$ \cfload. \Nobs that the fact that for all $y \in \cM$ it holds that $(\nabla f)(y) =0$ proves that
\begin{equation} 
0 = \tfrac{\d}{\d t}(\nabla f)(\gamma(t))|_{t=0} = ((\Hs f)(x)) \gamma'(0) = ((\Hs f)(x)) v.
\end{equation}
\Hence for all $\fv \in \tangentspace{x}{\cM}$ that $((\Hs f)(x)) \fv = 0$. \cref{lemma:tangent_space} \hence demonstrates that $\rk ( (\Hs  f ) ( x ) ) \leq \dimension - k$.
\end{cproof}

\subsection{Necessary and sufficient conditions for local extremum and saddle points}
\label{subsec:necessary_and_sufficient_conditions}

\cfclear
\begin{lemma}\label{lem:mult:taylor}
Let $\dimension \in \N$, let $U \subseteq \R^{\dimension}$ be open, let $f \in C^2(U, \R)$, and let $K \subseteq U$ be compact. Then it holds for all $\varepsilon \in (0, \infty)$ that there exists $\delta \in (0, \infty)$ such that for all $x, y \in K$ with $\norm{x - y} \allowbreak < \delta$ and $\cup_{t \in [0, 1]} \{(1 - t) x + t y\} \allowbreak \subseteq K$ it holds that
\begin{equation}\label{eq:lem:mult:taylor:quadratic}
\textstyle \abs{f(y) - f(x) - \scalprod{(\nabla f)(x)}{y - x} - \frac{1}{2}\scalprod{y - x}{(\Hs f)(x)(y - x)}} \leq \varepsilon \norm{x - y}^2 \ifnocf.
\end{equation}
\cfout[.]
\end{lemma}
\begin{cproof}{lem:mult:taylor}
Throughout this proof for every $x, y \in K$ let $J_{x, y} \allowbreak \subseteq \R$ satisfy $J_{x,y} = \{t \in \R \colon (1 - t) x + t y \allowbreak \in U\}$ and let $\varphi_{x,y} \colon J_{x,y} \to \R$ satisfy for all $t \in J_{x,y}$ that $\varphi_{x,y}(t) = f((1 - t) x + t y)$. \Nobs that the assumption that $U$ is open, the assumption that $f \in C^2(U,\R)$, and the chain rule ensure for all $x,y \in K$ that $J_{x,y}$ is open, that $\varphi_{x,y} \in C^2(J_{x,y},\R)$, and that for all $t \in J_{x,y}$ it holds that
\begin{equation}\label{eq:lem:mult:taylor:first:derivative}
(\varphi_{x,y})'(t)=\scalprod{(\nabla f)((1-t)x+ty)}{y-x}
\end{equation}
and
\begin{equation}\label{eq:lem:mult:taylor:second:derivative}
(\varphi_{x,y})''(t)=\scalprod{y-x}{(\Hs f)((1-t)x+ty)(y-x)}.
\end{equation}
Moreover, \nobs that Taylor's theorem proves for all $x,y \in K$ with $\cup_{t \in [0,1]}\{(1-t)x+ty\} \subseteq K$ that
\begin{equation}\label{eq:lem:mult:taylor:phi}
\textstyle \varphi_{x,y}(1) = \varphi_{x,y}(0) + (\varphi_{x,y})'(0) + \int_0^1(1-t)[(\varphi_{x,y})''(t)]\,\d t.
\end{equation}
Combining this with \cref{eq:lem:mult:taylor:first:derivative,eq:lem:mult:taylor:second:derivative} shows for all $x, y \allowbreak \in K$ with $\cup_{t \in [0, 1]} \{(1 - t) x + t y\} \subseteq K$ that
\begin{equation}\label{eq:lem:mult:taylor:remainder}
\begin{split}
& \textstyle f(y)-f(x)-\scalprod{(\nabla f)(x)}{y-x}-\frac{1}{2}\scalprod{y-x}{(\Hs f)(x)(y-x)} \\
& \textstyle = \varphi_{x,y}(1) - \varphi_{x,y}(0) - (\varphi_{x,y})'(0) - \frac{1}{2} (\varphi_{x,y})''(0) \\
& \textstyle = \int_0^1(1-t)[(\varphi_{x,y})''(t) - (\varphi_{x,y})''(0)]\,\d t \\
& \textstyle = \int_0^1(1-t)\scalprod{y-x}{[(\Hs f)((1-t)x+ty)-(\Hs f)(x)](y-x)}\,\d t.
\end{split}
\end{equation}
In the next step let $\varepsilon \in (0,\infty)$. \Nobs that the fact that the second-order partial derivatives of $f$ are continuous and the fact that $K$ is compact ensure that there exists $\delta \in (0, \infty)$ which satisfies for all $x \allowbreak = (x_1, \allowbreak \dots, \allowbreak x_{\dimension}), \allowbreak y \allowbreak = (y_1, \allowbreak \dots, \allowbreak y_{\dimension}) \allowbreak \in K$, $i, \allowbreak j \allowbreak \in \{1, \allowbreak 2, \allowbreak \dots, \allowbreak \dimension\}$ with $\norm{x - y} \allowbreak < \delta$ that
\begin{equation}\label{eq:lem:mult:taylor:uniform}
\textstyle \abs{(\frac{\partial^2}{\partial x_i\partial x_j}f)(x)-(\frac{\partial^2}{\partial x_i\partial x_j}f)(y)}\leq \frac{2\varepsilon}{\dimension}.
\end{equation}
The Cauchy-Schwarz inequality, \cref{eq:lem:mult:taylor:remainder}, and \cref{eq:lem:mult:taylor:uniform} \hence demonstrate for all $x, y \allowbreak \in K$ with $\norm{x - y} \allowbreak < \delta$ and $\cup_{t \in [0, 1]} \{(1 - t) x + t y\} \allowbreak \subseteq K$ that
\begin{equation}
\begin{split}
&\textstyle \abs{f(y)-f(x)-\scalprod{(\nabla f)(x)}{y-x}-\frac{1}{2}\scalprod{y-x}{(\Hs f)(x)(y-x)}}\\
&\textstyle \leq 2\varepsilon\norm{x-y}^2\int_0^1(1-t)\,\d t = \varepsilon \norm{x - y}^2.
\end{split}
\end{equation}
\end{cproof}

\cfclear
\begin{definition}[\resname{Real spectrum}]\label{def:spectrum}
\cfconsiderloaded{def:spectrum} Let $\dimension \in \N$, $A \in \R^{\dimension \times \dimension}$. Then we denote by $\spectrum{A} \subseteq \R$ the set given by
\begin{equation}
\spectrum{A} = \{\lambda \in \R \colon (\exists \, v \in \R^{\dimension} \backslash \{ 0 \} \colon Av = \lambda v)\}\ifnocf.
\end{equation}
\cfload[.]
\end{definition}

\cfclear
\begin{lemma}[\resname{Sufficient condition for local minima}]\label{lem:eigenvalues:locmin:point}
Let $\dimension, \allowbreak k \allowbreak \in \N$, let $U \allowbreak \subseteq \R^{\dimension}$ be open, let $f \allowbreak \in C^2(U, \allowbreak \R)$, let $\cM \allowbreak \subseteq U$ satisfy $\cM \allowbreak \subseteq \{x \allowbreak \in U \colon (\nabla f)(x) \allowbreak = 0\}$, assume that $\cM$ is a $k$-dimensional $C^2$-\submanifold of $\R^{\dimension}$, and let $x \in \cM$ satisfy $\rk((\Hs f)(x)) \allowbreak = \dimension \allowbreak - k $ and $\spectrum{(\Hs  f)(x)} \allowbreak \subseteq [0, \allowbreak \infty)$ \cfload. Then it holds that $x$ is a \locmin point of $f$ \cfout.
\end{lemma}
\begin{cproof}{lem:eigenvalues:locmin:point}
Throughout this proof for every $(R, \delta) \in (0, \infty)^2$ let $\tubular{}{R, \delta} \subseteq \R^{\dimension}$ satisfy
\begin{multline}
\textstyle \tubular{}{R, \delta} = \bigl\{ z \in \R^{\dimension} \colon [\, \Exists y \in \cM \cap \{w \in \R^{\dimension} \colon \norm{x - w} \le R\} \colon \\ 
\textstyle (\Exists v \in (\tangentspace{y}{\cM})^{\perp} \cap \{w \in \R^{\dimension} \colon \norm{w} < \delta\} \colon z = y + v)] \bigr\}\ifnocf.
\end{multline}
\cfload[.]In the following we distinguish between the case $ k = \dimension $ and the case $ k < \dimension $. We first prove in the case
\begin{equation}\label{eq:lem:eigenvalues:locmin:point:case1}
k = \dimension
\end{equation}
that $ x $ is a \locmin point of $ f $ \cfload. \Nobs that \cref{eq:lem:eigenvalues:locmin:point:case1} and the assumption that $ \cM $ is a $ k $-dimensional $ C^2 $-\submanifold of $ \R^{\dimension} $ ensure that there exists $ r \in (0,\infty) $ which satisfies $ \{ w \in \R^{ \dimension } \colon \norm{ x - w } < r \} \subseteq \cM $. The fact that for all $ y \in \cM $ it holds that $ ( \nabla f )(y) = 0 $ \hence implies that for all $ y \in \{ w \in \R^{ \dimension } \colon \norm{ x - w } < r \} $ it holds that $ f(y) = f(x) $. This proves in the case $ k = \dimension $ that $ x $ is a \locmin point of $ f $. In the next step we prove in the case
\begin{equation}\label{eq:lem:eigenvalues:locmin:point:case2}
k < \dimension
\end{equation}
that $ x $ is a \locmin point of $ f $. \Nobs that the assumption that $ \rk((\Hs f)(x)) = \dimension - k $ and \cref{eq:lem:eigenvalues:locmin:point:case2} demonstrate that $ \rk((\Hs f)(x)) > 0 $. Moreover, \nobs that the assumption that $ f \in C^2( U, \R ) $ ensures that for all $ y \in U $ it holds that $ ( \Hs f )(y) $ is symmetric. \cref{lemma:tangent_space} and \cref{lem:rank:bound:item1} in \cref{lem:rank:bound} \hence establish that there exist $\lambda_i \colon \cM \allowbreak \to \R$, $i \allowbreak \in \{1, \allowbreak 2, \allowbreak \dots, \allowbreak \dimension \}$, and $\mf v_i \colon \cM \allowbreak \to \R^{\dimension}$, $i \allowbreak \in \{1, \allowbreak 2, \allowbreak \dots, \allowbreak \dimension \}$, which satisfy that
\begin{enumerate}[label=(\roman*)]
\item\label{item:lem:eigenvalues:locmin:point:1} it holds for all $ y \in \cM $ that $ \{ \mf v_{\dimension - k + 1}(y), \mf v_{\dimension - k + 2}(y), \dots, \mf v_{\dimension}(y) \} $ is a Hamel basis of $ \tangentspace{y}{\cM} $,

\item\label{item:lem:eigenvalues:locmin:point:2} it holds for all $ y \in \cM $, $ i \in \{ 1, 2, \dots, \dimension \} $ that $ ( \Hs f )( y ) \mf v_i(y) = \lambda_i(y) \mf v_i(y) $, and
    
\item\label{item:lem:eigenvalues:locmin:point:3} it holds for all $ y \in \cM $, $ i, j \in \{ 1, 2, \dots, \dimension \} $ that
\begin{equation}
\scalprod{ \mf v_i(y) }{ \mf v_j(y) } =
\begin{cases}
1 &\colon i = j \\
0 &\colon i \neq j.
\end{cases}
\end{equation}
\end{enumerate}
\Nobs that \cref{item:lem:eigenvalues:locmin:point:1,item:lem:eigenvalues:locmin:point:2,item:lem:eigenvalues:locmin:point:3}, the assumption that $ \rk((\Hs f)(x)) \allowbreak = \dimension \allowbreak - k $, the assumption that $\spectrum{(\Hs f)(x)}\allowbreak \subseteq [0, \allowbreak \infty)$, and \cref{lem:rank:bound:item1} in \cref{lem:rank:bound} show that $\min_{i \in \{1, 2, \dots, \dimension - k \}} \allowbreak \lambda_i(x) \allowbreak > 0$. Moreover, \nobs that, e.g., \cite[Proposition~34]{Adrian2021ConvergencePiecewise} ensures that $\cM \allowbreak \subseteq \bfPP_{\cM}$. This, the fact that $\bfPP_{\cM} \allowbreak \subseteq \R^\dimension$ is open, the fact that $\min_{ i \in \{ 1, 2, \dots, \dimension - k \} } \lambda_i(x) \allowbreak > 0 $, the continuity of the tangent spaces of $\cM$, and the fact that the eigenvalues depend continuously on a symmetric matrix (cf., e.g., Kato~\cite[Theorem 5.2]{Kato1995_1980}) demonstrate that there exist $ \varepsilon \in (0, \infty) $ and an open set $ W \subseteq \bfPP_{\cM} $ which satisfy that
\begin{enumerate}[label=(\alph*)]
\item it holds that $ x \in W $,
\item it holds for all $ y \in \cM \cap W $ that $ \min_{i \in \{1, 2, \dots, \dimension - k\}} \lambda_i(y) > \varepsilon $,
\item it holds that $ \cM \cap W $ is connected,
\item it holds that $ \overline{W} $ is compact, and
\item it holds that $ \overline{W} \subseteq U$ \ifnocf.
\end{enumerate}
\cfload[.]\Nobs that \cref{item:lem:eigenvalues:locmin:point:1,item:lem:eigenvalues:locmin:point:3} ensure that for all $ y \in \cM \cap W $, $ v \in (\tangentspace{y}{\cM})^\perp $ there exist $ u_1, \dots, u_{\dimension - k} \in \R $ such that $ v = \sum_{i=1}^{\dimension - k} u_i \mf v_i(y) $. This implies that for all $ y \in \cM \cap W $, $ v \in (\tangentspace{y}{\cM})^\perp $ there exist $ u_1, \dots, u_{\dimension - k} \in \R $ such that
\begin{equation} \label{lem:critical:proof:eq0}
\begin{split}
\scalprod{ v }{ (\Hs f)(y) v } & \textstyle = \sum_{i=1}^{\dimension - k} \sum_{j=1}^{\dimension - k} u_i u_j \lambda_j(y) \scalprod{ \mf v_i(y) }{ \mf v_j(y) } = \sum_{i=1}^{\dimension - k} ( \lambda_i(y) \abs{u_i}^2 ) \\
& \textstyle \geq \varepsilon \sum_{i=1}^{\dimension - k} \abs{u_i}^2 = \varepsilon \sum_{i=1}^{\dimension - k} \sum_{j=1}^{\dimension - k} u_i u_j \scalprod{ \mf v_i(y) }{ \mf v_j(y) } = \varepsilon \norm{v}^2.
\end{split}
\end{equation}
Next \nobs that \cref{lem:mult:taylor} (applied with $K \with \overline{W}$, $\varepsilon \with \varepsilon/4$ in the notation of \cref{lem:mult:taylor}) proves that there exists $\eta \in (0,\infty)$ which satisfies for all $y,z \in \overline{W}$ with $\norm{y-z}<\eta$ and $\cup_{t\in[0,1]}\{(1-t)y+tz\}\subseteq\overline{W}$ that
\begin{equation} \label{lem:critical:proof:eq1}
\textstyle \abs{f(y)-f(z)-\scalprod{(\nabla f)(z)}{y-z}-\frac{1}{2}\scalprod{y-z}{(\Hs f)(z)(y-z)}}\leq\frac{\varepsilon}{4}\norm{y-z}^2.
\end{equation}
Furthermore, \nobs that the assumption that $ \cM \cap W $ is connected and the fact that for all $ y \in \cM $ it holds that $ (\nabla f)(y) = 0 $ show that for all $ y \in \cM \cap W $ it holds that $ f(y) = f(x) $. Next \nobs that \cite[Proposition~39]{Adrian2021ConvergencePiecewise} ensures that there exist $ R, \delta \in (0, \infty) $ which satisfy for all $ y \in \cM \cap \{w \in \R^{\dimension} \colon \norm{x - w} \le R\} $, $ v \in (\tangentspace{y}{\cM})^\perp \cap \{w \in \R^{\dimension} \colon \norm{w} < \delta\} $ that $ \pp_{\cM} (y + v) = y $, $\delta<\eta$, $ \tubular{}{R, \delta} \subseteq W $, and
\begin{equation} \label{lem:critical:proof:eq2}
\textstyle \tubular{}{R, \delta} = \{ a \in \R^{\dimension} \colon \inf\nolimits_{b \in \cM} \norm{a - b} = \inf\nolimits_{b \in \cM \cap \{w \in \R^{\dimension} \colon \norm{x - w} \le R\}} \norm{a - b} < \delta \}\ifnocf.
\end{equation}
\cfload[.]\Nobs that \cref{lem:critical:proof:eq2} implies that for all $ y \in \tubular{}{R, \delta} $, $ t \in [0,1] $ it holds that $ y - \pp_{\cM}(y) \in (\tangentspace{\pp_{\cM}(y)}{\cM})^\perp $ and $ \pp_{\cM}(y) + t (y - \pp_{\cM}(y)) \in \tubular{}{R, \delta} $. Combining this with \cref{lem:critical:proof:eq0}, \cref{lem:critical:proof:eq1}, the fact that for all $ y \in \tubular{}{R, \delta} $ it holds that $ (\nabla f)(\pp_{\cM}(y)) = 0 $, and the fact that for all $ y \in \tubular{}{R, \delta} $ it holds that $ f(\pp_{\cM}(y)) = f(x) $ establishes that for all $ y \in \tubular{}{R, \delta} $ it holds that
\begin{equation}\label{eqn:lem:eigenvalues:locmin:point:final}
\begin{split}
\textstyle f(y) & \geq f( \pp_{\cM}(y)) + \spro{(\nabla f)(\pp_{\cM}(y)), y - \pp_{\cM}(y)} \\
& \textstyle \quad+ \frac{1}{2} \spro{y - \pp_{\cM}(y), (\Hs f)(\pp_{\cM}(y))(y - \pp_{\cM}(y))}  - \frac{\varepsilon}{4} \norm{y - \pp_{\cM}(y)}^2 \\
& \textstyle \geq f(x) + \frac{\varepsilon}{2} \norm{y - \pp_{\cM}(y)}^2 - \frac{\varepsilon}{4} \norm{y - \pp_{\cM}(y)}^2 \\
& \textstyle \geq f(x) + \frac{\varepsilon}{4}\norm{y - \pp_{\cM}(y)}^2 \geq  f(x).
\end{split}
\end{equation}
Moreover, \nobs that the continuity of $\pp_{\cM}$ on $W$, the fact that $\pp_{\cM}(x)=x$, and \cite[Item~$($iv$)$ in Proposition~39]{Adrian2021ConvergencePiecewise} ensure that there exists $\rho\in(0,\infty)$ which satisfies $\{y\in\R^{\dimension}\colon\norm{y-x}<\rho\}\subseteq\tubular{}{R,\delta}$. This and \cref{eqn:lem:eigenvalues:locmin:point:final} prove in the case $k<\dimension$ that $x$ is a \locmin point of $f$\cfload.
\end{cproof}

\cfclear
\begin{lemma}[\resname{Necessary condition for local minima}]\label{lem:eigenvalues:nonnegative}
Let $\dimension \in \N$, let $U \subseteq \R^{\dimension}$ be open, let $f \allowbreak \in \allowbreak C^2(U, \allowbreak \R)$, $x \allowbreak \in \allowbreak U$, and assume that $x$ is a \locmin point of $f$ \cfload. Then
\begin{equation}\label{eq:lem:eigenvalues:nonnegative}
\spectrum{(\Hs  f)(x)} \subseteq [0, \infty)
\end{equation}
\cfout.
\end{lemma}
\begin{cproof}{lem:eigenvalues:nonnegative}
We prove \cref{eq:lem:eigenvalues:nonnegative} by contradiction. In the following we thus assume that there exist $ \lambda \in (-\infty,0) $, $ v \in \R^{ \dimension } \backslash \{ 0 \} $ which satisfy
\begin{equation}\label{eq:lem:eigenvalues:nonnegative:Hessian}
(\Hs f)(x) v = \lambda v.
\end{equation}
\Nobs that the assumption that $U \subseteq \R^{\dimension}$ is open and the assumption that $f \in C^2(U, \R)$ ensure that there exist $\varepsilon \in (0,\infty)$, $g \in C^2( (-\varepsilon,\varepsilon), \R )$ which satisfy for all $t \in (-\varepsilon,\varepsilon)$ that
\begin{equation}\label{eq:lem:eigenvalues:nonnegative:g}
x + tv \in U \qqandqq g(t) = f( x + tv ).
\end{equation}
\Nobs that \cref{eq:lem:eigenvalues:nonnegative:g}, the assumption that $ f \in C^2(U, \R) $, and the chain rule show that for all $t \in (-\varepsilon,\varepsilon) $ it holds that
\begin{equation}\label{eq:lem:eigenvalues:nonnegative:g:prime}
g'(t) = \scalprod{ v }{ (\nabla f)( x + tv ) } \qqandqq g''(t) = \scalprod{ v }{ (\Hs f)( x + tv ) v }\ifnocf.
\end{equation}
\cfload[.]Moreover, \nobs that the assumption that $ x $ is a \locmin point of $ f $ demonstrates that $ (\nabla f) (x) = 0 $. Combining this, \cref{eq:lem:eigenvalues:nonnegative:Hessian}, and \cref{eq:lem:eigenvalues:nonnegative:g:prime} establishes that
\begin{equation}\label{eq:lem:eigenvalues:nonnegative:g:prime:0}
g'(0) = \scalprod{ v }{ (\nabla f) (x) } = 0 \qqandqq g''(0) = \scalprod{ v }{ (\Hs f) (x) v } = \scalprod{ v }{ \lambda v } = \lambda \norm{ v }^2 < 0.
\end{equation}
\Hence that there exists $ \delta \in (0,\varepsilon) $ which satisfies for all $ t \in (-\delta,\delta) $ that
\begin{equation}\label{eq:lem:eigenvalues:nonnegative:g:prime:negative}
g''(t) < \nicefrac{g''(0)}{2}.
\end{equation}
\Nobs that \cref{eq:lem:eigenvalues:nonnegative:g:prime:0}, \cref{eq:lem:eigenvalues:nonnegative:g:prime:negative}, and the fundamental theorem of calculus establish that for all $ t \in (-\delta,\delta) \backslash \{ 0 \} $ it holds that
\begin{equation}
\begin{split}
g(t) & = \textstyle g(0) + \int_0^t g'(s) \, \d s = g(0) + \int_0^t \bigl( g'(0) + \int_0^s g''(r) \, \d r \bigr) \,\d s \\
& = \textstyle g(0) + \int_0^t \int_0^s g''(r) \, \d r \, \d s < g(0).
\end{split}
\end{equation}
\Hence for all $t \in (-\delta,\delta) \backslash \{ 0 \}$ that $f( x + tv) = g(t) < g(0) = f(x)$. This is a contradiction to the assumption that $ x $ is a \locmin point of $f$.
\end{cproof}

\cfclear
\begin{proposition}[\resname{Classification of a critical point}]\label{prop:class:crit:points}
Let $\dimension, \allowbreak k \allowbreak \in \N$, let $U \allowbreak \subseteq \R^{\dimension}$ be open, let $f \allowbreak \in \allowbreak C^2(U, \allowbreak \R)$, let $\cM \allowbreak \subseteq U$ satisfy $\cM \allowbreak \subseteq \allowbreak \{x \allowbreak \in \allowbreak U \colon \allowbreak (\nabla f)(x) \allowbreak = 0\}$, assume that $\cM$ is a $k$-dimensional $C^2 $-\submanifold of $\R^{\dimension}$, and let $x \in \cM$ satisfy $\rk((\Hs f)(x)) \geq \dimension - k$ \cfload. Then
\begin{enumerate}[label=(\roman*)]
\item \label{prop:class:crit:points:item1} it holds that $ \rk((\Hs f)(x)) = \dimension - k $,
\item \label{prop:class:crit:points:item2} it holds that $ x $ is a \locmin point of $ f $ if and only if $ \spectrum{(\Hs  f)(x)} \subseteq [0, \infty) $,
\item \label{prop:class:crit:points:item3} it holds that $ x $ is a \locmax point of $ f $ if and only if $ \spectrum{(\Hs  f)(x)} \subseteq (- \infty, 0] $, and
\item \label{prop:class:crit:points:item4} it holds that $ x $ is a \saddle point of $ f $ if and only if
\begin{equation}
\min \{\# (\spectrum{(\Hs  f)(x)} \cap (0, \infty) ), \# (\spectrum{(\Hs  f)(x)} \cap ( - \infty, 0) )\} > 0 \ifnocf.
\end{equation}
\end{enumerate}
\cfout[.]
\end{proposition}
\begin{cproof}{prop:class:crit:points}
\Nobs that \cref{lem:rank:bound} and the assumption that
\begin{equation}
\textstyle \rk((\Hs f)(x)) \allowbreak \ge \dimension - k
\end{equation}
demonstrate that $\rk((\Hs f)(x)) = \dimension - k $. This establishes \cref{prop:class:crit:points:item1}. \Nobs that \cref{prop:class:crit:points:item1} and \cref{lem:eigenvalues:locmin:point} prove that
\begin{equation}\label{eq:prop:class:crit:points:implication1}
\bigl( \bigl[ \spectrum{(\Hs f)(x)} \subseteq [0, \infty) \bigr] \rightarrow \bigl[ x \text{ is a \locmin point of } f \bigr] \bigr).
\end{equation}
Moreover, \nobs that \cref{lem:eigenvalues:nonnegative} establishes that
\begin{equation}\label{eq:prop:class:crit:points:implication2}
    \bigl( \bigl[ x \text{ is a \locmin point of } f \bigr]
    \rightarrow
    \bigl[ \spectrum{(\Hs f)(x)} \subseteq [0, \infty) \bigr] \bigr).
  \end{equation}
  Combining this and \cref{eq:prop:class:crit:points:implication1} establishes \cref{prop:class:crit:points:item2}.
  \Nobs that \cref{prop:class:crit:points:item1} and \cref{lem:eigenvalues:locmin:point}
  (applied with $ f \with - f $ in the notation of \cref{lem:eigenvalues:locmin:point})
  ensure that
  \begin{equation}\label{eq:prop:class:crit:points:implication3}
    \bigl( \bigl[ \spectrum{(\Hs f)(x)} \subseteq (-\infty, 0] \bigr]
    \rightarrow
    \bigl[ x \text{ is a \locmax point of } f \bigr] \bigr).
  \end{equation}
  In addition, \nobs that \cref{lem:eigenvalues:nonnegative}
  (applied with $ f \with - f $ in the notation of \cref{lem:eigenvalues:nonnegative})
  demonstrates that
  \begin{equation}\label{eq:prop:class:crit:points:implication4}
    \bigl( \bigl[ x \text{ is a \locmax point of } f \bigr]
    \rightarrow
    \bigl[ \spectrum{(\Hs f)(x)} \subseteq (-\infty, 0] \bigr] \bigr).
  \end{equation}
  This and \cref{eq:prop:class:crit:points:implication3} establish \cref{prop:class:crit:points:item3}.
  \Nobs that \cref{prop:class:crit:points:item2,prop:class:crit:points:item3} prove \cref{prop:class:crit:points:item4}.
\end{cproof}

\section{On infinitely many realization functions of non-global local minimum points}
\label{sec:negative_result}

In this section we employ \cref{prop:class:crit:points} from \cref{sec:differential_geometric_prelimiaries} above to establish in \cref{cor:uncountably_many_realizations_of_non_global_local_minimas_with_varying_delta} in \cref{subsec:existence_of_inf_many_realizations} below that there exist uncountably many Lipschitz continuous target functions such that the associated risk functions have infinitely many realization functions of non-global local minimum points. \cref{cor:uncountably_many_realizations_of_non_global_local_minimas_with_varying_delta} is a simple consequence of \cref{cor:uncountably_many_realizations_of_non_global_local_minimas_with_fixed_delta} in \cref{subsec:existence_of_inf_many_realizations}. Our proof of \cref{cor:uncountably_many_realizations_of_non_global_local_minimas_with_fixed_delta}, in turn, makes use of \cref{lemma:zero_derivative} in \cref{subsec:risks_for_suitable_ANN_realization_functions} below, \cref{lemma:rank&Hessian} in \cref{subsec:rank_of_hessian_of_risk_function} below, and \cref{lemma:uncountably_many_realizations_of_critical_points} in \cref{subsec:submanifold_of_local_minimum_points_of_ANN_parameter_space} below.
\begin{itemize}
\item In \cref{lemma:zero_derivative} we calculate the risks for suitable \ANN\ realization functions,

\item in \cref{lemma:rank&Hessian} we establish some properties of the Hessian matrix of the risk function at points in $\cM \subseteq \R^{\dimension}$, and

\item in \cref{lemma:uncountably_many_realizations_of_critical_points} we prove that every $\theta \in \cM$ has the same risk value and is a critical point of $\riskR \colon \R^{ \dimension } \to \R$.
\end{itemize}
Our proof of \cref{lemma:rank&Hessian} uses a well-known criterion for positive definiteness of symmetric matrices -- known as the Schur complement criterion in the literature -- which we recall in \cref{lemma:Schur_complement} in \cref{subsec:rank_of_hessian_of_risk_function} below. Only for the sake of completeness we include the proof of \cref{lemma:Schur_complement} in this section.

In \cref{setting0} in \cref{subsec:ANNs_onedim_input_and_multdim_hidden_layer} below we introduce the mathematical objects considered in this section such as the specific target function $f \colon [\scra, \allowbreak \scrb]^d \allowbreak \to \R$ (see \cref{0_eq:f}), the realization functions $\functionSNN^{\theta} \colon [\scra, \allowbreak \scrb]^d \allowbreak \to \R$, $\theta \in \R^{\dimension}$, the risk function $\riskR \colon \R^{\dimension} \to \R$, the $(\dimension - d - 1)$-dimensional $C^{\infty}$-submanifold $\cM$ of $\R^{\dimension}$ (see~\cref{0_eq:M} and \cref{lemma:set_of_critical_points_submanifold}), the suitable open set $\fV \subseteq \R^{\dimension}$ which contains $\cM$ and on which the risk function $\riskR$ is twice continuously differentiable (see~\cref{0_eq:fV} and ~\cref{lemma:Hessian_of_Risk_cM}), and the realization functions $\bfM \subseteq C([\scra, \scrb]^d, \R)$ associated to $\cM \allowbreak \subseteq \R^{\dimension}$. Some of the calculations in this section were performed using {\sc Wolfram Mathematica} (see~\cite{Mathematica}).

In \cref{fig:local:minimum:points} in \cref{subsec:existence_of_inf_many_realizations} we present numerical simulations associated to \cref{cor:uncountably_many_realizations_of_non_global_local_minimas_with_fixed_delta} in the case where $d = 1$, $\width = 4$, $\dimension = 13$, $\scra = 0$, $\scrb = 1$, $\alpha = \nicefrac{1}{3}$, and $\beta = \nicefrac{2}{3}$. In these simulations we randomly initialize 50 \ANNs\ with the Xavier initialization, then we approximately train these \ANNs\ with the \GD\ optimization method using a learning rate of $\nicefrac{1}{20}$ until the maximum norm of the generalized gradient function evaluated at the current position of the \GD\ process is strictly less than $10^{- 4}$, and, thereafter, we gradually plot the realization functions of the resulting \ANNs\ whereby a realization function is not drawn if a realization function with an $L^2$-distance strictly less than $10^{- 2}$ has already been drawn. We also refer to \cref{list:python} for the {\sc Python} source code used to create \cref{fig:local:minimum:points}.

For intuition about how we found the specific target function in \cref{0_eq:f} in \cref{setting0} resulting in an objective function with infinitely many different realization functions of non-global local minimum points we refer the reader to \cref{remark:intuition} in \cref{subsec:main_positive_results} below.

\subsection{ANNs and special target functions}
\label{subsec:ANNs_onedim_input_and_multdim_hidden_layer}

\cfclear
\begin{setting}\label{setting0}
Let $d, \width, \dimension \in \N$, $\alpha, \scra, \beta \in \R$, $\scrb \in (\scra, \infty)$ satisfy $\dimension = d \width + 2 \width + 1$ and $0 < \alpha \allowbreak < \beta \allowbreak < 1$, let $e_1, \allowbreak e_2, \allowbreak \dots, \allowbreak e_d \in \R^d$ satisfy $e_1 \allowbreak = (1, \allowbreak 0, \allowbreak 0, \allowbreak \dots, \allowbreak 0)$, $e_2 \allowbreak = (0, \allowbreak 1, \allowbreak 0, \allowbreak \dots, \allowbreak 0)$, $\dots$, $e_d \allowbreak = (0, \allowbreak 0, \allowbreak \dots, \allowbreak 0, \allowbreak 1) \allowbreak \in \R^d$, for every $n \in \N$ let $\ones{n} \in \R^n$ satisfy $\ones{n} = (1, 1, \dots, 1)$, let $\ff \allowbreak \in \allowbreak C([0, \allowbreak 1]^d, \allowbreak \R)$, $f \allowbreak \allowbreak \in \allowbreak C([\scra, \allowbreak \scrb]^d, \allowbreak \R)$ satisfy for all $x = \allowbreak (x_1, \allowbreak \dots, \allowbreak x_d) \allowbreak \in [0, \allowbreak 1]^d$ that
\begin{equation}\label{0_eq:f}
\ff(x) =
\begin{cases}
\frac{4(x_1 - \alpha) + 3 \alpha^2 - 1}{4 (1 - \alpha)^{1/2} (1 + 3 \alpha)^{3/2}} & \colon x_1 \in [0, \alpha] \\[1ex]
\frac{3 [x_1]^2 - 1}{4 (1 - x_1)^{1/2} (1 + 3 x_1)^{3/2}} & \colon x_1 \in (\alpha, \beta] \\[1ex]
\frac{12 \beta [x_1]^2 - (18 \beta^2 + 8 \beta - 2) x_1 + 3 \beta^4 + 10 \beta^2 - 1}{4 (1 - \beta)^{5/2} (1 + 3 \beta)^{3/2}} & \colon x_1 \in (\beta, 1]
\end{cases}
\end{equation}
and $f((\scrb - \scra) x + \scra \ones{d}) = \ff(x)$, let $\Rect_r \allowbreak \colon \R \allowbreak \to \R$, $r \in \N \cup \{\infty\}$, satisfy for all $x \in \R$ that $(\cup_{r \in \N} \{\Rect_r \}) \allowbreak \subseteq C^1 (\R, \allowbreak \R)$, $\Rect_\infty(x) = \max \{x, 0\}$, $\sup_{r \in \N} \sup_{y \in [- \abs{x}, \abs{x}]} \abs{(\Rect_r)'(y)} < \infty$, and
\begin{equation}
\limsup\nolimits_{r \to \infty} (\abs{\Rect_r(x) - \Rect_\infty(x)} + \abs{(\Rect_r)'(x) - \indicator{(0, \infty)}(x)}) = 0,
\end{equation}
for every $\theta = \allowbreak (\theta_1, \allowbreak \dots, \allowbreak \theta_{\dimension}) \allowbreak \in \R^{\dimension}$, $i \in \allowbreak \{1, \allowbreak 2, \allowbreak \dots, \allowbreak \width\}$ let $\fw_i^{\theta} \allowbreak = (\fw_{i, 1}^{\theta}, \allowbreak \dots, \allowbreak \fw_{i, d}^{\theta}) \allowbreak \in \R^d$, $\fb_i^{\theta}, \fv_i^{\theta}, \fc^{\theta} \in \R$ satisfy $\wedge_{j = 1}^d [\fw_{i, j}^{\theta} \allowbreak = \theta_{(i - 1)d + j}]$, $\fb_i^{\theta} \allowbreak = \theta_{d\width + i}$, $\fv_i^{\theta} \allowbreak = \theta_{d\width + \width + i}$, and $\fc^{\theta} \allowbreak = \theta_{\dimension}$ and let $I_i^\theta \allowbreak \subseteq \R^d$ and $\partial I_i^\theta \allowbreak \subseteq \R^d$ satisfy $I_i^\theta = \allowbreak \{x \allowbreak \in [\scra, \scrb]^d \colon \allowbreak \fb_i^{\theta} + \allowbreak \spro{\fw_i^{\theta}, x} \allowbreak > 0\}$ and $\partial I_i^\theta = \allowbreak \{x \allowbreak \in [\scra, \allowbreak \scrb]^d \colon \allowbreak \fb_i^{\theta} + \allowbreak \spro{\fw_i^{\theta}, x} \allowbreak = 0\}$, for every $r \in \allowbreak \N \allowbreak \cup \{\infty\}$ let $\cN_r \allowbreak = (\realapprox{\theta}{r})_{\theta \in \R^{\dimension}} \colon \R^{\dimension} \allowbreak \to C([\scra, \allowbreak \scrb]^d, \allowbreak \R)$ and $\cR_r \colon \R^{\dimension} \to \R$ satisfy for all $\theta \allowbreak \in \R^{\dimension}$, $x \allowbreak \in [\scra, \allowbreak \scrb]^d$ that 
\begin{equation}\label{setting0:eq:realization_risk}
\textstyle \realapprox{\theta}{r}(x) = \fc^{\theta} + \smallsum_{i = 1}^\width \fv_i^{\theta} \bigl[\Rect_r (\fb_i^{\theta} + \spro{\fw_i^{\theta}, x})\bigr] \qandq \cR_r(\theta) = \int_{[\scra, \scrb]^d} \abs{\realapprox{\theta}{r}(y) - f(y)}^2 \, \d y,
\end{equation}
let $\pi \colon \N \allowbreak \to \N$ satisfy $\pi(1) = 1$, $\pi(2) \allowbreak = d \width \allowbreak + 1$, $\pi(3) \allowbreak = d \width \allowbreak + \width \allowbreak + 1$, $\pi(4) \allowbreak = \dimension$, and $\wedge_{k = 5}^{d + 3} [\pi(k) = k - 3]$, let $\cM \subseteq \allowbreak \R^{\dimension}$, $\fV \subseteq \allowbreak \R^{\dimension}$, and $\bfM \subseteq C([\scra, \scrb]^d, \R)$ satisfy
\begin{multline}\label{0_eq:M}
\textstyle \cM = \Bigl\{\theta \in \R^{\dimension} \colon \textstyle \bigl(\Exists (\fx, \fy) \in (\alpha, \beta) \times (0, \infty) \colon \bigl[\bigl(\sum_{j = 2}^d \abs{\fw_{1, j}^{\theta}} = - \abs{\fw_{1, 1}^{\theta} - \frac{\fy}{\scrb - \scra}}\bigr) \\
\textstyle \wedge\! \bigl(\frac{\fb_1^{\theta}}{\fy} = \frac{\scra}{\scra - \scrb} - \fx\bigr) \!\wedge\! \bigl(\fv_1^{\theta} = \frac{(1 + 3 \fx)^{- 1/2}}{2 \fy (1 - \fx)^{3/2}}\bigr)\! \wedge\! \bigl(\fc^{\theta} = \frac{- (1 - \fx)^{1/2}}{4 (1 + 3 \fx)^{1/2}}\bigr) \!\wedge \!\bigl(\bigwedge_{i = 2}^\width [I_i^{\theta} = \emptyset = \partial I_i^{\theta}]\bigr)\bigr]\bigr)\!\Bigr\},
\end{multline}
\begin{multline}\label{0_eq:fV}
\textstyle \fV = \Bigl\{\theta \in \R^{\dimension} \colon \bigl(\Exists (\fx, \fy) \in (\alpha, \beta) \times (0, \infty) \colon \bigl[(\fw_{1, 1}^{\theta} = \frac{\fy}{\scrb - \scra}) \\
\textstyle \wedge\! \bigl(\sum_{j = 2}^d \abs{\fw_{1, j}^{\theta}} < \frac{\fy \min\{\fx - \alpha, \beta - \fx\}}{\abs{\scrb} + \abs{\scra}}\bigr) \!\wedge\! \bigl(\frac{\fb_1^{\theta}}{\fy} = \frac{\scra}{\scra - \scrb} - \fx\bigr) \!\wedge\! \bigl(\bigwedge_{i = 2}^\width [I_i^{\theta} = \emptyset = \partial I_i^{\theta}]\bigr)\bigr]\bigr)\!\Bigr\},
\end{multline}
and $\bfM = \cup_{\theta \in \cM} \{\functionSNN^{\theta}\}$, and for every $\theta \allowbreak \in \fV$ let $\scrq_{\theta}, \fq_{\theta} \in \R$ satisfy $\scrq_{\theta} \allowbreak = \allowbreak - \frac{\fb_1^{\theta}}{\fw_{1, 1}^{\theta}}$ and $\fq_{\theta} = \frac{\scrq_{\theta} - \scra}{\scrb - \scra}$.
\end{setting}

\subsection{Risks for suitable ANN realization functions}
\label{subsec:risks_for_suitable_ANN_realization_functions}

\cfclear
\begin{lemma}\label{lemma:zero_derivative}
Assume \cref{setting0} and let $q \in (\alpha, \beta)$, $g, \fg \in C([0, 1], \R)$ satisfy for all $x \in [0, 1]$ that $\fg(x) = \ff(x, 0, \dots, 0)$ and
\begin{equation}
\textstyle g(x) = - \frac{(1 - q)^{1/2}}{4(1 + 3q)^{1/2}} + \frac{\max\{x - q, 0\}}{2(1 - q)^{3/2} (1 + 3q)^{1/2}} = \frac{2 \max\{x - q, 0\} - (1 - q)^2}{4(1-q)^{3/2}(1+3q)^{1/2}}.
\end{equation}
Then
\begin{enumerate}[label=(\roman *)]
\item
\label{item1:lemma:zero_derivative} it holds that 
\begin{equation}
\textstyle \int_0^q (g(x) - \fg(x)) \, \d x = \int_q^1 (g(x) - \fg(x)) \, \d x = \int_q^1 x (g(x) - \fg(x)) \, \d x = 0
\end{equation}

and
\item
\label{item2:lemma:zero_derivative} it holds that
\begin{equation}
\textstyle \int_0^1 (g(x) - \fg(x))^2 \, \d x = \int_0^1 [\fg(x)]^2 \, \d x - \frac{1}{48}.
\end{equation}
\end{enumerate}
\end{lemma}
\begin{cproof}{lemma:zero_derivative}
\Nobs that the chain rule ensures that for all $x \in (0, 1)$ it holds that
\begin{equation}\label{eqn:lemma:zero_derivative:anti_derivative_f}
\begin{split}
\textstyle \bigl[- \frac{x (1 - x)^{1/2}}{4 (1 + 3x)^{1/2}}\bigr]' & \textstyle = - \frac{[(1 - x)^{1/2} - (x/2) (1 - x)^{-1/2}] (1 + 3x)^{1/2} - x (1 - x)^{1/2} [(3/2)(1 + 3x)^{-1/2}]}{4 (1 + 3x)} \\
& \textstyle = - \frac{(2 (1 - x) - x)(1 + 3x) - 3 x (1 - x)}{8 (1 - x)^{1/2} (1 + 3x)^{3/2}} = \frac{3x^2 - 1}{4 (1 - x)^{1/2} (1 + 3x)^{3/2}} 
\end{split}
\end{equation}
and
\begin{equation}\label{eqn:lemma:zero_derivative:anti_derivative_xf}
\begin{split}
& \textstyle \bigl[-\frac{(3x^2 + 2x + 1) (1 - x)^{1/2}}{24 (1 + 3x)^{1/2}}\bigr]' \\
& \textstyle = - \frac{[(6x + 2)(1 - x)^{1/2} - (3x^2 + 2x + 1)(1/2)(1 - x)^{-1/2}] (1 + 3x)^{1/2} - (3x^2 + 2x + 1)(1 - x)^{1/2} [(3/2)(1 + 3x)^{-1/2}]}{24 (1 + 3x)} \\
& \textstyle = - \frac{[(6x + 2)(2 - 2x) - (3x^2 + 2x + 1)] (1 + 3x) - (3x^2 + 2x + 1)(3 - 3x)}{48 (1 - x)^{1/2} (1 + 3x)^{3/2}} = - \frac{(3x + 1)^2 (1 - x) - (3x^2 + 2x + 1)}{12 (1 - x)^{1/2} (1 + 3x)^{3/2}} \\
& \textstyle = \frac{x (3x^2 - 1)}{4 (1 - x)^{1/2} (1 + 3x)^{3/2}}.
\end{split}
\end{equation}
\Hence that
\begin{equation}
\begin{split}
& \textstyle \int_0^q (g(x) - \fg(x)) \, \d x \\
& \textstyle = \int_0^{\alpha} (g(x) - \fg(x)) \, \d x + \int_{\alpha}^q (g(x) - \fg(x)) \, \d x \\
& \textstyle = \int_0^{\alpha} \bigl(- \frac{(1-q)^{1/2}}{4 (1 + 3q)^{1/2}} - \frac{4(x - \alpha) + 3 \alpha^2 - 1}{4 (1 - \alpha)^{1/2} (1 + 3 \alpha)^{3/2}}\bigr) \, \d x + \int_{\alpha}^q \bigl(- \frac{(1-q)^{1/2}}{4 (1 + 3q)^{1/2}} - \frac{3 x^2 - 1}{4 (1 - x)^{1/2} (1 + 3x)^{3/2}}\bigr) \, \d x \\
& = \textstyle \bigl[- \frac{(1 - q)^{1/2} x}{4 (1 + 3 q)^{1/2}} - \frac{2x^2 + (3 \alpha^2 - 4 \alpha - 1) x}{4 (1 - \alpha)^{1/2}(1 + 3 \alpha)^{3/2}}\bigr]_{x = 0}^{x = \alpha} + \bigl[- \frac{(1 - q)^{1/2} x}{4 (1 + 3q)^{1/2}} + \frac{x (1 - x)^{1/2}}{4 (1 + 3x)^{1/2}}\bigr]_{x = \alpha}^{x = q} \\
& \textstyle = \bigl[- \frac{(1 - q)^{1/2} \alpha}{4 (1 + 3 q)^{1/2}} - \frac{3 \alpha^3 - 2 \alpha^2 - \alpha}{4 (1 - \alpha)^{1/2}(1 + 3 \alpha)^{3/2}}\bigr] + \bigl[\frac{(1 - q)^{1/2} \alpha}{4 (1 + 3q)^{1/2}} - \frac{\alpha (1 - \alpha)^{1/2}}{4 (1 + 3\alpha)^{1/2}}\bigr] \\
& \textstyle = - \frac{3 \alpha^3 - 2 \alpha^2 - \alpha}{4 (1 - \alpha)^{1/2}(1 + 3 \alpha)^{3/2}} - \frac{\alpha (1 - \alpha)^{1/2}}{4 (1 + 3\alpha)^{1/2}} = \frac{\alpha (1 + 3 \alpha) (1 - \alpha)}{4 (1 - \alpha)^{1/2}(1 + 3 \alpha)^{3/2}} - \frac{\alpha (1 - \alpha)^{1/2}}{4 (1 + 3\alpha)^{1/2}} = 0, 
\end{split}
\end{equation}
\begin{equation}
\begin{split}
& \textstyle \int_q^1 (g(x) - \fg(x)) \, \d x \\
& \textstyle = \int_q^{\beta} (g(x) - \fg(x)) \, \d x + \int_{\beta}^1 (g(x) - \fg(x)) \, \d x \\
& \textstyle = \int_q^{\beta} \bigl(\frac{2 x - 1 - q^2}{4(1 - q)^{3/2}(1 + 3q)^{1/2}} - \frac{3 x^2 - 1}{4 (1 - x)^{1/2} (1 + 3x)^{3/2}}\bigr) \, \d x \\
& \quad \textstyle + \int_{\beta}^1 \bigl(\frac{2 x - 1 - q^2}{4(1 - q)^{3/2}(1 + 3q)^{1/2}} - \frac{12 \beta x^2 - (18 \beta^2 + 8 \beta - 2) x + 3 \beta^4 + 10 \beta^2 - 1}{4 (1 - \beta)^{5/2} (1 + 3 \beta)^{3/2}}\bigr) \, \d x \\
& = \textstyle \bigl[\frac{x^2 - (1 + q^2) x}{4 (1 - q)^{3/2}(1 + 3q)^{1/2}} + \frac{x(1-x)^{1/2}}{4 (1 + 3x)^{1/2}}\bigr]_{x = q}^{x = \beta} \\
& \quad \textstyle + \bigl[\frac{x^2 - (1 + q^2) x}{4 (1 - q)^{3/2}(1 + 3q)^{1/2}} - \frac{4 \beta x^3 - (9 \beta^2 + 4 \beta - 1) x^2 + (3 \beta^4 + 10 \beta^2 - 1) x}{4 (1 - \beta)^{5/2}(1 + 3 \beta)^{3/2}}\bigr]_{x = \beta}^{x = 1} \\
& = \textstyle \bigl[\frac{\beta^2 - (1 + q^2) \beta}{4 (1 - q)^{3/2}(1 + 3q)^{1/2}} + \frac{\beta (1-\beta)^{1/2}}{4 (1 + 3 \beta)^{1/2}} - \frac{q^2 - (1 + q^2) q}{4 (1 - q)^{3/2}(1 + 3q)^{1/2}} - \frac{q (1-q)^{1/2}}{4 (1 + 3 q)^{1/2}}\bigr] \\
& \quad \textstyle + \bigl[\frac{- q^2}{4 (1 - q)^{3/2}(1 + 3q)^{1/2}} - \frac{4 \beta - 9 \beta^2 - 4 \beta + 1 + 3 \beta^4 + 10 \beta^2 - 1}{4 (1 - \beta)^{5/2}(1 + 3 \beta)^{3/2}} \\
& \quad \quad \textstyle - \frac{\beta^2 - (1 + q^2) \beta}{4 (1 - q)^{3/2}(1 + 3q)^{1/2}} + \frac{4 \beta^4 - 9 \beta^4 - 4 \beta^3 + \beta^2 + 3 \beta^5 + 10 \beta^3 - \beta}{4 (1 - \beta)^{5/2}(1 + 3 \beta)^{3/2}}\bigr] \\
& \textstyle = \bigl[\frac{\beta (1-\beta)^{1/2}}{4 (1 + 3 \beta)^{1/2}} - \frac{3 \beta^4 + \beta^2}{4 (1 - \beta)^{5/2}(1 + 3 \beta)^{3/2}} + \frac{3 \beta^5 - 5 \beta^4 + 6 \beta^3 + \beta^2 - \beta}{4 (1 - \beta)^{5/2}(1 + 3 \beta)^{3/2}}\bigr] \\
& \quad \textstyle + \bigl[\frac{q(1 + q^2) - q^2}{4 (1 - q)^{3/2} (1 + 3q)^{1/2}} + \frac{- q^2}{4 (1 - q)^{3/2} (1 + 3q)^{1/2}} - \frac{q (1 - q)^{1/2}}{4 (1 + 3q)^{1/2}}\bigr] \\
& \textstyle = \bigl[\frac{\beta (1-\beta)^{1/2}}{4 (1 + 3 \beta)^{1/2}} + \frac{3 \beta^5 - 8 \beta^4 + 6 \beta^3 - \beta}{4 (1 - \beta)^{5/2}(1 + 3 \beta)^{3/2}}\bigr] + \bigl[\frac{q(1 + q^2 - 2q)}{4 (1 - q)^{3/2} (1 + 3q)^{1/2}} - \frac{q (1 - q)^{2}}{4 (1 - q)^{3/2} (1 + 3q)^{1/2}}\bigr] \\
& \textstyle = \frac{\beta (1 - \beta)^3 (1 + 3 \beta) + 3 \beta^5 - 8 \beta^4 + 6 \beta^3 - \beta}{4 (1 - \beta)^{5/2} (1 + 3 \beta)^{3/2}} = \frac{\beta (1 - \beta)^3 (1 + 3 \beta) + (\beta + 3 \beta^2)(\beta^3 - 3 \beta^2 + 3 \beta - 1)}{4 (1 - \beta)^{5/2} (1 + 3 \beta)^{3/2}} = 0,
\end{split}
\end{equation}
and
\begin{align}
& \textstyle \int_q^1 x (g(x) - \fg(x)) \, \d x \notag \\
& \textstyle = \int_q^1 x [g(x)] \, \d x - \int_q^{\beta} x [\fg(x)] \, \d x - \int_{\beta}^1 [x \fg(x)] \, \d x \notag \\
& \textstyle = \int_q^1 \frac{x(2 x - 1 - q^2)}{4(1 - q)^{3/2}(1 + 3q)^{1/2}} \, \d x - \int_q^{\beta} \frac{x (3 x^2 - 1)}{4 (1 - x)^{1/2} (1 + 3x)^{3/2}} \, \d x - \int_{\beta}^1 \frac{12 \beta x^3 - (18 \beta^2 + 8 \beta - 2) x^2 + (3 \beta^4 + 10 \beta^2 - 1) x}{4 (1 - \beta)^{5/2} (1 + 3 \beta)^{3/2}} \, \d x \notag \\
& \textstyle = \bigl[\frac{4 x^3 - 3 x^2 (1 + q^2)}{24 (1 - q)^{3/2} (1 + 3q)^{1/2}}\bigr]_{x = q}^{x = 1} + \bigl[\frac{(1 - x)^{1/2}(1 + 2x + 3x^2)}{24 (1 + 3x)^{1/2}} \bigr]_{x = q}^{x = \beta} - \bigl[\frac{18 \beta x^4 - 2x^3 (18 \beta^2 + 8 \beta - 2) + 3x^2 (3 \beta^4 + 10 \beta^2 - 1)}{24 (1 - \beta)^{5/2} (1 + 3 \beta)^{3/2}} \bigr]_{x = \beta}^{x = 1} \notag \\
& \textstyle = \frac{4 (1 - q^3) - 3 (1 - q^2)(1 + q^2)}{24 (1 - q)^{3/2} (1 + 3q)^{1/2}} + \bigl[\frac{(1 - \beta)^{1/2}(1 + 2 \beta + 3 \beta^2)}{24 (1 + 3 \beta)^{1/2}} - \frac{(1 - q)^{1/2}(1 + 2 q + 3 q^2)}{24 (1 + 3 q)^{1/2}}\bigr] \notag \\
& \quad \textstyle - \frac{18 \beta (1 - \beta^4) - 2 (1 - \beta^3) (18 \beta^2 + 8 \beta - 2) + 3 (1 - \beta^2) (3 \beta^4 + 10 \beta^2 - 1)}{24 (1 - \beta)^{5/2} (1 + 3 \beta)^{3/2}} \notag \\
& \textstyle = \bigl[\frac{4 (1 + q + q^2) - 3 (1 + q)(1 + q^2)}{24 (1 - q)^{1/2} (1 + 3q)^{1/2}} - \frac{(1 - q)^{1/2}(1 + 2 q + 3 q^2)}{24 (1 + 3 q)^{1/2}}\bigr] + \bigl[\frac{(1 - \beta)^{1/2}(1 + 2 \beta + 3 \beta^2)}{24 (1 + 3 \beta)^{1/2}} \notag \\
& \quad \textstyle - \frac{18 \beta (1 + \beta + \beta^2 + \beta^3) - 2 (1 + \beta + \beta^2) (18 \beta^2 + 8 \beta - 2) + 3 (1 + \beta) (3 \beta^4 + 10 \beta^2 - 1)}{24 (1 - \beta)^{3/2} (1 + 3 \beta)^{3/2}}\bigr] \notag \\
& \textstyle = \frac{[1 + q + q^2 - 3q^3] - (1 - q) (1 + 2q + 3q^2)}{24 (1 - q)^{1/2} (1 + 3q)^{1/2}} + \frac{(1 - \beta)^2 (1 + 3 \beta) (1 + 2 \beta + 3 \beta^2) - [9 \beta^5 - 9 \beta^4 - 4 \beta^3 + 3 \beta + 1]}{24 (1 - \beta)^{3/2}(1 + 3 \beta)^{3/2}} \notag \\
& \textstyle = \frac{(1 - \beta) (1 + 5 \beta + 9 \beta^2 + 9 \beta^3) + [9 \beta^4 - 4 \beta^2 - 4 \beta - 1]}{24 (1 - \beta)^{1/2}(1 + 3 \beta)^{3/2}} = 0.
\end{align}
This establishes \cref{item1:lemma:zero_derivative}. Moreover, \nobs that \cref{item1:lemma:zero_derivative} shows that
\begin{equation}
\begin{split}
\textstyle \int_0^1 [g(x) \fg(x)] \, \d x & \textstyle = \int_0^1 [g(x)]^2 \, \d x - \int_0^1 g(x) [g(x) - \fg(x)] \, \d x \\
& \textstyle = \int_0^1 [g(x)]^2 \, \d x + \frac{(1 - q)^{1/2}}{4 (1 + 3q)^{1/2}}\int_0^q (g(x) - \fg(x)) \, \d x \\
& \quad \textstyle - \frac{1}{4(1 - q)^{3/2}(1 + 3q)^{1/2}} \int_q^1 (2x - 1 - q^2)(g(x) - \fg(x)) \, \d x \\
& \textstyle = \int_0^1 [g(x)]^2 \, \d x
\end{split}
\end{equation}

and
\begin{equation}
\begin{split}
\textstyle \int_0^1 [g(x)]^2 \, \d x & \textstyle = \int_0^q [-\frac{(1 - q)^{1/2}}{4(1 + 3q)^{1/2}}]^2 \, \d x + \int_q^1 [\frac{2x - 1 - q^2}{4(1 - q)^{3/2}(1 + 3q)^{1/2}}]^2 \, \d x \\
& \textstyle = \frac{q (1 - q)}{16 (1 + 3q)} + \bigl[\frac{(2x - 1 - q^2)^3}{96 (1 - q)^3(1 + 3q)}\bigr]_{x = q}^{x = 1} = \frac{1}{48}.
\end{split}
\end{equation}
\Hence that
\begin{equation}
\begin{split}
\textstyle \int_0^1 (g(x) - \fg(x))^2 \, \d x & \textstyle = \int_0^1 [\fg(x)]^2 \, \d x - 2 \int_0^1 [g(x) \fg(x)] \, \d x + \int_0^1 [g(x)]^2 \, \d x \\
& \textstyle = \int_0^1 [\fg(x)]^2 \, \d x - \int_0^1 [g(x)]^2 \, \d x = \int_0^1 [\fg(x)]^2 \, \d x - \frac{1}{48}.
\end{split}
\end{equation}
This establishes \cref{item2:lemma:zero_derivative}.
\end{cproof}

\subsection{On a submanifold of the ANN parameter space}

\cfclear
\begin{lemma}\label{lemma:set_of_critical_points_submanifold}
Assume \cref{setting0} \cfload. Then it holds that $\cM$ is a $(\dimension - d - 1)$-dimensional $C^{\infty}$-\submanifold of $\R^{\dimension}$ \cfout.
\end{lemma}
\begin{cproof}{lemma:set_of_critical_points_submanifold}
Throughout this proof let $z \allowbreak \in \cM$, $\eps \allowbreak \in (0, \allowbreak \infty)$, let $\scrU \subseteq \R^{\dimension - d - 1}$ satisfy
\begin{multline}\label{eqn:lemma:set_of_critical_points_submanifold:fU}
\textstyle \scrU = \bigl\{x = (x_1, \dots, x_{\dimension - d - 1}) \in \R^{\dimension - d - 1} \colon \bigl[(\alpha < x_1 < \beta) \wedge (x_2 > 0) \\
\textstyle \wedge\! \bigl(\Forall t_1, t_2, \dots, t_d \in [\scra, \scrb] \colon\! \bigwedge_{i = 2}^\width [x_{d (\width - 1) + i + 1} + \sum_{j = 1}^d x_{(i - 2)d + j + 2} t_j < 0]\bigr)\bigr]\bigr\},
\end{multline}
let $U \subseteq \R^{\dimension - d - 1}$ satisfy
\begin{multline}\label{eqn:lemma:set_of_critical_points_submanifold:scrU}
\textstyle U = \Bigl\{x = (x_1, \dots, x_{\dimension - d - 1}) \in \scrU \colon \bigl(\Exists y \in \allowbreak \R^{\dimension} \colon [\norm{y - z} < \eps]\! \wedge\! \bigl[\fc^y = - \frac{(1 - x_1)^{1/2}}{4 (1 + 3x_1)^{1/2}}\bigr] \\
\textstyle \wedge\! \bigl[\sum_{j = 2}^d \abs{\fw_{1, j}^y} = - \abs{\fw_{1, 1}^y - \frac{x_2}{\scrb - \scra}}\bigr]\! \wedge\! [\fb_1^y = - x_1 x_2 - \frac{x_2 \scra}{\scrb - \scra}]\! \wedge\! \bigl[\fv_1^y = \frac{(1 + 3x_1)^{- 1/2}}{2 x_2 (1 - x_1)^{3/2}}\bigr] \\
\textstyle \wedge\! \bigl[\bigwedge_{i = 2}^\width \bigwedge_{j = 1}^d [(\fv_i^y = x_{d(\width - 1) + \width + i})\! \wedge\! (\fb_i^y = x_{d(\width - 1) + i + 1})\! \wedge\! (\fw_{i, j}^y = x_{(i - 2)d + j + 2})]\bigr]\bigr)\!\Bigr\},
\end{multline}
and let $\varphi \colon U \to \R^{\dimension}$ satisfy for all $x = (x_1, \dots, \allowbreak x_{\dimension - d - 1}) \in U$ that
\begin{equation}\label{eqn:lemma:set_of_critical_points_submanifold:phi}
\begin{gathered}
\textstyle \fw_{1, 1}^{\varphi(x)} = \frac{x_2}{\scrb - \scra}, \qquad \sum_{j = 2}^d \abs{\fw_{1, j}^{\varphi(x)}} = 0, \qquad \fc^{\varphi(x)} = - \frac{(1 - x_1)^{1/2}}{4 (1 + 3x_1)^{1/2}}, \\
\textstyle \bigwedge_{i = 2}^\width \bigwedge_{j = 1}^d [(\fv_i^{\varphi(x)} = x_{d(\width - 1) + \width + i})\! \wedge\! (\fb_i^{\varphi(x)} = x_{d(\width - 1) + i + 1})\! \wedge\! (\fw_{i, j}^{\varphi(x)} = x_{(i - 2)d + j + 2})], \\
\textstyle \fb_1^{\varphi(x)} = - x_1 x_2 - \frac{x_2 \scra}{\scrb - \scra}, \qqandqq \fv_1^{\varphi(x)} = \frac{(1 + 3x_1)^{- 1/2}}{2 x_2 (1 - x_1)^{3/2}}\ifnocf.
\end{gathered}
\end{equation}
\cfload[.]\Nobs that \cref{eqn:lemma:set_of_critical_points_submanifold:fU,eqn:lemma:set_of_critical_points_submanifold:scrU} assure that $U$ is open in $\R^{\dimension - d - 1}$. Next \nobs that \cref{eqn:lemma:set_of_critical_points_submanifold:phi} ensures that for all $x = (x_1, \dots, x_{\dimension - d - 1}) \in U$ it holds that
\begin{multline}\label{eqn:lemma:set_of_critical_points_submanifold:phi_sub_Jacobian}
\textstyle \tfrac{\d}{\d x}(\varphi_{d \width + 1}, \varphi_1, \varphi_{d + 1}, \dots, \varphi_{d \width}, \varphi_{d \width + 2}, \dots, \varphi_{d \width + \width}, \varphi_{d \width + \width + 2}, \dots, \varphi_{d \width + 2 \width})(x) \\
=
\begin{pmatrix}
(-x_2) & (-x_1 - \frac{\scra}{\scrb - \scra}) & 0 & 0 & \cdots & 0 \\
   0 &   \frac{1}{\scrb - \scra}  & 0 & 0 & \cdots & 0 \\
   0 &   0  & 1 & 0 & \cdots & 0 \\
   0 &   0  & 0 & 1 & \cdots & 0 \\
\vdots & \vdots & \vdots & \vdots & \ddots & \vdots \\
   0 &   0  & 0 & 0 & \cdots & 1
\end{pmatrix}\! \in \R^{(\dimension - d - 1) \times (\dimension - d - 1)}.
\end{multline}
This, \cref{eqn:lemma:set_of_critical_points_submanifold:fU}, and \cref{eqn:lemma:set_of_critical_points_submanifold:scrU} show that for all $x \in U$ it holds that $\rk(\varphi'(x)) = \dimension - d - 1$. Combining this with the fact that $\varphi \in C^{\infty} (U, \R^{\dimension})$ implies that $U \ni x \mapsto \varphi(x) \in \R^{\dimension}$ is a $C^{\infty}$-immersion from $U$ to $\R^{\dimension}$ \cfadd{def:immersion}\cfload. Next \nobs that \cref{eqn:lemma:set_of_critical_points_submanifold:scrU,eqn:lemma:set_of_critical_points_submanifold:phi} ensure that for all $x \in U$ it holds that $\norm{\varphi(x) - z} < \eps$ \cfload. Combining this with \cref{0_eq:M,eqn:lemma:set_of_critical_points_submanifold:fU,eqn:lemma:set_of_critical_points_submanifold:scrU,eqn:lemma:set_of_critical_points_submanifold:phi} assures that 
\begin{equation}\label{eqn:lemma:set_of_critical_points_submanifold:phi_cM}
\varphi(U) = \cM \cap \{x \in \R^{\dimension} \colon \norm{x - z} < \eps\}.
\end{equation}
Next \nobs that \cref{eqn:lemma:set_of_critical_points_submanifold:fU,eqn:lemma:set_of_critical_points_submanifold:scrU,eqn:lemma:set_of_critical_points_submanifold:phi} show that for all $x = \allowbreak (x_1, \allowbreak \dots, \allowbreak x_{\dimension - d - 1}), y = \allowbreak (y_1, \allowbreak \dots, \allowbreak y_{\dimension - d - 1}) \allowbreak \in U$ with $\varphi(x) = \allowbreak \varphi(y)$ it holds that $x = y$. This shows that $U \ni x \mapsto \varphi(x) \in \varphi(U)$ is bijective. Moreover, \nobs that \cref{eqn:lemma:set_of_critical_points_submanifold:phi} shows for all $x \allowbreak = (x_1, \allowbreak \dots, \allowbreak x_{\dimension - d - 1}) \allowbreak \in U$, $\theta \in \R^{\dimension}$ with $\theta \allowbreak = \varphi(x)$ that $x_1 \allowbreak = (- [\fb_1^{\theta} / \fw_{1,1}^{\theta}] - \scra) / (\scrb - \scra)$, $x_2 \allowbreak = (\scrb - \scra) \fw_{1,1}^{\theta}$, and that the remaining coordinates of $x$ are the corresponding inner weights, inner biases, and outer weights of the neurons of $\theta$ with indices in $\{2, \allowbreak 3, \allowbreak \dots, \allowbreak \width\}$. The fact that for all $\theta \allowbreak \in \varphi(U)$ it holds that $\fw_{1,1}^{\theta} \allowbreak > 0$ \hence ensures that the inverse of $U \ni x \allowbreak \mapsto \allowbreak \varphi(x) \allowbreak \in \varphi(U)$ is continuous. Combining this with the fact that $\varphi\in C^{\infty}(U,\R^{\dimension})$ demonstrates that 
\begin{equation}
U \ni x \mapsto \varphi(x) \in \varphi(U)
\end{equation}
is a homeomorphism. The fact that $U \ni \allowbreak x \allowbreak \mapsto \allowbreak \varphi(x) \allowbreak \in \R^{\dimension}$ is a $C^{\infty}$-immersion from $U$ to $\R^{\dimension}$ and \cref{eqn:lemma:set_of_critical_points_submanifold:phi_cM} \hence show that $\cM$ is a $(\dimension - d - 1)$-dimensional $C^{\infty}$-\submanifold of $\R^{\dimension}$ \cfload.
\end{cproof}

\subsection{Explicit representations for the Hessian matrix of the risk function}
\label{subsec:Hessian_explicit_cM}

\cfclear
\begin{lemma}\label{lemma:Hessian_of_Risk_cM}
Assume \cref{setting0}. Then
\begin{enumerate}[label=(\roman*)]
\item
\label{item_I:lemma:Hessian_of_Risk_cM} it holds for all $\theta \in \R^{\dimension}$, $i \in \{1, 2, \dots, \width\}$ with $I_i^{\theta} = \emptyset = \partial I_i^{\theta}$ that
\begin{equation}
\textstyle \Exists \eps \in (0, \infty) \colon \bigl[\Forall \vartheta \in [- \eps, \eps]^{\dimension} \colon (I_i^{\theta + \vartheta} = \emptyset = \partial I_i^{\theta + \vartheta}) \wedge (\abs{\fb_i^{\theta + \vartheta}} + \norm{\fw_{i}^{\theta + \vartheta}} > 0)\bigr],
\end{equation}

\item
\label{item_fV_M:lemma:Hessian_of_Risk_cM} it holds that $\fV \subseteq \R^{\dimension}$ is open and $\cM \subseteq \fV$,

\item
\label{item_fV_I:lemma:Hessian_of_Risk_cM} it holds for all $\theta = (\theta_1, \dots, \theta_{\dimension}) \in \fV$, $y_2, y_3, \dots, y_d \in [\scra, \scrb]$ that
\begin{gather}
\textstyle \theta_1 > 0, \quad \alpha < \fq_{\theta} - \frac{1}{(\scrb - \scra) \theta_1} \sum_{j = 2}^d \theta_j y_j < \beta, \quad \scra < \scrq_{\theta} - \frac{1}{\theta_1} \sum_{j = 2}^d \theta_j y_j < \scrb, \\
\textstyle I_1^{\theta} = \{x = (x_1, \dots, x_d) \in [\scra, \scrb]^d \colon \scrq_{\theta} - [\theta_1]^{-1}\sum_{j = 2}^d \theta_j x_j < x_1 \le \scrb\}, \\
\textstyle \andqq \partial I_1^{\theta} = \{x = (x_1, \dots, x_d) \in [\scra, \scrb]^d \colon x_1 = \scrq_{\theta} - [\theta_1]^{-1}\sum_{j = 2}^d \theta_j x_j\},
\end{gather}

\item
\label{item_M_I:lemma:Hessian_of_Risk_cM} it holds for all $\theta \in \cM$ that $\theta_1 > 0$, $\scra < \scrq_{\theta} < \scrb$, $\alpha < \fq_{\theta} < \beta$, and\footnote{\Nobs that for all $A, B \subseteq \R$ it holds that $B^0 = \times_{k = 1}^0 B = \{\emptyset\}$ and $A \times B^0 = A$ and it holds for all $f \colon A \allowbreak \times \allowbreak B^0 \allowbreak \to \R$ that $\Forall x \allowbreak \in A, \allowbreak y \allowbreak \in B^0 \colon [f(x, y) = f(x)] \wedge [\int_{B^0} f(x, t) \d t = f(x)]$.}
\begin{equation}
\textstyle (I_1^{\theta} = (\scrq_{\theta}, \scrb] \times [\scra, \scrb]^{d - 1}) \wedge (\partial I_1^{\theta} = \{\scrq_{\theta}\} \times [\scra, \scrb]^{d - 1}) \wedge \bigl(\bigwedge_{i = 2}^\width [I_i^{\theta} = \emptyset = \partial I_i^{\theta}]\bigr),
\end{equation}

\item
\label{item_fV_C2:lemma:Hessian_of_Risk_cM} it holds that $(\riskR)|_{\fV} \in C^2 (\fV, \R)$, and

\item
\label{item_Hess:lemma:Hessian_of_Risk_cM} it holds for all $\theta = (\theta_1, \dots, \theta_{\dimension}) \in \fV$, $j, \allowbreak \fj \allowbreak \in \{1, \allowbreak 2, \allowbreak \dots, \allowbreak d\}$ that\footnote{\Nobs that we denote by $\sigma$ the $(d - 1)$-dimensional surface measure on smooth hypersurfaces in $\R^d$.}
\begin{gather}
\textstyle \bigwedge_{k \in \{1, 2, \dots, \dimension\} \backslash \{\pi(1), \pi(2), \dots, \pi(d + 3)\}} \bigwedge_{\ell \in \{1, 2, \dots, \dimension\}} \bigl[(\frac{\partial}{\partial \theta_k} \riskR)(\theta) = 0 = (\frac{\partial^2}{\partial \theta_k \partial \theta_{\ell}} \riskR)(\theta)\bigr], \\
\textstyle \norm{\fw_1^{\theta}} \ge \theta_1 > 0, \quad (\frac{\partial^2}{\partial \theta_j \partial \theta_{\dimension}} \riskR)(\theta) = 2 \fv_1^{\theta} \int_{I_1^{\theta}} \spro{e_j, x} \d x, \quad (\frac{\partial^2}{\partial \theta_{\dimension}^2} \riskR)(\theta) = 2 (\scrb - \scra)^d, \\
\textstyle (\frac{\partial^2}{\partial \theta_{d \width + 1} \partial \theta_{\dimension}} \riskR)(\theta) = 2 \fv_1^{\theta} \int_{I_1^{\theta}} \d x, \qquad (\frac{\partial^2}{\partial \theta_{d \width + \width + 1} \partial \theta_{\dimension}} \riskR)(\theta) = 2 \int_{I_1^{\theta}} (\fb_1^{\theta} + \spro{\fw_1^{\theta}, x}) \d x, \\
\textstyle (\frac{\partial^2}{\partial \theta_j \partial \theta_{d \width + \width + 1}} \riskR)(\theta) = 2 \fv_1^{\theta} \int_{I_1^{\theta}} \spro{e_j, x} (\fb_1^{\theta} + \spro{\fw_1^{\theta}, x}) \d x + 2 \int_{I_1^{\theta}} \spro{e_j, x} (\functionSNN^{\theta}(x) - f(x)) \d x, \\
\textstyle (\frac{\partial^2}{\partial \theta_{d \width + 1} \partial \theta_{d \width + \width + 1}} \riskR)(\theta) = 2 \fv_1^{\theta} \int_{I_1^{\theta}} (\fb_1^{\theta} + \spro{\fw_1^{\theta}, x}) \d x + 2 \int_{I_1^{\theta}} (\functionSNN^{\theta}(x) - f(x)) \d x, \\
\textstyle (\frac{\partial^2}{\partial \theta_{d \width + \width + 1}^2} \riskR)(\theta) = 2 \int_{I_1^{\theta}} (\fb_1^{\theta} + \spro{\fw_1^{\theta}, x})^2 \d x, \\
\textstyle (\frac{\partial^2}{\partial \theta_j \partial \theta_{\fj}} \riskR)(\theta) = 2 [\fv_1^{\theta}]^2 \int_{I_1^{\theta}} \spro{e_j, x} \spro{e_{\fj}, x} \d x + 2 \fv_1^{\theta} \int_{\partial I_1^{\theta}} \frac{\spro{e_j, x} \spro{e_{\fj}, x} (\functionSNN^{\theta}(x) - f(x))}{\norm{\fw_1^{\theta}}} \d \sigma(x), \\
\textstyle (\frac{\partial^2}{\partial \theta_{d \width + 1} \partial \theta_j} \riskR)(\theta) = 2 [\fv_1^{\theta}]^2 \int_{I_1^{\theta}} \spro{e_j, x} \d x + 2 \fv_1^{\theta} \int_{\partial I_1^{\theta}} \frac{\spro{e_j, x} (\functionSNN^{\theta}(x) - f(x))}{\norm{\fw_1^{\theta}}} \d \sigma(x), \\
\textstyle \andqq (\frac{\partial^2}{\partial \theta_{d \width + 1}^2} \riskR)(\theta) = 2 [\fv_1^{\theta}]^2 \int_{I_1^{\theta}} \d x + 2 \fv_1^{\theta} \int_{\partial I_1^{\theta}} \frac{(\functionSNN^{\theta}(x) - f(x))}{\norm{\fw_1^{\theta}}} \d \sigma(x).
\end{gather}
\end{enumerate}
\cfload[.]
\end{lemma}
\begin{cproof}{lemma:Hessian_of_Risk_cM}
\Nobs that \cref{0_eq:M}, \cref{0_eq:fV}, and the fact that $[\scra, \scrb]^d$ is a compact subset of $\R^d$ establish \cref{item_I:lemma:Hessian_of_Risk_cM,item_fV_M:lemma:Hessian_of_Risk_cM}. Next \nobs that \cref{0_eq:fV}, the triangle inequality, and the fact that for all $t \in [\scra, \scrb]$ it holds that $\abs{t} \le \abs{\scrb} + \abs{\scra}$ demonstrate for all $\theta \allowbreak = (\theta_1, \allowbreak \dots, \allowbreak \theta_{\dimension}) \allowbreak \in \fV$, $y_2, \allowbreak y_3, \allowbreak \dots, \allowbreak y_d \allowbreak \in [\scra, \allowbreak \scrb]$ that $\theta_1 > 0$, $\alpha < \fq_{\theta} < \beta$, $\scra < \scrq_{\theta} < \scrb$, and
\begin{equation}
\textstyle \Abs{\sum_{j = 2}^d \theta_j y_j} \le (\abs{\scrb} + \abs{\scra}) \sum_{j = 2}^d \abs{\theta_j} < \theta_1 (\scrb - \scra) \min\{\fq_{\theta} - \alpha, \beta - \fq_{\theta}\}.
\end{equation}
\Hence for all $\theta \allowbreak = (\theta_1, \allowbreak \dots, \allowbreak \theta_{\dimension}) \allowbreak \in \fV$, $y_2, \allowbreak y_3, \allowbreak \dots, \allowbreak y_d \allowbreak \in [\scra, \allowbreak \scrb]$ that
\begin{equation}
\textstyle \theta_1 > 0, \qquad \alpha < \fq_{\theta} - \frac{1}{(\scrb - \scra) \theta_1} \sum_{j = 2}^d \theta_j y_j < \beta, \qqandqq \scra < \scrq_{\theta} - \frac{1}{\theta_1} \sum_{j = 2}^d \theta_j y_j < \scrb.
\end{equation}
This establishes \cref{item_fV_I:lemma:Hessian_of_Risk_cM}. Moreover, \nobs that \cref{item_fV_M:lemma:Hessian_of_Risk_cM}, \cref{item_fV_I:lemma:Hessian_of_Risk_cM}, \cref{0_eq:M}, and \cref{0_eq:fV} establish \cref{item_M_I:lemma:Hessian_of_Risk_cM}. For every $k \allowbreak \in \{0, \allowbreak 1\}$, $j \in \{1, 2, \dots, d\}$ let $F_{k, j} \allowbreak \colon \fV \times [\scra, \allowbreak \scrb]^{d - 1} \allowbreak \to \R$ and $\varphi \allowbreak \colon \fV \times [\scra, \allowbreak \scrb]^{d - 1} \allowbreak \to \R$ satisfy for all $\theta \allowbreak = (\theta_1, \allowbreak \dots, \allowbreak \theta_{\dimension}) \allowbreak \in \fV$, $x \allowbreak = (x_2, \allowbreak \dots, \allowbreak x_d) \allowbreak \in [\scra, \allowbreak \scrb]^{d - 1}$ that
\begin{equation}\label{eqn:lemma:Hessian_of_Risk_cM:F}
\textstyle \varphi(\theta, x) = \scrq_{\theta} - \sum_{i = 2}^d \frac{\theta_i x_i}{\theta_1} \qandq F_{k, j}(\theta, x) = \int_{\varphi(\theta, x)}^{\scrb} (k x_j + 1 - k) (\functionSNN^{\theta}(x_1, x) - f(x_1,x)) \d x_1
\end{equation}
(cf.\ \Cref{item_fV_I:lemma:Hessian_of_Risk_cM}). \Nobs that \cref{item_fV_I:lemma:Hessian_of_Risk_cM}, \cref{eqn:lemma:Hessian_of_Risk_cM:F}, the fact that for all $\theta \allowbreak = (\theta_1, \allowbreak \dots, \allowbreak \theta_{\dimension}) \allowbreak \in \fV$, $x \allowbreak = (x_2, \allowbreak \dots, \allowbreak x_d) \allowbreak \in [\scra, \allowbreak \scrb]^{d - 1}$, $\fj \allowbreak \in \{1, \allowbreak 2, \allowbreak \dots, \allowbreak d\}$ it holds that
\begin{gather}
\textstyle \varphi \in C^{\infty}(\fV \times (\scra, \scrb)^{d - 1}, \R), \quad \scra < \varphi(\theta, x) < \scrb, \quad (\nabla_x \varphi)(\theta, x) = - [\theta_1]^{-1} (\theta_2, \allowbreak \dots, \allowbreak \theta_d), \\
\textstyle (\frac{\partial}{\partial \theta_{d \width + 1}} \varphi)(\theta, x) = - [\theta_1]^{-1}, \qqandqq (\frac{\partial}{\partial \theta_{\fj}} \varphi)(\theta, x) =
    \begin{cases}
    - [\theta_1]^{-1} x_{\fj} & \colon \fj > 1 \\
    - [\theta_1]^{-1} \varphi(\theta, x) & \colon \fj = 1,
    \end{cases}
\end{gather}
and Leibniz's integral rule ensure for all $k \allowbreak \in \{0, \allowbreak 1\}$, $j, \allowbreak \fj \allowbreak \in \{1, \allowbreak 2, \allowbreak \dots, \allowbreak d\}$, $\theta \allowbreak = (\theta_1, \allowbreak \dots, \allowbreak \theta_{\dimension}) \allowbreak \in \fV$, $x \allowbreak = (x_2, \allowbreak \dots, \allowbreak x_d) \allowbreak \in [\scra, \allowbreak \scrb]^{d - 1}$, $g \in C([\scra, \scrb]^d, \R)$ that
\begin{equation}\label{eqn:lemma:Hessian_of_Risk_cM:int_transform}
\textstyle \int_{\partial I_1^{\theta}} \frac{g(x)}{\norm{\fw_1^{\theta}}} \d \sigma(x) = \int_{[\scra, \scrb]^{d - 1}} \frac{g(\varphi(\theta, y), y)}{\norm{\fw_1^{\theta}}} \sqrt{1 + \norm{(\nabla_y \varphi)(\theta,y)}^2} \, \d y = \int_{[\scra, \scrb]^{d - 1}} \frac{g(\varphi(\theta, y), y)}{\theta_1} \d y,
\end{equation}
\begin{multline}\label{eqn:lemma:Hessian_of_Risk_cM:F_grad_w}
\textstyle (\frac{\partial}{\partial \theta_{\fj}} F_{k, j})(\theta, x) = \fv_1^{\theta} \int_{\varphi(\theta, x)}^{\scrb} (k x_j + 1 - k) x_{\fj} \d x_1 \\
\textstyle + [k \spro{e_j, (\varphi(\theta, x), x)} + 1 - k] \spro{e_{\fj}, (\varphi(\theta, x), x)} [\theta_1]^{-1} [\functionSNN^{\theta}(\varphi(\theta, x), x) - f(\varphi(\theta,x),x)],
\end{multline}
\begin{equation}\label{eqn:lemma:Hessian_of_Risk_cM:F_grad_b}
\textstyle \andqq (\frac{\partial}{\partial \theta_{d \width + 1}} F_{0, 1})(\theta, x) = \fv_1^{\theta} \int_{\varphi(\theta, x)}^{\scrb} \d x_1 + \bigl([\theta_1]^{-1} [\functionSNN^{\theta}(\varphi(\theta, x), x) - f(\varphi(\theta,x),x)]\bigr).
\end{equation}
Next \nobs that \cref{item_I:lemma:Hessian_of_Risk_cM}, \cref{item_fV_I:lemma:Hessian_of_Risk_cM}, \cref{prop:limit:lr}, \cref{prop:loss:diff:vc}, \cref{prop:loss:diff:wb}, \cref{0_eq:f}, \cref{0_eq:fV}, \cref{eqn:lemma:Hessian_of_Risk_cM:F}, and Fubini's theorem assure for all $\theta \allowbreak = (\theta_1, \allowbreak \dots, \allowbreak \theta_{\dimension}) \allowbreak \in \fV$, $j \allowbreak \in \{1, \allowbreak 2, \allowbreak \dots, \allowbreak d\}$ that
\begin{gather}
\textstyle (\frac{\partial}{\partial \theta_j} \riskR)(\theta) = 2 \fv_1^{\theta} \int_{I_1^\theta} \scalprod{e_j}{x} (\functionSNN^{\theta}(x) - f (x)) \d x = 2 \fv_1^{\theta} \int_{[\scra, \scrb]^{d - 1}} F_{1, j}(\theta, t) \d t, \\
\textstyle (\frac{\partial}{\partial \theta_{d \width + 1}} \riskR)(\theta) = 2 \fv_1^{\theta} \int_{I_1^\theta} (\functionSNN^{\theta}(x) - f(x)) \d x = 2 \fv_1^{\theta} \int_{[\scra, \scrb]^{d - 1}} F_{0, 1}(\theta, t) \d t, \\
\textstyle (\frac{\partial}{\partial \theta_{d \width + \width + 1}} \riskR)(\theta) = 2 \int_{I_1^{\theta}} (\fb_1^{\theta} + \spro{\fw_1^{\theta}, x}) (\functionSNN^{\theta}(x) - f(x)) \d x, \\
\textstyle (\frac{\partial}{\partial \theta_{\dimension}} \riskR)(\theta) = 2 \int_{[\scra, \scrb]^d} (\functionSNN^{\theta}(x) - f(x)) \d x, \\ \label{eqn:lemma:Hessian_of_Risk_cM:grad_risk_rest}
\textstyle \andqq \bigwedge_{k \in \{1, 2, \dots, \dimension\} \backslash \{\pi(1), \pi(2), \dots, \pi(d + 3)\}} \bigwedge_{\ell \in \{1, 2, \dots, \dimension\}} \bigl[(\frac{\partial}{\partial \theta_k} \riskR)(\theta) = 0 = (\frac{\partial^2}{\partial \theta_k \partial \theta_{\ell}} \riskR)(\theta)\bigr].
\end{gather}
Combining \cref{cor:interchange}, \cref{item_fV_I:lemma:Hessian_of_Risk_cM}, \cref{0_eq:f}, \cref{eqn:lemma:Hessian_of_Risk_cM:F}, \cref{eqn:lemma:Hessian_of_Risk_cM:int_transform}, \cref{eqn:lemma:Hessian_of_Risk_cM:F_grad_w}, \cref{eqn:lemma:Hessian_of_Risk_cM:F_grad_b}, and Fubini's theorem \hence shows for all $\theta \allowbreak = (\theta_1, \allowbreak \dots, \allowbreak \theta_{\dimension}) \allowbreak \in \fV$, $j, \allowbreak \fj \allowbreak \in \{1, \allowbreak 2, \allowbreak \dots, \allowbreak d\}$ that
\begin{gather}
\textstyle (\frac{\partial^2}{\partial \theta_j \partial \theta_{\dimension}} \riskR)(\theta) = 2 \fv_1^{\theta} \int_{I_1^{\theta}} \spro{e_j, x} \d x, \qquad (\frac{\partial^2}{\partial \theta_{d \width + 1} \partial \theta_{\dimension}} \riskR)(\theta) = 2 \fv_1^{\theta} \int_{I_1^{\theta}} \d x, \\
\textstyle (\frac{\partial^2}{\partial \theta_{d \width + \width + 1} \partial \theta_{\dimension}} \riskR)(\theta) = 2 \int_{I_1^{\theta}} (\fb_1^{\theta} + \spro{\fw_1^{\theta}, x}) \d x, \qquad (\frac{\partial^2}{\partial \theta_{\dimension}^2} \riskR)(\theta) = 2 (\scrb - \scra)^d, \\
\textstyle (\frac{\partial^2}{\partial \theta_j \partial \theta_{d \width + \width + 1}} \riskR)(\theta) = 2 \fv_1^{\theta} \int_{I_1^{\theta}} \spro{e_j, x} (\fb_1^{\theta} + \spro{\fw_1^{\theta}, x}) \d x + 2 \int_{I_1^{\theta}} \spro{e_j, x} (\functionSNN^{\theta}(x) - f(x)) \d x, \\
\textstyle (\frac{\partial^2}{\partial \theta_{d \width + 1} \partial \theta_{d \width + \width + 1}} \riskR)(\theta) = 2 \fv_1^{\theta} \int_{I_1^{\theta}} (\fb_1^{\theta} + \spro{\fw_1^{\theta}, x}) \d x + 2 \int_{I_1^{\theta}} (\functionSNN^{\theta}(x) - f(x)) \d x, \\
\textstyle (\frac{\partial^2}{\partial \theta_{d \width + \width + 1}^2} \riskR)(\theta) = 2 \int_{I_1^{\theta}} (\fb_1^{\theta} + \spro{\fw_1^{\theta}, x})^2 \d x,
\end{gather}
\begin{align}
\textstyle (\frac{\partial^2}{\partial \theta_j \partial \theta_{\fj}} \riskR)(\theta) & \textstyle = 2 \fv_1^{\theta} \int_{[\scra, \scrb]^{d - 1}} (\frac{\partial}{\partial \theta_{\fj}} F_{1, j})(\theta, t) \d t \notag \\
& \textstyle = 2 [\fv_1^{\theta}]^2 \int_{[\scra, \scrb]^{d - 1}}\bigl(\int_{\varphi(\theta, x_2, \dots, x_d)}^{\scrb} x_j x_{\fj} \d x_1\bigr) \d x_2 \dots \d x_d \notag \\
& \textstyle \quad + 2 \fv_1^{\theta} \int_{[\scra, \scrb]^{d - 1}} \bigl([\functionSNN^{\theta}(\varphi(\theta, x), x) - f(\varphi(\theta,x),x)] \frac{\spro{e_j, (\varphi(\theta, x), x)} \spro{e_{\fj}, (\varphi(\theta, x), x)}}{\theta_1}\bigr) \d x \notag \\
& \textstyle = 2 [\fv_1^{\theta}]^2 \int_{I_1^{\theta}} \spro{e_j, x} \spro{e_{\fj}, x} \d x + 2 \fv_1^{\theta} \int_{\partial I_1^{\theta}} \frac{\spro{e_j, x} \spro{e_{\fj}, x} (\functionSNN^{\theta}(x) - f(x))}{\norm{\fw_1^{\theta}}} \d \sigma(x),
\end{align}
\begin{equation}
\begin{split}
\textstyle (\frac{\partial^2}{\partial \theta_{d \width + 1} \partial \theta_j} \riskR)(\theta) & \textstyle = 2 \fv_1^{\theta} \int_{[\scra, \scrb]^{d - 1}} (\frac{\partial}{\partial \theta_j} F_{0, 1})(\theta, t) \d t \\
& \textstyle = 2 [\fv_1^{\theta}]^2 \int_{[\scra, \scrb]^{d - 1}}\bigl(\int_{\varphi(\theta, x_2, \dots, x_d)}^{\scrb} x_j \d x_1\bigr) \d x_2 \dots \d x_d \\
& \textstyle \quad + 2 \fv_1^{\theta} \int_{[\scra, \scrb]^{d - 1}} \bigl([\functionSNN^{\theta}(\varphi(\theta, x), x) - f(\varphi(\theta,x),x)] \frac{\spro{e_j, (\varphi(\theta, x), x)}}{\theta_1}\bigr) \d x \\
& \textstyle = 2 [\fv_1^{\theta}]^2 \int_{I_1^{\theta}} \spro{e_j, x} \d x + 2 \fv_1^{\theta} \int_{\partial I_1^{\theta}} \frac{\spro{e_j, x} (\functionSNN^{\theta}(x) - f(x))}{\norm{\fw_1^{\theta}}} \d \sigma(x),
\end{split}
\end{equation}
\begin{equation}
\begin{split}
\textstyle (\frac{\partial^2}{\partial \theta_{d \width + 1}^2} \riskR)(\theta) & \textstyle = 2 \fv_1^{\theta} \int_{[\scra, \scrb]^{d - 1}} (\frac{\partial}{\partial \theta_{d \width + 1}} F_{0, 1})(\theta, t) \d t \\
& \textstyle = 2 [\fv_1^{\theta}]^2 \int_{[\scra, \scrb]^{d - 1}}\bigl(\int_{\varphi(\theta, x_2, \dots, x_d)}^{\scrb} \d x_1\bigr) \d x_2 \dots \d x_d \\
& \textstyle \quad + 2 \fv_1^{\theta} \int_{[\scra, \scrb]^{d - 1}} \bigl([\theta_1]^{-1} [\functionSNN^{\theta}(\varphi(\theta, x), x) - f(\varphi(\theta,x),x)]\bigr) \d x \\
& \textstyle = 2 [\fv_1^{\theta}]^2 \int_{I_1^{\theta}} \d x + 2 \fv_1^{\theta} \int_{\partial I_1^{\theta}} \frac{(\functionSNN^{\theta}(x) - f(x))}{\norm{\fw_1^{\theta}}} \d \sigma(x),
\end{split}
\end{equation}
and $(\riskR)|_{\fV} \in C^2(\fV, \R)$. This and \cref{eqn:lemma:Hessian_of_Risk_cM:grad_risk_rest} establish \cref{item_fV_C2:lemma:Hessian_of_Risk_cM,item_Hess:lemma:Hessian_of_Risk_cM}.
\end{cproof}

\subsection{On the rank and the spectrum of the Hessian matrix of the risk function}
\label{subsec:rank_of_hessian_of_risk_function}

\cfclear
\begin{lemma}[\resname{Schur complement criterion for positive definiteness}]\label{lemma:Schur_complement}
Let $m, n \in \N$ and let $A \in \R^{n \times n}$, $B \in \R^{n \times m}$, $C \allowbreak \in \allowbreak \R^{m \times m}$, $D \in \allowbreak \R^{(n + m) \times (n + m)}$ satisfy\footnote{\Nobs that for every matrix $A \in \cup_{m \in \N} \allowbreak (\cup_{n \in \N} \allowbreak \R^{m \times n})$ we denote by $A^{\top}$ the transpose of $A$ and for every $n \allowbreak \in \allowbreak \N$ we denote by $\idMatrix_n$ the $n \times n$ identity matrix. Moreover, \nobs that for every $n \in \N$ and every symmetric matrix $A \allowbreak \in \allowbreak \R^{n \times n}$ we say that $A$ is positive definite (we write $A \succ 0$) if and only if $\Forall x \allowbreak \in \allowbreak \R^{1 \times n} \backslash \{0\} \colon x A x^{\top} \allowbreak > 0$ and we say that $A$ is positive semi-definite (we write $A \succeq 0$) if and only if $\Forall x \allowbreak \in \allowbreak \R^{1 \times n} \colon x A x^{\top} \ge 0$.}
\begin{equation}\label{eqn:lemma:Schur_complement:D}
\textstyle A^{\top} = A, \qquad C^{\top} = C \succ 0, \qqandqq D =
\begin{psmallmatrix}
A & B \\[1ex]
B^{\top} & C
\end{psmallmatrix} \!\in \R^{(n + m) \times (n + m)}.
\end{equation}
Then
\begin{enumerate}[label=(\roman*)]
\item it holds that $D \succ 0$ if and only if $A - B C^{-1} B^{\top} \succ 0$,

\item it holds that $D \succeq 0$ if and only if $A - B C^{-1} B^{\top} \succeq 0$, and

\item it holds that $\rk(D) = \rk(C) + \rk(A - B C^{-1} B^{\top})$.
\end{enumerate}
\end{lemma}
\begin{cproof}{lemma:Schur_complement}
\Nobs that the assumption that $C \succ 0$ ensures that $C$ is invertible. Combining this with \cref{eqn:lemma:Schur_complement:D} implies that
\begin{equation}\label{eqn:lemma:Schur_complement:D_factorized}
\textstyle D =
\begin{psmallmatrix}
A & B \\[1ex]
B^{\top} & C
\end{psmallmatrix} =
\begin{psmallmatrix}
\idMatrix_n & B C^{-1} \\[1ex]
0 & \idMatrix_m
\end{psmallmatrix}
\begin{psmallmatrix}
A - B C^{-1} B^{\top} & 0 \\[1ex]
0 & C
\end{psmallmatrix}
\begin{psmallmatrix}
\idMatrix_n & B C^{-1} \\[1ex]
0 & \idMatrix_m
\end{psmallmatrix}^{\top}.
\end{equation}
Next \nobs that the fact that
\begin{equation}
\begin{psmallmatrix}
\idMatrix_n & B C^{-1} \\[1ex]
0 & \idMatrix_m
\end{psmallmatrix}
\begin{psmallmatrix}
\idMatrix_n & -B C^{-1} \\[1ex]
0 & \idMatrix_m
\end{psmallmatrix} = \idMatrix_{m + n}
\end{equation}
demonstrates that
\begin{equation}
\R^{(n + m) \times (n + m)} \ni \!
\begin{psmallmatrix}
\idMatrix_n & B C^{-1} \\[1ex]
0 & \idMatrix_m
\end{psmallmatrix}
\end{equation}
is invertible. This and \cref{eqn:lemma:Schur_complement:D_factorized} prove that it holds that $D \succ 0$ if and only if $A - B C^{-1} B^{\top} \succ 0$, it holds that $D \succeq 0$ if and only if $A - B C^{-1} B^{\top} \succeq 0$, and it holds that
\begin{equation}
\textstyle \rk(D) = \rk
\begin{psmallmatrix}
A - B C^{-1} B^{\top} & 0 \\[1ex]
0 & C
\end{psmallmatrix}
= \rk(C) + \rk(A - B C^{-1} B^{\top}).
\end{equation}
\end{cproof}

\cfclear
\begin{lemma}\label{lemma:MatrixFullRank_and_positivity}
Let $n \in \N$, $\fa, \fb \in \R$, $v \in \R^{1 \times n}$, $\fA = \allowbreak (\fA_{i, j})_{(i, j) \in \{1, 2, \dots, n\}^2} \allowbreak \in \R^{n \times n}$ satisfy
\begin{equation}\label{eqn:lemma:MatrixFullRank_and_positivity:A}
\textstyle \fa > \fb \geq 0, \quad v = (1, 1, \dots, 1), \qandq \Forall i, j \in \{1, 2, \dots, n\} \colon \fA_{i, j} = \fa \indicator{\{j\}}(i) + \fb \indicator{\R \backslash\{j\}}(i).
\end{equation}
Then it holds for all $c \in (0, \infty)$ that $(c \fA)^{\top} = c \fA \succ 0$, $\rk(c \fA) = n$, and
\begin{equation}\label{eqn:lemma:MatrixFullRank_and_positivity:A_inverse}
\textstyle (c \fA)^{-1} = c^{-1} (\fa - \fb)^{-1} [\idMatrix_n - \fb (\fa + (n - 1) \fb)^{-1} v^{\top} v].
\end{equation}
\end{lemma}
\begin{cproof}{lemma:MatrixFullRank_and_positivity}
\Nobs that \cref{eqn:lemma:MatrixFullRank_and_positivity:A} ensures that $\fA$ is symmetric and
\begin{equation}\label{eqn:lemma:MatrixFullRank_and_positivity:A_simple}
\textstyle \fA = (\fa - \fb) \idMatrix_n + \fb v^{\top} v.
\end{equation}
The fact that $\fa > \fb \geq 0$ \hence proves for all $x \in \R^{1 \times n} \backslash \{0\}$ that
\begin{multline}
\textstyle x \fA x^{\top} = x ((\fa - \fb) \idMatrix_n + \fb v^{\top} v) x^{\top} = (\fa - \fb) x x^{\top} + \fb x v^{\top} v x^{\top} \\
\textstyle = (\fa - \fb) \norm{x}^2 + \fb (x v^{\top}) (x v^{\top})^{\top} = (\fa - \fb) \norm{x}^2 + \fb \abs{x v^{\top}}^2 \ge (\fa - \fb) \norm{x}^2 > 0.
\end{multline}
This shows for all $c \in (0, \infty)$ that $c \fA \succ 0$ and $\rk(c \fA) = n$. Moreover, \nobs that \cref{eqn:lemma:MatrixFullRank_and_positivity:A_simple}, the fact that $\fa > \fb \geq 0$, and the fact that $v^{\top} v v^{\top} v = v^{\top} n v = n (v^{\top} v)$ show that
\begin{multline}
\textstyle \fA \bigl(\frac{\idMatrix_n}{\fa - \fb} - \frac{\fb v^{\top} v}{(\fa - \fb) (\fa + (n - 1) \fb)}\bigr) = ((\fa - \fb) \idMatrix_n + \fb v^{\top} v) \bigl(\frac{\idMatrix_n}{\fa - \fb} - \frac{\fb v^{\top} v}{(\fa - \fb) (\fa + (n - 1) \fb)}\bigr) = \idMatrix_n - \frac{\fb v^{\top} v}{\fa + (n - 1) \fb} \\
\textstyle + \frac{\fb v^{\top} v}{\fa - \fb} - \frac{\fb^2 v^{\top} v v^{\top} v}{(\fa - \fb) (\fa + (n - 1) \fb)} = \idMatrix_n + \bigl(\frac{\fb}{\fa - \fb} - \frac{\fb}{\fa + (n - 1) \fb} - \frac{\fb^2 n}{(\fa - \fb) (\fa + (n - 1) \fb)}\bigr) v^{\top} v = \idMatrix_n.
\end{multline}
This establishes \cref{eqn:lemma:MatrixFullRank_and_positivity:A_inverse}.
\end{cproof}

\begin{lemma}\label{lemma:rank&Hessian}
Assume \cref{setting0} and let $\theta \allowbreak = \allowbreak (\theta_1, \allowbreak \dots, \allowbreak \theta_{\dimension}) \allowbreak \in \cM$, $\fH \allowbreak = (\fH^{i, j})_{(i, j) \in \{1, 2, \dots, d + 3\}^2} \allowbreak \in \R^{(d + 3) \times (d + 3)}$ satisfy for all $i, j \in \{1, 2, \dots, d + 3\}$ that $\fH^{i, j} = (\scrb - \scra)^{- d} (\frac{\partial^2}{\partial \theta_{\pi(i)} \partial \theta_{\pi(j)}} \riskR)(\theta)$ (cf.\ \cref{lemma:Hessian_of_Risk_cM}). Then
\begin{enumerate}[label=(\roman *)]
\item
\label{item1:lemma:rank&Hessian} it holds for all $x  = (x_1, \dots, x_d) \in [\scra, \scrb]^d$ that 
\begin{equation}
\textstyle \functionSNN^{\theta}(x) = - \frac{(1 - \fq_{\theta})^{1/2}}{4 (1 + 3\fq_{\theta})^{1/2}} + \frac{1}{2 (1 - \fq_{\theta})^{3/2}(1 + 3 \fq_{\theta})^{1/2}} \max\{\frac{x_1 - \scra}{\scrb - \scra} - \fq_{\theta}, 0\},
\end{equation}

\item
\label{item2:lemma:rank&Hessian} it holds for all $j \in \{1, 2, \dots, d\}$ that 
\begin{equation}
\textstyle \int_{I_1^{\theta}} (\functionSNN^{\theta}(x) - f(x)) \, \d x = \int_{[\scra, \scrb]^d} (\functionSNN^{\theta}(x) - f(x)) \, \d x = \int_{I_1^{\theta}} \spro{e_j, x} (\functionSNN^{\theta}(x) - f(x)) \, \d x = 0,
\end{equation}

\item
\label{item3:lemma:rank&Hessian} it holds that $(\nabla \riskR)(\theta) = 0$,

\item
\label{item4:lemma:rank&Hessian} it holds that $\rk(\fH) = d + 1$,

\item
\label{item5:lemma:rank&Hessian} it holds that $\spectrum{\fH} \subseteq [0, \infty)$,

\item
\label{item6:lemma:rank&Hessian} it holds that $\rk((\Hs \riskR)(\theta)) = d + 1$, and

\item
\label{item7:lemma:rank&Hessian} it holds that $\spectrum{(\Hs \riskR)(\theta)} \subseteq [0, \infty)$\ifnocf.
\end{enumerate}
\cfout.
\end{lemma}
\begin{cproof}{lemma:rank&Hessian}
Throughout this proof let $\scrd \in \N$, $\fv, \fq \in \R$ satisfy $\scrd \allowbreak = \max\{d \allowbreak - 1, \allowbreak 1\}$, $\fv = \fv_1^{\theta}$, and $\fq = \fq_{\theta}$, let $A \allowbreak \in \R^{4 \times 4}$, $B \allowbreak \in \R^{4 \times \scrd}$, $C \in \R^{\scrd \times \scrd}$, $v \allowbreak \in \R^{1 \times \scrd}$ satisfy
\begin{equation}\label{eqn:lemma:rank&Hessian:A_B_C}
\begin{gathered}
\textstyle A = (\fH^{i, j})_{(i, j) \in \{1, 2, 3, 4\}^2}, \qquad B = [\indicator{(1, \infty)}(d)] \bigl(\fH^{i, j - \indicator{\{1\}}(d)}\bigr)_{(i, j) \in \{1, 2, 3, 4\} \times \{5, 6, \dots, \scrd + 4\}}, \\
\textstyle C = \bigl(\fH^{i - \indicator{\{1\}}(d), j - \indicator{\{1\}}(d)}\bigr)_{(i, j) \in \{5, 6, \dots, \scrd + 4\}^2}, \qqandqq v = \ones{\scrd} [\indicator{(1, \infty)}(d)],
\end{gathered}
\end{equation}
and let $g, \fg \in C([0, 1], \R)$ satisfy for all $x \in [0, 1]$ that
\begin{equation}\label{eq:lemma:rank&Hessian:g_fg}
\textstyle g(x) = \functionSNN^{\theta}((\scrb - \scra) x + \scra, \scra, \dots, \scra) \qqandqq \fg(x) = \ff(x, 0, \dots, 0).
\end{equation}
\Nobs that \cref{eqn:lemma:rank&Hessian:A_B_C} shows that
\begin{equation}\label{eqn:lemma:rank&Hessian:fH_A_B_C}
\textstyle \R^{(d + 3) \times (d + 3)} \ni \fH =
    \begin{cases}
    A & \colon d = 1 \\
    \begin{psmallmatrix}
    A & B \\[1ex]
    B^{\top} & C
    \end{psmallmatrix} & \colon d > 1.
    \end{cases}
\end{equation}
Next \nobs that \cref{0_eq:M} ensures for all $x = (x_1, \dots, x_d) \in [\scra, \scrb]^d$ that
\begin{equation}\label{eqn:lemma:rank&Hessian:realization}
\textstyle \functionSNN^{\theta} (x) = \fc^{\theta} + \fv_1^{\theta} \max\{\fb_1^{\theta} + \fw_{1, 1}^{\theta} x_1, 0\}.
\end{equation}
Furthermore, \nobs that \cref{0_eq:M} demonstrates that there exists $y \in (0, \infty)$ which satisfies
\begin{gather}\label{eqn:lemma:rank&Hessian:thetas0}
\textstyle \scra < \scrq_{\theta} < \scrb, \qquad 0 < \alpha < \fq < \beta < 1, \qquad \scrq_{\theta} = (\scrb - \scra) \fq + \scra, \qquad \sum_{j = 2}^d \abs{\fw_{1, j}^{\theta}} = 0, \\ \label{eqn:lemma:rank&Hessian:thetas}
\textstyle \fw_{1, 1}^{\theta} = \theta_1 = \frac{y}{\scrb - \scra}, \quad \fb_1^{\theta} = -\scrq_{\theta} \fw_{1, 1}^{\theta}, \quad \fv = \fv_1^{\theta} = \frac{(1 + 3 \fq)^{- 1/2}}{2 y (1 - \fq)^{3/2}}, \qandq \fc^{\theta} = - \frac{(1 - \fq)^{1/2}}{4 (1 + 3 \fq)^{1/2}}.
\end{gather}
\Hence that
\begin{equation}\label{eqn:lemma:rank&Hessian:theta_2H+1}
\textstyle \scrb - \scrq_{\theta} = (\scrb - \scra) (1 - \fq), \quad \fb_1^{\theta} = - \theta_1 [(\scrb - \scra) \fq + \scra], \qandq \theta_1 = \frac{(1 + 3 \fq)^{- 1/2}}{2 \fv (\scrb - \scra) (1 - \fq)^{3/2}}.
\end{equation}
Combining this with \cref{eqn:lemma:rank&Hessian:realization,eqn:lemma:rank&Hessian:thetas0,eqn:lemma:rank&Hessian:thetas} ensures for all $x = (x_1, \allowbreak \dots, \allowbreak x_d) \allowbreak \in \allowbreak [\scra, \scrb]^d$ that
\begin{equation}\label{eqn:lemma:rank&Hessian:realization_simple}
\textstyle \functionSNN^{\theta}(x) = - \frac{(1 - \fq)^{1/2}}{4 (1 + 3 \fq)^{1/2}} + \frac{1}{2 (1 - \fq)^{3/2}(1 + 3 \fq)^{1/2}} \max\{\frac{x_1 - \scra}{\scrb - \scra} - \fq, 0\}.
\end{equation}
This establishes \cref{item1:lemma:rank&Hessian}. \Nobs that \cref{eq:lemma:rank&Hessian:g_fg,eqn:lemma:rank&Hessian:realization_simple} assure for all $x \in [0, \allowbreak 1]$ that
\begin{equation}\label{eqn:lemma:rank&Hessian:realization_trans}
\textstyle g(x) = \functionSNN^{\theta}((\scrb - \scra) x + \scra, \scra, \dots, \scra) = - \frac{(1 - \fq)^{1/2}}{4 (1 + 3 \fq)^{1/2}} + \frac{\max\{x - \fq, 0\}}{2 (1 - \fq)^{3/2} (1 + 3 \fq)^{1/2}}.
\end{equation}
\Cref{item1:lemma:zero_derivative} in \cref{lemma:zero_derivative}, \cref{lemma:Hessian_of_Risk_cM}, \cref{0_eq:f}, \cref{eq:lemma:rank&Hessian:g_fg}, \cref{eqn:lemma:rank&Hessian:thetas0}, \cref{eqn:lemma:rank&Hessian:realization_simple}, the fact that $\alpha < \fq < \beta$, the integral transformation theorem, and Fubini's theorem \hence ensure for all $j \in \{1, 2, \dots, d\}$ that
\begin{equation}
\begin{split}
& \textstyle \frac{1}{(\scrb - \scra)^d} \int_{I_1^{\theta}} (\functionSNN^{\theta} (x) - f(x)) \, \d x = \frac{1}{(\scrb - \scra)^d} \int_{(\scrq_{\theta}, \scrb] \times [\scra, \scrb]^{d - 1}} (\functionSNN^{\theta} (x) - f(x)) \, \d x \\
& \textstyle = \int_{(\frac{\scrq_{\theta} - \scra}{\scrb - \scra}, 1] \times [0, 1]^{d - 1}} (g(\spro{e_1, x}) - \ff(x)) \, \d x = \int_{(\fq, 1] \times [0, 1]^{d - 1}} (g(\spro{e_1, x}) - \ff(x)) \, \d x \\
& \textstyle = \int_{\fq}^1 (g(y) - \fg(y)) \, \d y = 0,
\end{split}
\end{equation}
\begin{equation}
\begin{split}
& \textstyle \frac{1}{(\scrb - \scra)^d} \int_{[\scra, \scrb]^d \backslash I_1^{\theta}} (\functionSNN^{\theta} (x) - f(x)) \, \d x = \frac{1}{(\scrb - \scra)^d} \int_{[\scra, \scrq_{\theta}) \times [\scra, \scrb]^{d - 1}} (\functionSNN^{\theta} (x) - f(x)) \, \d x \\
& \textstyle = \int_{[0, \frac{\scrq_{\theta} - \scra}{\scrb - \scra}) \times [0, 1]^{d - 1}} (g(\spro{e_1, x}) - \ff(x)) \, \d x = \int_{[0, \fq) \times [0, 1]^{d - 1}} (g(\spro{e_1, x}) - \ff(x)) \, \d x \\
& \textstyle = \int_0^{\fq} (g(y) - \fg(y)) \, \d y = 0,
\end{split}
\end{equation}
\begin{equation}
\textstyle \int_{[\scra, \scrb]^d} (\functionSNN^{\theta} (x) - f(x)) \, \d x = \int_{I_1^{\theta}} (\functionSNN^{\theta} (x) - f(x)) \, \d x + \int_{[\scra, \scrb]^d \backslash I_1^{\theta}} (\functionSNN^{\theta} (x) - f(x)) \, \d x = 0,
\end{equation}
and
\begin{equation}
\begin{split}
& \textstyle \frac{1}{(\scrb - \scra)^d} \int_{I_1^{\theta}} \spro{e_j, x} (\functionSNN^{\theta} (x) - f(x)) \, \d x = \frac{1}{(\scrb - \scra)^d} \int_{(\scrq_{\theta}, \scrb] \times [\scra, \scrb]^{d - 1}} \spro{e_j, x} (\functionSNN^{\theta} (x) - f(x)) \, \d x \\
& \textstyle = \int_{(\fq, 1] \times [0, 1]^{d - 1}} ((\scrb - \scra) \spro{e_j, x} + \scra) (g(\spro{e_1, x}) - \ff(x)) \, \d x \\
& \textstyle = \indicator{\{1\}}(j) \int_{\fq}^1 ((\scrb - \scra) y + \scra) (g(y) - \fg(y)) \, \d y \\
& \textstyle \quad + \indicator{(1, d]}(j) \int_0^1 ((\scrb - \scra) y + \scra) \d y \int_{\fq}^1 (g(y) - \fg(y)) \, \d y = 0.
\end{split}
\end{equation}
This establishes \cref{item2:lemma:rank&Hessian}. Next \nobs that \cref{prop:limit:lr}, \cref{prop:loss:diff:vc}, \cref{prop:loss:diff:wb}, \cref{item_fV_M:lemma:Hessian_of_Risk_cM,item_Hess:lemma:Hessian_of_Risk_cM} in \cref{lemma:Hessian_of_Risk_cM}, \cref{item2:lemma:rank&Hessian}, and \cref{0_eq:M} prove that $(\nabla \riskR)(\theta) = 0$. This establishes \cref{item3:lemma:rank&Hessian}. In the next step \nobs that \cref{0_eq:f}, \cref{eq:lemma:rank&Hessian:g_fg}, \cref{eqn:lemma:rank&Hessian:thetas0}, and \cref{eqn:lemma:rank&Hessian:realization_trans} show that
\begin{equation}
\textstyle g(\fq) = - \frac{(1 - \fq)^{1/2}}{4 (1 + 3 \fq)^{1/2}} \qqandqq \fg(\fq) = \frac{3 \fq^2 - 1}{4 (1 - \fq)^{1/2} (1 + 3 \fq)^{3/2}}.
\end{equation}
This and \cref{eqn:lemma:rank&Hessian:theta_2H+1} imply that
\begin{multline}\label{eqn:lemma:rank&Hessian:g_minus_fg}
\textstyle \frac{g(\fq) - \fg(\fq)}{\theta_1} = - \frac{(1 - \fq)^{1/2}}{4 \theta_1 (1 + 3 \fq)^{1/2}} - \frac{3 \fq^2 - 1}{4 \theta_1 (1 - \fq)^{1/2} (1 + 3 \fq)^{3/2}} \\
\textstyle = \frac{- (1 - \fq) (1 + 3 \fq) - (3 \fq^2 - 1)}{4 \theta_1 (1 - \fq)^{1/2} (1 + 3 \fq)^{3/2}} = - \frac{2 \fq (2 \fv) (\scrb - \scra) (1 - \fq)^{3/2} (1 + 3 \fq)^{1/2}}{4 (1 - \fq)^{1/2} (1 + 3 \fq)^{3/2}} = - \frac{\fq \fv (\scrb - \scra) (1 - \fq)}{1 + 3 \fq}.
\end{multline}
Moreover, \nobs that \cref{lemma:Hessian_of_Risk_cM}, \cref{item1:lemma:rank&Hessian}, \cref{0_eq:f}, \cref{0_eq:M}, \cref{eq:lemma:rank&Hessian:g_fg}, \cref{eqn:lemma:rank&Hessian:thetas0}, \cref{eqn:lemma:rank&Hessian:thetas}, and Fubini's theorem ensure for all $y \allowbreak = \allowbreak (y_1, \allowbreak \dots, \allowbreak y_d) \allowbreak \in \partial I_1^{\theta} \allowbreak = \{\scrq_{\theta}\} \allowbreak \times [\scra, \allowbreak \scrb]^{d - 1}$, $j, \allowbreak \fj \allowbreak \in \{1, \allowbreak 2, \allowbreak \dots, \allowbreak d\}$ with $j \allowbreak \le \fj$ that
\begin{gather}\label{eqn:lemma:rank&Hessian:N_F_q}
\textstyle y_1 = \scrq_{\theta}, \qquad \functionSNN^{\theta}(y) = \functionSNN^{\theta}(\scrq_{\theta}, \scra, \dots, \scra) = g(\fq), \qquad f(y) = \allowbreak f(\scrq_{\theta}, \allowbreak \scra, \allowbreak \dots, \allowbreak \scra) \allowbreak = \fg(\fq), \\
\label{eqn:lemma:rank&Hessian:partial_I_1}
\textstyle \int_{\partial I_1^{\theta}} \d \sigma(x) = (\scrb - \scra)^{d - 1}, \quad
\int_{\partial I_1^{\theta}} \spro{e_j, x} \d \sigma(x) =
    \begin{cases}
    (\scrb - \scra)^{d - 1} \scrq_{\theta} & \colon 1 = j, \\
    (\scrb - \scra)^{d - 1} [\frac{\scrb + \scra}{2}] & \colon 1 < j,
    \end{cases} \\
\label{eqn:lemma:rank&Hessian:I_1_w}
\textstyle \int_{I_1^{\theta}} \spro{e_j, x} \d x = \int_{(\scrq_{\theta}, \scrb] \times [\scra, \scrb]^{d - 1}} \spro{e_j, x} \d x =
\begin{cases}
(\scrb - \scra)^{d - 1} \frac{\scrb^2 - [\scrq_{\theta}]^2}{2} & \colon 1 = j,\\
(\scrb - \scra)^{d - 1} \frac{(\scrb + \scra)(\scrb - \scrq_{\theta})}{2} & \colon 1 < j,
\end{cases} \\
\label{eqn:lemma:rank&Hessian:partial_I_1_w_w}
\textstyle \int_{\partial I_1^{\theta}} \spro{e_j, x} \spro{e_{\fj}, x} \d \sigma(x) =
\begin{cases}
(\scrb - \scra)^{d - 1} [\scrq_{\theta}]^2 & \colon 1 = j = \fj, \\
(\scrb - \scra)^{d - 1} [\frac{\scrq_{\theta} (\scrb + \scra)}{2}] & \colon 1 = j < \fj, \\
(\scrb - \scra)^{d - 1} [\frac{\scrb^2 + \scrb \scra + \scra^2}{3}] & \colon 1 < j = \fj, \\
(\scrb - \scra)^{d - 1} [\frac{\scrb + \scra}{2}]^2 & \colon 1 < j < \fj,
\end{cases} \\
\label{eqn:lemma:rank&Hessian:I_1_w_w}
\textstyle \int_{I_1^{\theta}} \spro{e_j, x} \spro{e_{\fj}, x} \d x =
\begin{cases}
(\scrb - \scra)^{d - 1} \frac{\scrb^3 - [\scrq_{\theta}]^3}{3} & \colon 1 = j = \fj, \\
(\scrb - \scra)^{d - 1} [\frac{(\scrb^2 - [\scrq_{\theta}]^2) (\scrb + \scra)}{4}] & \colon 1 = j < \fj, \\
(\scrb - \scra)^{d - 1} [\frac{(\scrb - \scrq_{\theta})(\scrb^2 + \scrb \scra + \scra^2)}{3}] & \colon 1 < j = \fj, \\
(\scrb - \scra)^{d - 1} [\frac{\scrb + \scra}{2}]^2 (\scrb - \scrq_{\theta}) & \colon 1 < j < \fj,
\end{cases} \\
\label{eqn:lemma:rank&Hessian:I_1_v}
\textstyle \int_{I_1^{\theta}} (\fb_1^{\theta} + \spro{\fw_1^{\theta}, x}) \d x = \int_{(\scrq_{\theta}, \scrb] \times [\scra, \scrb]^{d - 1}} \theta_1 (\spro{e_1, x} - \scrq_{\theta}) \d x = \theta_1 (\scrb - \scra)^{d - 1} \frac{(\scrb - \scrq_{\theta})^2}{2}, \\
\label{eqn:lemma:rank&Hessian:I_1_v_v}
\textstyle \int_{I_1^{\theta}} (\fb_1^{\theta} + \spro{\fw_1^{\theta}, x})^2 \d x = \int_{(\scrq_{\theta}, \scrb] \times [\scra, \scrb]^{d - 1}} [\theta_1]^2 (\spro{e_1, x} - \scrq_{\theta})^2 \d x = [\theta_1]^2 (\scrb - \scra)^{d - 1} \frac{(\scrb - \scrq_{\theta})^3}{3}, \\
\label{eqn:lemma:rank&Hessian:I_1_w_v}
\textstyle \andqq \int_{I_1^{\theta}} \spro{e_j, x} (\fb_1^{\theta} + \spro{\fw_1^{\theta}, x}) \d x =
\begin{cases}
\theta_1 (\scrb - \scra)^{d - 1} \frac{(\scrb - \scrq_{\theta})^2 (2 \scrb + \scrq_{\theta})}{6} & \colon 1 = j,\\
\theta_1 (\scrb - \scra)^{d - 1} \frac{(\scrb - \scrq_{\theta})^2 (\scrb + \scra)}{4} & \colon 1 < j.
\end{cases}
\end{gather}
Combining \cref{eqn:lemma:rank&Hessian:thetas0,eqn:lemma:rank&Hessian:thetas,eqn:lemma:rank&Hessian:g_minus_fg} \hence proves for all $j, \allowbreak \fj \allowbreak \in \{1, \allowbreak 2, \allowbreak \dots, \allowbreak d\}$ with $j \allowbreak \le \fj$ that
\begin{equation}
\textstyle \int_{\partial I_1^{\theta}} \frac{(\functionSNN^{\theta}(x) - f(x))}{\norm{\fw_1^{\theta}}} \d \sigma(x) = \frac{g(\fq) - \fg(\fq)}{\theta_1} \int_{\partial I_1^{\theta}} \d \sigma(x) = - \frac{\fq \fv (\scrb - \scra)^d (1 - \fq)}{1 + 3 \fq},
\end{equation}
\begin{equation}
\textstyle \int_{\partial I_1^{\theta}} \frac{\spro{e_j, x} (\functionSNN^{\theta}(x) - f(x))}{\norm{\fw_1^{\theta}}} \d \sigma(x) = \frac{g(\fq) - \fg(\fq)}{\theta_1}\! \int_{\partial I_1^{\theta}} \spro{e_j, x} \d \sigma(x) =
    \begin{cases}
    - \frac{\fq \fv (\scrb - \scra)^d \scrq_{\theta} (1 - \fq)}{1 + 3 \fq} \!\!\!& \colon 1 = j, \\
    - \frac{\fq \fv (\scrb - \scra)^d (\scrb + \scra) (1 - \fq)}{2 (1 + 3 \fq)} \!\!\!& \colon 1 < j,
    \end{cases}
\end{equation}
\begin{multline}
\textstyle \andqq \int_{\partial I_1^{\theta}} \frac{\spro{e_j, x} \spro{e_{\fj}, x} (\functionSNN^{\theta}(x) - f(x))}{\norm{\fw_1^{\theta}}} \d \sigma(x) = \frac{g(\fq) - \fg(\fq)}{\theta_1} \int_{\partial I_1^{\theta}} \spro{e_j, x} \spro{e_{\fj}, x} \d \sigma(x) \\
\textstyle =
    \begin{cases}
    - \frac{\fq \fv (\scrb - \scra)^d [\scrq_{\theta}]^2 (1 - \fq)}{1 + 3 \fq} & \colon 1 = j = \fj, \\
    - \frac{\fq \fv (\scrb - \scra)^d (\scrb + \scra) \scrq_{\theta} (1 - \fq)}{2 (1 + 3 \fq)} & \colon 1 = j < \fj, \\
    - \frac{\fq \fv (\scrb - \scra)^d (\scrb^2 + \scrb \scra + \scra^2) (1 - \fq)}{3 (1 + 3 \fq)} & \colon 1 < j = \fj,\\
    - \frac{\fq \fv (\scrb - \scra)^d (\scrb + \scra)^2 (1 - \fq)}{4 (1 + 3 \fq)} & \colon 1 < j < \fj.
    \end{cases}
\end{multline}
\cref{prop:limit:lr}, \cref{prop:loss:diff:vc}, \cref{prop:loss:diff:wb}, \cref{lemma:Hessian_of_Risk_cM}, \cref{item2:lemma:rank&Hessian}, \cref{0_eq:M}, \cref{eqn:lemma:rank&Hessian:thetas0,eqn:lemma:rank&Hessian:thetas,eqn:lemma:rank&Hessian:theta_2H+1,eqn:lemma:rank&Hessian:N_F_q,eqn:lemma:rank&Hessian:partial_I_1,eqn:lemma:rank&Hessian:I_1_w,eqn:lemma:rank&Hessian:partial_I_1_w_w,eqn:lemma:rank&Hessian:I_1_w_w,eqn:lemma:rank&Hessian:I_1_v,eqn:lemma:rank&Hessian:I_1_v_v,eqn:lemma:rank&Hessian:I_1_w_v} \hence show that
\begin{equation}\label{eqn:2nd_derivative_ww_1}
\begin{split}
\textstyle (\frac{\partial^2}{\partial \theta_1^2} \riskR)(\theta) & \textstyle = 2 [\fv_1^{\theta}]^2 \int_{I_1^{\theta}} \spro{e_1, x}^2 \d x + 2 \fv_1^{\theta} \int_{\partial I_1^{\theta}} \frac{\spro{e_1, x}^2 (\functionSNN^{\theta}(x) - f(x))}{\norm{\fw_1^{\theta}}} \d \sigma(x) \\
& \textstyle = 2 \fv^2 (\scrb - \scra)^{d - 1} \frac{\scrb^3 - [\scrq_{\theta}]^3}{3} - \frac{2 \fq \fv^2 (\scrb - \scra)^d [\scrq_{\theta}]^2 (1 - \fq)}{1 + 3 \fq} \\
& \textstyle = 2 \fv^2 (\scrb - \scra)^d (1 - \fq) \frac{\scrb^2  + \scrb \scrq_{\theta} + [\scrq_{\theta}]^2}{3} - \frac{2 \fq \fv^2 (\scrb - \scra)^d [\scrq_{\theta}]^2 (1 - \fq)}{1 + 3 \fq} \\
& \textstyle = 2 \fv^2 (\scrb - \scra)^d (1 - \fq) \bigl[\frac{\scrb^2  + \scrb \scrq_{\theta} + [\scrq_{\theta}]^2}{3} - \frac{\fq [\scrq_{\theta}]^2}{1 + 3 \fq}\bigr] \\
& \textstyle = \frac{2 \fv^2 (\scrb - \scra)^d (1 - \fq)}{3 (1 + 3 \fq)} [(1 + 3 \fq) (\scrb^2 + \scrb \scrq_{\theta}) + [\scrq_{\theta}]^2] \\
& \textstyle = \frac{2 \fv^2 (\scrb - \scra)^d (1 - \fq)}{3 (1 + 3 \fq)} [(1 + 3 \fq) (\scrb^2 + \scrb [(\scrb - \scra) \fq + \scra]) + [(\scrb - \scra) \fq + \scra]^2] \\
& \textstyle = \frac{2 (\scrb - \scra)^d [\scra^2 (1 - \fq)^2 + \scrb^2 (1 + 2 \fq)^2 + \scrb \scra (1 + 4 \fq - 5 \fq^2)]}{3 \fv^{-2} (1 - \fq)^{-1} (1 + 3 \fq)},
\end{split}
\end{equation}
\begin{equation}\label{eqn:2nd_derivative_bb}
\begin{split}
\textstyle (\frac{\partial^2}{\partial \theta_{d \width + 1}^2} \riskR)(\theta) & \textstyle = 2 [\fv_1^{\theta}]^2 \int_{I_1^{\theta}} \d x + 2 \fv_1^{\theta} \int_{\partial I_1^{\theta}} \frac{\functionSNN^{\theta}(x) - f(x)}{\norm{\fw_1^{\theta}}} \d \sigma(x) \\
& \textstyle = 2 \fv^2 (\scrb - \scrq_{\theta}) (\scrb - \scra)^{d - 1} - \frac{2 \fq \fv^2 (\scrb - \scra)^d (1 - \fq)}{1 + 3 \fq} \\
& \textstyle = 2 \fv^2 (\scrb - \scra)^d (1 - \fq) - \frac{2 \fq \fv^2 (\scrb - \scra)^d (1 - \fq)}{1 + 3 \fq} \\
& \textstyle = 2 \fv^2 (\scrb - \scra)^d (1 - \fq) (1 - \frac{\fq}{1 + 3 \fq}) = \frac{2 \fv^2 (\scrb - \scra)^d (1 - \fq) (1 + 2 \fq)}{1 + 3 \fq},
\end{split}
\end{equation}
\begin{equation}\label{eqn:2nd_derivative_vv}
\begin{split}
\textstyle (\frac{\partial^2}{\partial \theta_{d \width + \width + 1}^2} \riskR)(\theta) & \textstyle = 2 \int_{I_1^{\theta}} (\fb_1^{\theta} + \spro{\fw_1^{\theta}, x})^2 \d x \\
& \textstyle = \frac{2}{3} [\theta_1]^2 (\scrb - \scrq_{\theta})^3 (\scrb - \scra)^{d - 1} = \frac{2}{3} [\theta_1]^2 (\scrb - \scra)^{d + 2} (1 - \fq)^3 \\
& \textstyle = \frac{2}{3} (\scrb - \scra)^{d + 2} (1 - \fq)^3 \frac{(1 + 3 \fq)^{- 1}}{4 \fv^2 (\scrb - \scra)^2 (1 - \fq)^3} = \frac{(\scrb - \scra)^d}{6 \fv^2 (1 + 3 \fq)},
\end{split}
\end{equation}
\begin{equation}\label{eqn:2nd_derivative_cc}
\textstyle (\frac{\partial^2}{\partial \theta_{\dimension}^2} \riskR) (\theta) = 2 (\scrb - \scra)^d,
\end{equation}
\begin{equation}\label{eqn:2nd_derivative_bv}
\begin{split}
\textstyle (\frac{\partial^2}{\partial \theta_{d \width + 1} \partial \theta_{d \width + \width + 1}} \riskR)(\theta) & \textstyle = 2 \fv_1^{\theta} \int_{I_1^{\theta}} (\fb_1^{\theta} + \spro{\fw_1^{\theta}, x}) \d x + 2 \int_{I_1^{\theta}} (\functionSNN^{\theta}(x) - f(x)) \d x \\
& \textstyle = 2 \fv \int_{I_1^{\theta}} (\fb_1^{\theta} + \spro{\fw_1^{\theta}, x}) \d x = \fv \theta_1 (\scrb - \scrq_{\theta})^2 (\scrb - \scra)^{d - 1} \\
& \textstyle = \frac{(\scrb - \scra)^2 (1 - \fq)^2 (\scrb - \scra)^{d - 1}}{2 (\scrb - \scra) (1 - \fq)^{3/2} (1 + 3 \fq)^{1/2}} = \frac{(\scrb - \scra)^d (1 - \fq)^{1/2}}{2 (1 + 3 \fq)^{1/2}},
\end{split}
\end{equation}
\begin{equation}\label{eqn:2nd_derivative_bc}
\textstyle (\frac{\partial^2}{\partial \theta_{d \width + 1} \partial \theta_{\dimension}} \riskR)(\theta) = 2 \fv_1^{\theta} \int_{I_1^{\theta}} \d x = 2 \fv (\scrb - \scrq_{\theta}) (\scrb - \scra)^{d - 1} = 2 \fv (\scrb - \scra)^d (1 - \fq),
\end{equation}
\begin{equation}\label{eqn:2nd_derivative_vc}
\begin{split}
\textstyle (\frac{\partial^2}{\partial \theta_{d \width + \width + 1} \partial \theta_{\dimension}} \riskR)(\theta) & \textstyle = 2 \int_{I_1^{\theta}} (\fb_1^{\theta} + \spro{\fw_1^{\theta}, x}) \d x = \theta_1 (\scrb - \scrq_{\theta})^2 (\scrb - \scra)^{d - 1} \\
& \textstyle = \theta_1 (\scrb - \scra)^{d + 1} (1 - \fq)^2 = \frac{(\scrb - \scra)^{d + 1} (1 - \fq)^2 (1 + 3 \fq)^{- 1/2}}{2 \fv (\scrb - \scra) (1 - \fq)^{3/2}} = \frac{(\scrb - \scra)^d (1 - \fq)^{1/2}}{2 \fv (1 + 3 \fq)^{1/2}},
\end{split}
\end{equation}
\begin{equation}\label{eqn:2nd_derivative_wb_1}
\begin{split}
\textstyle (\frac{\partial^2}{\partial \theta_1 \partial \theta_{d \width + 1}} \riskR) (\theta) & \textstyle = 2 [\fv_1^{\theta}]^2 \int_{I_1^{\theta}} \spro{e_1, x} \d x + 2 \fv_1^{\theta} \int_{\partial I_1^{\theta}} \frac{\spro{e_1, x} (\functionSNN^{\theta}(x) - f(x))}{\norm{\fw_1^{\theta}}} \d \sigma(x) \\
& \textstyle = \fv^2 (\scrb - \scra)^{d - 1} (\scrb^2 - [\scrq_{\theta}]^2) - \frac{2 \fq \fv^2 (\scrb - \scra)^d \scrq_{\theta} (1 - \fq)}{1 + 3 \fq} \\
& \textstyle = \frac{\fv^2 (\scrb - \scra)^d (\scrb + \scrq_{\theta})}{(1 - \fq)^{-1}} - \frac{2 \fq \fv^2 (\scrb - \scra)^d \scrq_{\theta} (1 - \fq)}{1 + 3 \fq} = \frac{\fv^2 (\scrb - \scra)^d [(\scrb + \scrq_{\theta})(1 + 3 \fq) - 2 \fq \scrq_{\theta}]}{(1 - \fq)^{-1} (1 + 3 \fq)} \\
& \textstyle = \frac{\fv^2 [(\scrb + (\scrb - \scra) \fq + \scra)(1 + 3 \fq) - 2 \fq ((\scrb - \scra) \fq + \scra)]}{(\scrb - \scra)^{-d} (1 - \fq)^{-1} (1 + 3 \fq)} = \frac{(\scrb - \scra)^d [\scra (1 - \fq^2) + \scrb (1 + 4 \fq + \fq^2)]}{\fv^{-2} (1 - \fq)^{-1} (1 + 3 \fq)},
\end{split}
\end{equation}
\begin{equation}\label{eqn:2nd_derivative_wv_1}
\begin{split}
\textstyle (\frac{\partial^2}{\partial \theta_1 \partial \theta_{d \width + \width + 1}} \riskR) (\theta) & \textstyle = 2 \fv_1^{\theta} \int_{I_1^{\theta}} \spro{e_1, x} (\fb_1^{\theta} + \spro{\fw_1^{\theta}, x}) \d x + 2 \int_{I_1^{\theta}} \spro{e_1, x} (\functionSNN^{\theta}(x) - f(x)) \d x \\
& \textstyle = 2 \fv \int_{I_1^{\theta}} \spro{e_1, x} (\fb_1^{\theta} + \spro{\fw_1^{\theta}, x}) \d x = \frac{\fv \theta_1 (\scrb - \scra)^{d - 1} (\scrb - \scrq_{\theta})^2 (2 \scrb + \scrq_{\theta})}{3} \\
& \textstyle = \frac{\fv \theta_1 (\scrb - \scra)^{d + 1} (1 - \fq)^2 (2 \scrb + \scrq_{\theta})}{3} = \frac{(\scrb - \scra)^d (2 \scrb + \scrq_{\theta}) (1 - \fq)^{1/2}}{6 (1 + 3 \fq)^{1/2}} \\
& \textstyle = \frac{(\scrb - \scra)^d (2 \scrb + \scra + (\scrb - \scra) \fq) (1 - \fq)^{1/2}}{6 (1 + 3 \fq)^{1/2}} = \frac{(\scrb - \scra)^d (\scra (1 - \fq) + \scrb (2 + \fq))}{6 (1 - \fq)^{- 1/2} (1 + 3 \fq)^{1/2}},
\end{split}
\end{equation}
\begin{equation}\label{eqn:2nd_derivative_wc_1}
\begin{split}
\textstyle (\frac{\partial^2}{\partial \theta_1 \partial \theta_{\dimension}} \riskR) (\theta) & \textstyle = 2 \fv_1^{\theta} \int_{I_1^{\theta}} \spro{e_1, x} \d x = \fv (\scrb - \scra)^{d - 1} (\scrb^2 - [\scrq_{\theta}]^2) = \frac{(\scrb - \scra)^d (\scrb + \scrq_{\theta})}{\fv^{-1} (1 - \fq)^{-1}} \\
& \textstyle = \frac{(\scrb - \scra)^d (\scrb + (\scrb - \scra) \fq + \scra)}{\fv^{-1} (1 - \fq)^{-1}} = \frac{(\scrb - \scra)^d (\scra (1 - \fq) + \scrb (1 + \fq))}{\fv^{-1} (1 - \fq)^{-1}},
\end{split}
\end{equation}
\begin{equation}\label{eqn:2nd_derivative_the_rest}
\textstyle \andqq \bigwedge_{k \in \{1, 2, \dots, \dimension\} \backslash \{\pi(1), \pi(2), \dots, \pi(d + 3)\}} \bigwedge_{\ell \in \{1, 2, \dots, \dimension\}} \bigl[(\frac{\partial^2}{\partial \theta_k \partial \theta_{\ell}} \riskR)(\theta) = 0\bigr]
\end{equation}
and for all $j, \allowbreak \fj \allowbreak \in \{2, \allowbreak 3, \allowbreak \dots, \allowbreak d\}$ with $j \allowbreak \le \fj$ that
\begin{equation}\label{eqn:2nd_derivative_ww_1_j}
\begin{split}
\textstyle (\frac{\partial^2}{\partial \theta_1 \partial \theta_j} \riskR)(\theta) & \textstyle = 2 [\fv_1^{\theta}]^2 \int_{I_1^{\theta}} \spro{e_1, x} \spro{e_j, x} \d x + 2 \fv_1^{\theta} \int_{\partial I_1^{\theta}} \frac{\spro{e_1, x} \spro{e_j, x} (\functionSNN^{\theta}(x) - f(x))}{\norm{\fw_1^{\theta}}} \d \sigma(x) \\
& \textstyle = \frac{\fv^2 (\scrb - \scra)^{d - 1} (\scrb^2 - [\scrq_{\theta}]^2) (\scrb + \scra)}{2} - \frac{\fq \fv^2 (\scrb - \scra)^d (\scrb + \scra) \scrq_{\theta} (1 - \fq)}{1 + 3 \fq} \\
& \textstyle = \frac{\fv^2 (\scrb + \scra) (\scrb + \scrq_{\theta})}{2 (\scrb - \scra)^{- d} (1 - \fq)^{-1}} - \frac{\fv^2 (\scrb + \scra) \fq \scrq_{\theta}}{(\scrb - \scra)^{- d} (1 - \fq)^{-1} (1 + 3 \fq)} = \frac{\fv^2 (\scrb + \scra) [(\scrb + \scrq_{\theta}) (1 + 3 \fq) - 2 \fq \scrq_{\theta}]}{2 (\scrb - \scra)^{- d} (1 - \fq)^{-1} (1 + 3 \fq)} \\
& \textstyle = \frac{\fv^2 (\scrb + \scra) [(\scrb + (\scrb - \scra) \fq + \scra) (1 + 3 \fq) - 2 \fq ((\scrb - \scra) \fq + \scra)]}{2 (\scrb - \scra)^{- d} (1 - \fq)^{-1} (1 + 3 \fq)} = \frac{(\scrb - \scra)^d (\scrb + \scra) [\scra (1 - \fq^2) + \scrb (1 + 4 \fq + \fq^2)]}{2 \fv^{-2} (1 - \fq)^{-1} (1 + 3 \fq)},
\end{split}
\end{equation}
\begin{equation}\label{eqn:2nd_derivative_wb}
\begin{split}
\textstyle (\frac{\partial^2}{\partial \theta_j \partial \theta_{d \width + 1}} \riskR) (\theta) & \textstyle = 2 [\fv_1^{\theta}]^2 \int_{I_1^{\theta}} \spro{e_j, x} \d x + 2 \fv_1^{\theta} \int_{\partial I_1^{\theta}} \frac{\spro{e_j, x} (\functionSNN^{\theta}(x) - f(x))}{\norm{\fw_1^{\theta}}} \d \sigma(x) \\
& \textstyle = \fv^2 (\scrb - \scra)^{d - 1} (\scrb + \scra)(\scrb - \scrq_{\theta}) - \frac{\fv^2 (\scrb - \scra)^d (\scrb + \scra) (1 - \fq) \fq}{(1 + 3 \fq)} \\
& \textstyle = \frac{\fv^2 (\scrb - \scra)^d (\scrb + \scra)}{(1 - \fq)^{-1}} - \frac{\fv^2 (\scrb - \scra)^d (\scrb + \scra) \fq}{(1 - \fq)^{-1} (1 + 3 \fq)} = \frac{\fv^2 (\scrb - \scra)^d (\scrb + \scra) (1 + 2 \fq)}{(1 - \fq)^{-1} (1 + 3 \fq)},
\end{split}
\end{equation}
\begin{equation}\label{eqn:2nd_derivative_wv}
\begin{split}
\textstyle (\frac{\partial^2}{\partial \theta_j \partial \theta_{d \width + \width + 1}} \riskR) (\theta) & \textstyle = 2 \fv_1^{\theta} \int_{I_1^{\theta}} \spro{e_j, x} (\fb_1^{\theta} + \spro{\fw_1^{\theta}, x}) \d x + 2 \int_{I_1^{\theta}} \spro{e_j, x} (\functionSNN^{\theta}(x) - f(x)) \d x \\
& \textstyle = 2 \fv \int_{I_1^{\theta}} \spro{e_j, x} (\fb_1^{\theta} + \spro{\fw_1^{\theta}, x}) \d x = \frac{\fv \theta_1 (\scrb - \scra)^{d - 1} (\scrb - \scrq_{\theta})^2 (\scrb + \scra)}{2} \\
& \textstyle = \frac{\fv (\scrb - \scra)^{d + 1} (1 - \fq)^2 (\scrb + \scra)}{2} \cdot \frac{(1 + 3 \fq)^{- 1/2}}{2 \fv (\scrb - \scra) (1 - \fq)^{3/2}} = \frac{(\scrb - \scra)^d (\scrb + \scra) (1 - \fq)^{1/2}}{4 (1 + 3 \fq)^{1/2}},
\end{split}
\end{equation}
\begin{equation}\label{eqn:2nd_derivative_wc}
\begin{split}
\textstyle (\frac{\partial^2}{\partial \theta_j \partial \theta_{\dimension}} \riskR) (\theta) & \textstyle = 2 \fv_1^{\theta} \int_{I_1^{\theta}} \spro{e_j, x} \d x \\
& \textstyle = \fv (\scrb - \scra)^{d - 1} (\scrb + \scra)(\scrb - \scrq_{\theta}) = \fv (\scrb - \scra)^d (\scrb + \scra) (1 - \fq),
\end{split}
\end{equation}
and
\begin{equation}\label{eqn:2nd_derivative_ww}
\begin{split}
\textstyle (\frac{\partial^2}{\partial \theta_j \partial \theta_{\fj}} \riskR)(\theta) & \textstyle = 2 [\fv_1^{\theta}]^2 \int_{I_1^{\theta}} \spro{e_j, x} \spro{e_{\fj}, x} \d x + 2 \fv_1^{\theta} \int_{\partial I_1^{\theta}} \frac{\spro{e_j, x} \spro{e_{\fj}, x} (\functionSNN^{\theta}(x) - f(x))}{\norm{\fw_1^{\theta}}} \d \sigma(x) \\
& \textstyle =
    \begin{cases}
    \frac{2 \fv^2 (\scrb - \scra)^d (\scrb^2 + \scrb \scra + \scra^2) (1 + 2 \fq) (1 - \fq)}{3 (1 + 3 \fq)} & \colon 1 < j = \fj, \\
    \frac{\fv^2 (\scrb - \scra)^d (\scrb + \scra)^2 (1 + 2 \fq) (1 - \fq)}{2 (1 + 3 \fq)} & \colon 1 < j < \fj.
    \end{cases}
\end{split}
\end{equation}
This and \cref{eqn:lemma:rank&Hessian:A_B_C} demonstrate that
\begin{equation}\label{eqn:lemma:rank&Hessian:A_explicit}
\begin{aligned}[b]
A & \textstyle = (\fH^{i, j})_{(i, j) \in \{1, 2, 3, 4\}^2} = \frac{1}{(\scrb - \scra)^d} \bigl((\frac{\partial^2}{\partial \theta_{\pi(i)} \partial \theta_{\pi(j)}} \riskR) (\theta)\bigr)_{(i, j) \in \{1, 2, 3, 4\}^2} \\
& \textstyle = \!
\begin{psmallmatrix}
\frac{2 [\scra^2 (1 - \fq)^2 + \scrb^2 (1 + 2 \fq)^2 + \scrb \scra (1 + 4 \fq - 5 \fq^2)]}{3 \fv^{-2} (1 - \fq)^{-1} (1 + 3 \fq)} & \frac{\scra (1 - \fq^2) + \scrb (1 + 4 \fq + \fq^2)}{\fv^{-2} (1 - \fq)^{-1} (1 + 3 \fq)} & \frac{\scra (1 - \fq) + \scrb (2 + \fq)}{6 (1 - \fq)^{- 1/2} (1 + 3 \fq)^{1/2}} & \frac{\scra (1 - \fq) + \scrb (1 + \fq)}{\fv^{-1} (1 - \fq)^{-1}} \\[1ex]
\frac{\scra (1 - \fq^2) + \scrb (1 + 4 \fq + \fq^2)}{\fv^{-2} (1 - \fq)^{-1} (1 + 3 \fq)} & \frac{2 (1 - \fq) (1 + 2 \fq)}{\fv^{-2} (1 + 3 \fq)} & \frac{(1 - \fq)^{1/2}}{2 (1 + 3 \fq)^{1/2}} & 2 \fv (1 - \fq) \\[1ex]
\frac{\scra (1 - \fq) + \scrb (2 + \fq)}{6 (1 - \fq)^{- 1/2} (1 + 3 \fq)^{1/2}} & \frac{(1 - \fq)^{1/2}}{2 (1 + 3 \fq)^{1/2}} & \frac{1}{6 \fv^2 (1 + 3 \fq)} & \frac{(1 - \fq)^{1/2}}{2 \fv (1 + 3 \fq)^{1/2}} \\[1ex]
\frac{\scra (1 - \fq) + \scrb (1 + \fq)}{\fv^{-1} (1 - \fq)^{-1}} & 2 \fv (1 - \fq) & \frac{(1 - \fq)^{1/2}}{2 \fv (1 + 3 \fq)^{1/2}} & 2
\end{psmallmatrix} \!\in \R^{4 \times 4},
\end{aligned}
\end{equation}
\begin{align}
B & \textstyle = [\indicator{(1, \infty)}(d)] \bigl(\fH^{i, j - \indicator{\{1\}}(d)}\bigr)_{(i, j) \in \{1, 2, 3, 4\} \times \{5, 6, \dots, \scrd + 4\}} \notag \\
& \textstyle = \frac{\indicator{(1, \infty)}(d)}{(\scrb - \scra)^d} \bigl((\frac{\partial^2}{\partial \theta_{\pi(i)} \partial \theta_{\pi(j - \indicator{\{1\}}(d))}} \riskR) (\theta)\bigr)_{(i, j) \in \{1, 2, 3, 4\} \times \{5, 6, \dots, \scrd + 4\}} \notag \\ \label{eqn:lemma:rank&Hessian:B_explicit}
& \textstyle =\!
\begin{psmallmatrix}
\frac{(\scrb + \scra) [\scra (1 - \fq^2) + \scrb (1 + 4 \fq + \fq^2)]}{2 \fv^{-2} (1 - \fq)^{-1} (1 + 3 \fq)} \\[1ex]
\frac{\fv^2 (\scrb + \scra) (1 + 2 \fq)}{(1 - \fq)^{-1} (1 + 3 \fq)} \\[1ex]
\frac{(\scrb + \scra) (1 - \fq)^{1/2}}{4 (1 + 3 \fq)^{1/2}} \\[1ex]
\fv (\scrb + \scra) (1 - \fq)
\end{psmallmatrix} v = \fv (\scrb + \scra) (1 - \fq)
\begin{psmallmatrix}
\frac{\fv [\scra (1 - \fq^2) + \scrb (1 + 4 \fq + \fq^2)]}{2 (1 + 3 \fq)} \\[1ex]
\frac{\fv (1 + 2 \fq)}{1 + 3 \fq} \\[1ex]
\frac{(1 + 3 \fq)^{- 1/2}}{4 \fv (1 - \fq)^{1/2}} \\[1ex]
1
\end{psmallmatrix} v \in \R^{4 \times \scrd},
\end{align}
and
\begin{align}
C & \textstyle = \bigl(\fH^{i - \indicator{\{1\}}(d), j - \indicator{\{1\}}(d)}\bigr)_{(i, j) \in \{5, 6, \dots, \scrd + 4\}^2} \notag \\
& \textstyle = \frac{1}{(\scrb - \scra)^d} \bigl((\frac{\partial^2}{\partial \theta_{\pi(i - \indicator{\{1\}}(d))} \partial \theta_{\pi(j - \indicator{\{1\}}(d))}} \riskR) (\theta)\bigr)_{(i, j) \in \{5, 6, \dots, \scrd + 4\}^2} \notag \\
& \textstyle = 2 [\indicator{\{1\}}(d)] \ones{\scrd}^{\top} \ones{\scrd} + \indicator{(1, \infty)}(d) \!
\begin{psmallmatrix}
\frac{2 \fv^2 (\scrb^2 + \scrb \scra + \scra^2) (1 + 2 \fq)}{3 (1 - \fq)^{-1} (1 + 3 \fq)} & \frac{\fv^2 (\scrb + \scra)^2 (1 + 2 \fq)}{2 (1 - \fq)^{-1} (1 + 3 \fq)} & \cdots & \frac{\fv^2 (\scrb + \scra)^2 (1 + 2 \fq)}{2 (1 - \fq)^{-1} (1 + 3 \fq)} \\[1ex]
\frac{\fv^2 (\scrb + \scra)^2 (1 + 2 \fq)}{2 (1 - \fq)^{-1} (1 + 3 \fq)} & \frac{2 \fv^2 (\scrb^2 + \scrb \scra + \scra^2) (1 + 2 \fq)}{3 (1 - \fq)^{-1} (1 + 3 \fq)} & \cdots & \frac{\fv^2 (\scrb + \scra)^2 (1 + 2 \fq)}{2 (1 - \fq)^{-1} (1 + 3 \fq)} \\[1ex]
\vdots & \vdots & \ddots & \vdots \\[1ex]
\frac{\fv^2 (\scrb + \scra)^2 (1 + 2 \fq)}{2 (1 - \fq)^{-1} (1 + 3 \fq)} & \frac{\fv^2 (\scrb + \scra)^2 (1 + 2 \fq)}{2 (1 - \fq)^{-1} (1 + 3 \fq)} & \cdots & \frac{2 \fv^2 (\scrb^2 + \scrb \scra + \scra^2) (1 + 2 \fq)}{3 (1 - \fq)^{-1} (1 + 3 \fq)}
\end{psmallmatrix} \notag \\ \label{eqn:lemma:rank&Hessian:C_explicit}
& \textstyle = 2 [\indicator{\{1\}}(d)] \ones{\scrd}^{\top} \ones{\scrd} + \frac{2 \fv^2 (1 - \fq) (1 + 2 \fq) \indicator{(1, \infty)}(d)}{1 + 3 \fq}
\begin{psmallmatrix}
\frac{\scrb^2 + \scrb \scra + \scra^2}{3} & \frac{(\scrb + \scra)^2}{4} & \cdots & \frac{(\scrb + \scra)^2}{4} \\[1ex]
\frac{(\scrb + \scra)^2}{4} & \frac{\scrb^2 + \scrb \scra + \scra^2}{3} & \cdots & \frac{(\scrb + \scra)^2}{4} \\[1ex]
\vdots & \vdots & \ddots & \vdots \\[1ex]
\frac{(\scrb + \scra)^2}{4} & \frac{(\scrb + \scra)^2}{4} & \cdots & \frac{\scrb^2 + \scrb \scra + \scra^2}{3}
\end{psmallmatrix} \! \in \R^{\scrd \times \scrd}.
\end{align}
The fact that $\frac{2 \fv^2 (1 - \fq) (1 + 2 \fq)}{1 + 3 \fq} > 0$, the fact that
\begin{equation}
\textstyle \frac{\scrb^2 + \scrb \scra + \scra^2}{3} - \frac{(\scrb + \scra)^2}{4} = \frac{(\scrb - \scra)^2}{12} > 0,
\end{equation}
the fact that
\begin{equation}
\textstyle \bigl[\frac{2 \fv^2 (1 - \fq) (1 + 2 \fq)}{1 + 3 \fq}\bigr]^{-1} \bigl[\frac{\scrb^2 + \scrb \scra + \scra^2}{3} - \frac{(\scrb + \scra)^2}{4}\bigr]^{-1} = \frac{6 (1 + 3 \fq)}{\fv^2 (\scrb - \scra)^2 (1 - \fq) (1 + 2 \fq)},
\end{equation}
the fact that
\begin{equation}
\textstyle \frac{(\scrb + \scra)^2}{4} \bigl[\frac{\scrb^2 + \scrb \scra + \scra^2}{3} + \frac{(\scrd - 1) (\scrb + \scra)^2}{4}\bigr]^{-1} = \frac{3 (\scrb + \scra)^2}{3 (\scrd - 1) (\scrb + \scra)^2 + 4 (\scrb^2 + \scrb \scra + \scra^2)} = \frac{3 (\scrb + \scra)^2}{3 \scrd (\scrb + \scra)^2 + (\scrb - \scra)^2},
\end{equation}
and \cref{lemma:MatrixFullRank_and_positivity} (applied in the case $d = 1$ with $n \with 1$, $\fa \with 2$, $\fb \with 0$, $v \with 1$, $\fA \with 2$ and in the case $d > 1$ with $n \with \scrd$, $\fa \with 3^{- 1} (\scrb^2 + \scrb \scra + \scra^2)$, $\fb \with 4^{- 1} (\scrb + \scra)^2$, $v \with \ones{\scrd}$, $\fA \with \frac{(1 + 3 \fq) C}{2 \fv^2 (1 - \fq) (1 + 2 \fq)}$ in the notation of \cref{lemma:MatrixFullRank_and_positivity}) \hence assure that
\begin{gather}\label{eqn:lemma:rank&Hessian:C_rank}
\textstyle \rk(C) = \scrd, \qquad C^{\top} = C \succ 0, \\ \label{eqn:lemma:rank&Hessian:C_inv}
\textstyle \andqq C^{-1} = 2^{-1} [\indicator{\{1\}}(d)] \ones{\scrd}^{\top} \ones{\scrd} + \frac{6 (1 + 3 \fq) \indicator{(1, \infty)}(d)}{\fv^2 (\scrb - \scra)^2 (1 - \fq) (1 + 2 \fq)} \bigl[\idMatrix_{\scrd} - \frac{3 (\scrb + \scra)^2 \ones{\scrd}^{\top} \ones{\scrd}}{3 \scrd (\scrb + \scra)^2 + (\scrb - \scra)^2}\bigr].
\end{gather}
Moreover, \nobs that \cref{eqn:lemma:rank&Hessian:A_B_C} and the fact that $\scrd = \max\{d - 1, 1\}$ demonstrate that
\begin{equation}
\textstyle v \ones{\scrd}^{\top} \ones{\scrd} v^{\top} = v v^{\top} v v^{\top}, \qquad v v^{\top} = \indicator{(1, \infty)}(d) \ones{\scrd} \ones{\scrd}^{\top} = \indicator{(1, \infty)}(d) \scrd = d - 1,
\end{equation}
\begin{multline}
\textstyle \andqq v \bigl[\idMatrix_{\scrd} - \frac{3 (\scrb + \scra)^2 \ones{\scrd}^{\top} \ones{\scrd}}{3 \scrd (\scrb + \scra)^2 + (\scrb - \scra)^2}\bigr] v^{\top} = v v^{\top} - \frac{3 (\scrb + \scra)^2 v v^{\top} v v^{\top}}{3 \scrd (\scrb + \scra)^2 + (\scrb - \scra)^2} \\
\textstyle = (d - 1) - \frac{3 (\scrb + \scra)^2 (d - 1)^2}{3 (d - 1) (\scrb + \scra)^2 + (\scrb - \scra)^2} = \frac{(d - 1) (\scrb - \scra)^2}{3 (d - 1) (\scrb + \scra)^2 + (\scrb - \scra)^2}.
\end{multline}
This, \cref{eqn:lemma:rank&Hessian:A_B_C}, \cref{eqn:lemma:rank&Hessian:C_inv}, and the fact that $\indicator{\{1\}}(d) \indicator{(1, \infty)}(d) = 0$ imply for all $Y \in \cup_{k \in \N} \R^{k \times 1}$ that
\begin{multline}\label{eqn:lemma:rank&Hessian:Y_k}
\textstyle (Y v) C^{-1} (Y v)^{\top} = \frac{6 (1 + 3 \fq) \indicator{(1, \infty)}(d)}{\fv^2 (\scrb - \scra)^2 (1 - \fq) (1 + 2 \fq)} (Y v) \bigl[\idMatrix_{\scrd} - \frac{3 (\scrb + \scra)^2 \ones{\scrd}^{\top} \ones{\scrd}}{3 \scrd (\scrb + \scra)^2 + (\scrb - \scra)^2}\bigr] (Y v)^{\top} \\
\textstyle = \frac{6 (1 - \fq)^{-1} (1 + 3 \fq)}{\fv^2 (\scrb - \scra)^2 (1 + 2 \fq)} \cdot \frac{(d - 1) (\scrb - \scra)^2 Y Y^{\top}}{3 (d - 1) (\scrb + \scra)^2 + (\scrb - \scra)^2} = \frac{6 (d - 1) (1 - \fq)^{-1} (1 + 3 \fq) Y Y^{\top}}{\fv^2 (1 + 2 \fq) (3 (d - 1) (\scrb + \scra)^2 + (\scrb - \scra)^2)}.
\end{multline}
In the following let $P = (P_1, \dots, P_4) \in \R^{4 \times 1}$, $Q = (Q_1, Q_2) \in \R^{2 \times 1}$, $S \in \R$ satisfy
\begin{equation}\label{eqn:lemma:rank&Hessian:P_Q}
\textstyle P =
\begin{psmallmatrix}
\frac{\fv [\scra (1 - \fq^2) + \scrb (1 + 4 \fq + \fq^2)]}{2 (1 + 3 \fq)} \\[1ex]
\frac{\fv (1 + 2 \fq)}{1 + 3 \fq} \\[1ex]
\frac{(1 + 3 \fq)^{- 1/2}}{4 \fv (1 - \fq)^{1/2}} \\[1ex]
1
\end{psmallmatrix}\!, \qquad
Q =
\begin{psmallmatrix}
\frac{(1 + 3 \fq)^{- 1/2}}{4 \fv (1 - \fq)^{1/2}} \\[1ex]
1
\end{psmallmatrix}\!,
\qqandqq S = \fv (\scrb + \scra) (1 - \fq).
\end{equation}
\Nobs that \cref{eqn:lemma:rank&Hessian:Y_k}, \cref{eqn:lemma:rank&Hessian:P_Q}, and the fact that
\begin{multline}\label{eqn:lemma:rank&Hessian:fracs_temp0}
\textstyle \frac{6 (d - 1) (1 - \fq)^{-1} (1 + 3 \fq)}{\fv^2 (1 + 2 \fq) (3 (d - 1) (\scrb + \scra)^2 + (\scrb - \scra)^2)} \cdot [\fv (\scrb + \scra) (1 - \fq)]^2 = \frac{6 (d - 1) (\scrb + \scra)^2 (1 - \fq) (1 + 3 \fq)}{(1 + 2 \fq) (3 (d - 1) (\scrb + \scra)^2 + (\scrb - \scra)^2)} \\
\textstyle = \frac{2 (1 - \fq) (1 + 3 \fq)}{1 + 2 \fq} \cdot \frac{3 (d - 1) (\scrb + \scra)^2}{3 (d - 1) (\scrb + \scra)^2 + (\scrb - \scra)^2} = \frac{2 (1 - \fq) (1 + 3 \fq)}{1 + 2 \fq} \bigl[1 - \frac{(\scrb - \scra)^2}{3 (d - 1) (\scrb + \scra)^2 + (\scrb - \scra)^2}\bigr]
\end{multline}
demonstrate that
\begin{equation}
\textstyle Q Q^{\top} = \!
\begin{psmallmatrix}
\frac{(1 + 3 \fq)^{- 1/2}}{4 \fv (1 - \fq)^{1/2}} \\[1ex]
1
\end{psmallmatrix}
\begin{psmallmatrix}
\frac{(1 + 3 \fq)^{- 1/2}}{4 \fv (1 - \fq)^{1/2}} \\[1ex]
1
\end{psmallmatrix}^{\top}\! =\!
\begin{psmallmatrix}
\frac{(1 + 3 \fq)^{- 1}}{16 \fv^2 (1 - \fq)} & \frac{(1 + 3 \fq)^{- 1/2}}{4 \fv (1 - \fq)^{1/2}} \\[0.5ex]
\frac{(1 + 3 \fq)^{- 1/2}}{4 \fv (1 - \fq)^{1/2}} & 1
\end{psmallmatrix}\!,
\end{equation}
\begin{equation}
\begin{split}
\textstyle (S Q v) C^{-1} (S Q v)^{\top} & \textstyle = \frac{6 (d - 1) (1 - \fq)^{-1} (1 + 3 \fq)}{\fv^2 (1 + 2 \fq) (3 (d - 1) (\scrb + \scra)^2 + (\scrb - \scra)^2)} S^2 Q Q^{\top} \\
& \textstyle = \frac{6 (d - 1) (1 - \fq)^{-1} (1 + 3 \fq)}{\fv^2 (1 + 2 \fq) (3 (d - 1) (\scrb + \scra)^2 + (\scrb - \scra)^2)} \cdot [\fv (\scrb + \scra) (1 - \fq)]^2 Q Q^{\top} \\
& \textstyle = \frac{2 (1 - \fq) (1 + 3 \fq)}{1 + 2 \fq} \bigl[1 - \frac{(\scrb - \scra)^2}{3 (d - 1) (\scrb + \scra)^2 + (\scrb - \scra)^2}\bigr] Q Q^{\top} \in \R^{2 \times 2},
\end{split}
\end{equation}
and $Q_1 - P_3 = 0 = Q_2 - P_4$. Combining this with \cref{eqn:lemma:rank&Hessian:A_explicit}, the fact that
\begin{multline}\label{eqn:lemma:rank&Hessian:A_3_3_temp}
\textstyle \fH^{3, 3} - \frac{2 (1 - \fq) (1 + 3 \fq) [P_3]^2}{1 + 2 \fq} = \frac{1}{6 \fv^2 (1 + 3 \fq)} - \frac{2 (1 - \fq) (1 + 3 \fq)}{1 + 2 \fq} \cdot \frac{(1 + 3 \fq)^{- 1}}{16 \fv^2 (1 - \fq)} \\
\textstyle = \frac{1}{6 \fv^2 (1 + 3 \fq)} - \frac{1}{8 \fv^2 (1 + 2 \fq)} = \frac{1 - \fq}{24 \fv^2 (1 + 2 \fq) (1 + 3 \fq)},
\end{multline}
the fact that
\begin{multline}\label{eqn:lemma:rank&Hessian:A_3_4_temp}
\textstyle \fH^{3, 4} - \frac{2 (1 - \fq) (1 + 3 \fq) P_3 P_4}{1 + 2 \fq} = \frac{(1 - \fq)^{1/2}}{2 \fv (1 + 3 \fq)^{1/2}} - \frac{2 (1 - \fq) (1 + 3 \fq)}{1 + 2 \fq} \cdot \frac{(1 + 3 \fq)^{- 1/2}}{4 \fv (1 - \fq)^{1/2}} \\
\textstyle = \frac{(1 - \fq)^{1/2}}{2 \fv (1 + 3 \fq)^{1/2}} - \frac{(1 - \fq)^{1/2} (1 + 3\fq)^{1/2}}{2 \fv (1 + 2\fq)} = \frac{- \fq (1 - \fq)^{1/2}}{2 \fv (1 + 3 \fq)^{1/2} (1 + 2 \fq)},
\end{multline}
and the fact that
\begin{equation}\label{eqn:lemma:rank&Hessian:A_4_4_temp}
\textstyle \fH^{4, 4} - \frac{2 (1 - \fq) (1 + 3 \fq) [P_4]^2}{1 + 2 \fq} = 2 - \frac{2 (1 - \fq) (1 + 3 \fq)}{1 + 2 \fq} = \frac{2 (1 + 2 \fq - (1 - \fq)(1 + 3 \fq))}{1 + 2 \fq} = \frac{6 \fq^2}{1 + 2 \fq}
\end{equation}
implies that
\begin{align}
& \textstyle (\fH^{i, j})_{(i, j) \in \{3, 4\}^2} - (S Q v) C^{-1} (S Q v)^{\top} \notag\\
& \textstyle =\!
\begin{psmallmatrix}
\frac{1}{6 \fv^2 (1 + 3 \fq)} & \frac{(1 - \fq)^{1/2}}{2 \fv (1 + 3 \fq)^{1/2}} \\[0.5ex]
\frac{(1 - \fq)^{1/2}}{2 \fv (1 + 3 \fq)^{1/2}} & 2
\end{psmallmatrix}\! - \frac{2 (1 - \fq) (1 + 3 \fq)}{1 + 2 \fq} \bigl[1 - \frac{(\scrb - \scra)^2}{3 (d - 1) (\scrb + \scra)^2 + (\scrb - \scra)^2}\bigr] Q Q^{\top} \notag\\
& \textstyle =\!
\begin{psmallmatrix}
\frac{1}{6 \fv^2 (1 + 3 \fq)} & \frac{(1 - \fq)^{1/2}}{2 \fv (1 + 3 \fq)^{1/2}} \\[0.5ex]
\frac{(1 - \fq)^{1/2}}{2 \fv (1 + 3 \fq)^{1/2}} & 2
\end{psmallmatrix}\!
- \frac{2 (1 - \fq) (1 + 3 \fq)}{1 + 2 \fq} \!
\begin{psmallmatrix}
\frac{(1 + 3 \fq)^{- 1}}{16 \fv^2 (1 - \fq)} & \frac{(1 + 3 \fq)^{- 1/2}}{4 \fv (1 - \fq)^{1/2}} \\[0.5ex]
\frac{(1 + 3 \fq)^{- 1/2}}{4 \fv (1 - \fq)^{1/2}} & 1
\end{psmallmatrix}\! + \frac{2 (\scrb - \scra)^2 (1 - \fq) (1 + 3 \fq) Q Q^{\top}}{(1 + 2 \fq) (3 (d - 1) (\scrb + \scra)^2 + (\scrb - \scra)^2)} \notag\\
& \textstyle =\!
\begin{psmallmatrix}
\frac{1 - \fq}{24 \fv^2 (1 + 2 \fq) (1 + 3 \fq)} & \frac{- \fq (1 - \fq)^{1/2}}{2 \fv (1 + 3 \fq)^{1/2} (1 + 2 \fq)} \\[0.5ex]
\frac{- \fq (1 - \fq)^{1/2}}{2 \fv (1 + 3 \fq)^{1/2} (1 + 2 \fq)} & \frac{6 \fq^2}{1 + 2 \fq}
\end{psmallmatrix}\!
+ \frac{2 (\scrb - \scra)^2 (1 - \fq) (1 + 3 \fq) Q Q^{\top}}{(1 + 2 \fq) (3 (d - 1) (\scrb + \scra)^2 + (\scrb - \scra)^2)} \notag\\
& \textstyle = \frac{6 \fq^2}{1 + 2 \fq}
\begin{psmallmatrix}
\frac{1 - \fq}{144 \fv^2 \fq^2 (1 + 3 \fq)} & \frac{- (1 - \fq)^{1/2}}{12 \fv \fq (1 + 3 \fq)^{1/2}} \\[0.5ex]
\frac{- (1 - \fq)^{1/2}}{12 \fv \fq (1 + 3 \fq)^{1/2}} & 1
\end{psmallmatrix}\! + \frac{2 (\scrb - \scra)^2 (1 - \fq) (1 + 3 \fq) Q Q^{\top}}{(1 + 2 \fq) (3 (d - 1) (\scrb + \scra)^2 + (\scrb - \scra)^2)} \notag\\
& \textstyle = \frac{6 \fq^2}{1 + 2 \fq}
\Bigl(\begin{smallmatrix}
\frac{- (1 - \fq)^{1/2}}{12 \fv \fq (1 + 3 \fq)^{1/2}} \\[0.5ex]
1
\end{smallmatrix}\Bigr)\!
\Bigl(\begin{smallmatrix}
\frac{- (1 - \fq)^{1/2}}{12 \fv \fq (1 + 3 \fq)^{1/2}} \\[0.5ex]
1
\end{smallmatrix}\Bigr)^{\top} + \frac{2 (\scrb - \scra)^2 (1 - \fq) (1 + 3 \fq) Q Q^{\top}}{(1 + 2 \fq) (3 (d - 1) (\scrb + \scra)^2 + (\scrb - \scra)^2)} \in \R^{2 \times 2}.
\end{align}
The fact that $\scrb - \scra > 0 < \fv$, the fact that $0 < \alpha < \fq < \beta < 1$, the fact that $d \ge 1$, the fact that for all $\gamma, \delta \in (0, \infty)$ and a pair of linearly independent vectors $x, y \in \R^{2 \times 1}$ it holds that
\begin{equation}
\textstyle \R^{2 \times 2} \ni (\gamma x x^{\top} + \delta y y^{\top})^{\top} = \gamma x x^{\top} + \delta y y^{\top} \succ 0,
\end{equation}
and the fact that
\begin{equation}
\textstyle \det
\Bigl(\begin{smallmatrix}
\frac{(1 + 3 \fq)^{- 1/2}}{4 \fv (1 - \fq)^{1/2}} & \frac{- (1 - \fq)^{1/2}}{12 \fv \fq (1 + 3 \fq)^{1/2}} \\[0.5ex]
1 & 1
\end{smallmatrix}\Bigr)\! = \frac{(1 + 3 \fq)^{- 1/2}}{4 \fv (1 - \fq)^{1/2}} + \frac{(1 - \fq)^{1/2}}{12 \fv \fq (1 + 3 \fq)^{1/2}} = \frac{1 + 2 \fq}{12 \fv \fq (1 - \fq)^{1/2} (1 + 3 \fq)^{1/2}} > 0
\end{equation}
\hence prove that
\begin{multline}
\textstyle \R^{2 \times 2} \ni (\fH^{i, j})_{(i, j) \in \{3, 4\}^2} - (S Q v) C^{-1} (S Q v)^{\top} \\
\textstyle = \frac{6 \fq^2}{1 + 2 \fq}
\Bigl(\begin{smallmatrix}
\frac{- (1 - \fq)^{1/2}}{12 \fv \fq (1 + 3 \fq)^{1/2}} \\[0.5ex]
1
\end{smallmatrix}\Bigr)\!
\Bigl(\begin{smallmatrix}
\frac{- (1 - \fq)^{1/2}}{12 \fv \fq (1 + 3 \fq)^{1/2}} \\[0.5ex]
1
\end{smallmatrix}\Bigr)^{\top} + \frac{2 (\scrb - \scra)^2 (1 - \fq) (1 + 3 \fq) Q Q^{\top}}{(1 + 2 \fq) (3 (d - 1) (\scrb + \scra)^2 + (\scrb - \scra)^2)} \succ 0.
\end{multline}
Combining this, \cref{eqn:2nd_derivative_the_rest}, \cref{eqn:lemma:rank&Hessian:B_explicit}, \cref{eqn:lemma:rank&Hessian:C_explicit}, \cref{eqn:lemma:rank&Hessian:C_rank}, \cref{eqn:lemma:rank&Hessian:P_Q}, and the fact that
\begin{equation}
\textstyle \R^{2 \times \scrd} \ni [\indicator{(1, \infty)}(d)] (\fH^{i, j - \indicator{\{1\}}(d)})_{(i, j) \in \{3, 4\} \times \{5, 6, \dots, \scrd + 4\}} = S Q v
\end{equation}
with \cref{lemma:Schur_complement} demonstrates that
\begin{equation}
\textstyle \R^{(d + 1) \times (d + 1)} \ni (\fH^{i, j})_{(i, j) \in \{3, 4, \dots, d + 3\}^2} \succ 0
\end{equation}
\begin{equation}\label{eqn:lemma:rank&Hessian:rank_fH_lower}
\textstyle \andqq \rk((\Hs \riskR) (\theta)) = \rk(\fH) \ge \rk((\fH^{i, j})_{(i, j) \in \{3, 4, \dots, d + 3\}^2}) = d + 1.
\end{equation}
Next \nobs that \cref{eqn:lemma:rank&Hessian:B_explicit}, \cref{eqn:lemma:rank&Hessian:Y_k}, \cref{eqn:lemma:rank&Hessian:P_Q}, and \cref{eqn:lemma:rank&Hessian:fracs_temp0} ensure that
\begin{equation}\label{eqn:lemma:rank&Hessian:P_PT}
\textstyle B = S P v, \quad P P^{\top} =\!
\begin{psmallmatrix}
\frac{[\scra (1 - \fq^2) + \scrb (1 + 4 \fq + \fq^2)]^2}{4 \fv^{-2} (1 + 3 \fq)^2} & \frac{\scra (1 - \fq^2) + \scrb (1 + 4 \fq + \fq^2)}{2 \fv^{-2} (1 + 2 \fq)^{-1} (1 + 3 \fq)^2} & \frac{\scra (1 - \fq^2) + \scrb (1 + 4 \fq + \fq^2)}{8 (1 - \fq)^{1/2} (1 + 3 \fq)^{3/2}} & \frac{\scra (1 - \fq^2) + \scrb (1 + 4 \fq + \fq^2)}{2 \fv^{-1} (1 + 3 \fq)} \\[1ex]
\frac{\scra (1 - \fq^2) + \scrb (1 + 4 \fq + \fq^2)}{2 \fv^{-2} (1 + 2 \fq)^{-1} (1 + 3 \fq)^2} & \frac{\fv^2 (1 + 2 \fq)^2}{(1 + 3 \fq)^2} & \frac{(1 + 2 \fq) (1 + 3 \fq)^{- 3/2}}{4 (1 - \fq)^{1/2}} & \frac{\fv (1 + 2 \fq)}{1 + 3 \fq} \\[1ex]
\frac{\scra (1 - \fq^2) + \scrb (1 + 4 \fq + \fq^2)}{8 (1 - \fq)^{1/2} (1 + 3 \fq)^{3/2}} & \frac{(1 + 2 \fq) (1 + 3 \fq)^{- 3/2}}{4 (1 - \fq)^{1/2}} & \frac{(1 + 3 \fq)^{- 1}}{16 \fv^2 (1 - \fq)} & \frac{(1 + 3 \fq)^{- 1/2}}{4 \fv (1 - \fq)^{1/2}} \\[1ex]
\frac{\scra (1 - \fq^2) + \scrb (1 + 4 \fq + \fq^2)}{2 \fv^{-1} (1 + 3 \fq)} & \frac{\fv (1 + 2 \fq)}{1 + 3 \fq} & \frac{(1 + 3 \fq)^{- 1/2}}{4 \fv (1 - \fq)^{1/2}} & 1
\end{psmallmatrix}\!,
\end{equation}
and
\begin{equation}\label{eqn:lemma:rank&Hessian:B_C_BT}
\begin{split}
\textstyle B C^{-1} B^{\top} & \textstyle = (S P v) C^{-1} (S P v)^{\top} = \frac{6 (d - 1) (1 - \fq)^{-1} (1 + 3 \fq)}{\fv^2 (1 + 2 \fq) (3 (d - 1) (\scrb + \scra)^2 + (\scrb - \scra)^2)} S^2 P P^{\top} \\
& \textstyle = \frac{6 (d - 1) (1 - \fq)^{-1} (1 + 3 \fq)}{\fv^2 (1 + 2 \fq) (3 (d - 1) (\scrb + \scra)^2 + (\scrb - \scra)^2)} \cdot [\fv (\scrb + \scra) (1 - \fq)]^2 P P^{\top} \\
& \textstyle = \frac{2 (1 - \fq) (1 + 3 \fq)}{1 + 2 \fq} \bigl[1 - \frac{(\scrb - \scra)^2}{3 (d - 1) (\scrb + \scra)^2 + (\scrb - \scra)^2}\bigr] P P^{\top} \\
& \textstyle = \frac{2 (1 - \fq) (1 + 3 \fq) P P^{\top}}{1 + 2 \fq} - \frac{2 (\scrb - \scra)^2 (1 - \fq) (1 + 3 \fq) P P^{\top}}{(1 + 2 \fq) (3 (d - 1) (\scrb + \scra)^2 + (\scrb - \scra)^2)} \in \R^{4 \times 4}.
\end{split}
\end{equation}
Combining this with \cref{eqn:lemma:rank&Hessian:A_explicit}, \cref{eqn:lemma:rank&Hessian:A_3_3_temp}, \cref{eqn:lemma:rank&Hessian:A_3_4_temp}, \cref{eqn:lemma:rank&Hessian:A_4_4_temp}, the fact that
\begin{equation}
\textstyle 4 (1 + 2 \fq)(1 - \fq)^2 - 3 (1 - \fq^2)^2 = (1 - \fq)^3 (1 + 3 \fq) = 4 (1 + 2 \fq)^3 - 3 (1 + 4 \fq + \fq^2)^2,
\end{equation}
the fact that
\begin{equation}
\textstyle 4 (1 + 2 \fq) (1 + 4 \fq - 5 \fq^2) - 6 (1 - \fq^2) (1 + 4 \fq + \fq^2) = - 2 (1 - \fq)^3 (1 + 3 \fq),
\end{equation}
the fact that
\begin{equation}
\textstyle - 2 (1 + 2 \fq) (1 - \fq) + 3 (1 - \fq^2) = (1 - \fq)^2 = 2 (1 + 2 \fq) (2 + \fq) - 3 (1 + 4 \fq + \fq^2),
\end{equation}
the fact that
\begin{equation}
\textstyle (1 + 2 \fq)(1 - \fq) - (1 - \fq^2) = \fq (1 - \fq) = - (1 + 2 \fq)(1 + \fq) + (1 + 4 \fq + \fq^2),
\end{equation}
the fact that
\begin{equation}
\begin{split}
& \textstyle \fH^{1, 1} - \frac{2 (1 - \fq) (1 + 3 \fq) [P_1]^2}{1 + 2 \fq} \\
& \textstyle = \frac{2 [\scra^2 (1 - \fq)^2 + \scrb^2 (1 + 2 \fq)^2 + \scrb \scra (1 + 4 \fq - 5 \fq^2)]}{3 \fv^{-2} (1 - \fq)^{-1} (1 + 3 \fq)} - \frac{2 (1 - \fq) (1 + 3 \fq)}{1 + 2 \fq} \cdot \frac{[\scra (1 - \fq^2) + \scrb (1 + 4 \fq + \fq^2)]^2}{4 \fv^{-2} (1 + 3 \fq)^2} \\
& \textstyle = \frac{\fv^2 (1 - \fq) (4 (1 + 2 \fq) [\scra^2 (1 - \fq)^2 + \scrb^2 (1 + 2 \fq)^2 + \scrb \scra (1 + 4 \fq - 5 \fq^2)] - 3 [\scra (1 - \fq^2) + \scrb (1 + 4 \fq + \fq^2)]^2)}{6 (1 + 2 \fq) (1 + 3 \fq)} \\
& \textstyle = \frac{\fv^2 (1 - \fq)}{6 (1 + 2 \fq) (1 + 3 \fq)} \bigl[\scra^2 (4 (1 + 2 \fq)(1 - \fq)^2 - 3 (1 - \fq^2)^2) \\
& \textstyle \quad + \scra \scrb (4 (1 + 2 \fq) (1 + 4 \fq - 5 \fq^2) - 6 (1 - \fq^2) (1 + 4 \fq + \fq^2)) \\
& \textstyle \quad + \scrb^2 (4 (1 + 2 \fq)^3 - 3 (1 + 4 \fq + \fq^2)^2)\bigr] \\
& \textstyle = \frac{\fv^2 (1 - \fq)}{6 (1 + 2 \fq) (1 + 3 \fq)} \cdot (1 - \fq)^3 (1 + 3 \fq) (\scrb - \scra)^2 = \frac{\fv^2 (\scrb - \scra)^2 (1 - \fq)^4}{6 (1 + 2 \fq)},
\end{split}
\end{equation}
the fact that
\begin{equation}
\textstyle \fH^{1, 2} - \frac{2 (1 - \fq) (1 + 3 \fq) P_1 P_2}{1 + 2 \fq} = \frac{\scra (1 - \fq^2) + \scrb (1 + 4 \fq + \fq^2)}{\fv^{-2} (1 - \fq)^{-1} (1 + 3 \fq)} - \frac{2 (1 - \fq) (1 + 3 \fq)}{1 + 2 \fq} \cdot \frac{\scra (1 - \fq^2) + \scrb (1 + 4 \fq + \fq^2)}{2 \fv^{-2} (1 + 2 \fq)^{-1} (1 + 3 \fq)^2} = 0,
\end{equation}
the fact that
\begin{equation}
\begin{split}
& \textstyle \fH^{1, 3} - \frac{2 (1 - \fq) (1 + 3 \fq) P_1 P_3}{1 + 2 \fq} = \frac{\scra (1 - \fq) + \scrb (2 + \fq)}{6 (1 - \fq)^{- 1/2} (1 + 3 \fq)^{1/2}} - \frac{2 (1 - \fq) (1 + 3 \fq)}{1 + 2 \fq} \cdot \frac{\scra (1 - \fq^2) + \scrb (1 + 4 \fq + \fq^2)}{8 (1 - \fq)^{1/2} (1 + 3 \fq)^{3/2}} \\
& \textstyle = \frac{2 (1 + 2 \fq) (\scra (1 - \fq) + \scrb (2 + \fq)) - 3 (\scra (1 - \fq^2) + \scrb (1 + 4 \fq + \fq^2))}{12 (1 + 2 \fq) (1 - \fq)^{-1/2} (1 + 3 \fq)^{1/2}} \\
& \textstyle = \frac{\scra (2 (1 + 2 \fq) (1 - \fq) - 3 (1 - \fq^2)) + \scrb (2 (1 + 2 \fq) (2 + \fq) - 3 (1 + 4 \fq + \fq^2))}{12 (1 + 2 \fq) (1 - \fq)^{-1/2} (1 + 3 \fq)^{1/2}} \\
& \textstyle = \frac{(1 - \fq)^2 (\scrb - \scra)}{12 (1 + 2 \fq) (1 - \fq)^{-1/2} (1 + 3 \fq)^{1/2}} = \frac{(\scrb - \scra) (1 - \fq)^{5/2}}{12 (1 + 2 \fq) (1 + 3 \fq)^{1/2}},
\end{split}
\end{equation}
the fact that
\begin{equation}
\begin{split}
& \textstyle \fH^{1, 4} - \frac{2 (1 - \fq) (1 + 3 \fq) P_1 P_4}{1 + 2 \fq} = \frac{\scra (1 - \fq) + \scrb (1 + \fq)}{\fv^{-1} (1 - \fq)^{-1}} - \frac{2 (1 - \fq) (1 + 3 \fq)}{1 + 2 \fq} \cdot \frac{\scra (1 - \fq^2) + \scrb (1 + 4 \fq + \fq^2)}{2 \fv^{-1} (1 + 3 \fq)} \\
& \textstyle = \frac{(1 + 2 \fq) (\scra (1 - \fq) + \scrb (1 + \fq)) - (\scra (1 - \fq^2) + \scrb (1 + 4 \fq + \fq^2))}{\fv^{-1} (1 - \fq)^{-1} (1 + 2 \fq)} = \frac{\scra ((1 + 2 \fq)(1 - \fq) - (1 - \fq^2)) + \scrb ((1 + 2 \fq)(1 + \fq) - (1 + 4 \fq + \fq^2))}{\fv^{-1} (1 - \fq)^{-1} (1 + 2 \fq)} \\
& \textstyle = - \frac{\fq (1 - \fq) (\scrb - \scra)}{\fv^{-1} (1 - \fq)^{-1} (1 + 2 \fq)} = - \frac{(\scrb - \scra) (1 - \fq)^2 \fq \fv}{1 + 2 \fq},
\end{split}
\end{equation}
the fact that
\begin{multline}
\textstyle \fH^{2, 2} - \frac{2 (1 - \fq) (1 + 3 \fq) [P_2]^2}{1 + 2 \fq} = \frac{2 (1 - \fq) (1 + 2 \fq)}{\fv^{-2} (1 + 3 \fq)} - \frac{2 (1 - \fq) (1 + 3 \fq)}{1 + 2 \fq} \cdot \frac{\fv^2 (1 + 2 \fq)^2}{(1 + 3 \fq)^2} \\
\textstyle = \frac{2 \fv^2 (1 - \fq) (1 + 2 \fq)}{1 + 3 \fq} - \frac{2 \fv^2 (1 - \fq) (1 + 2 \fq)}{1 + 3 \fq} = 0,
\end{multline}
the fact that
\begin{multline}
\textstyle \fH^{2, 3} - \frac{2 (1 - \fq) (1 + 3 \fq) P_2 P_3}{1 + 2 \fq} = \frac{(1 - \fq)^{1/2}}{2 (1 + 3 \fq)^{1/2}} - \frac{2 (1 - \fq) (1 + 3 \fq)}{1 + 2 \fq} \cdot \frac{(1 + 2 \fq) (1 + 3 \fq)^{- 3/2}}{4 (1 - \fq)^{1/2}} \\
\textstyle = \frac{(1 - \fq)^{1/2}}{2 (1 + 3 \fq)^{1/2}} - \frac{(1 - \fq)^{1/2}}{2 (1 + 3 \fq)^{1/2}} = 0,
\end{multline}
and the fact that
\begin{equation}
\textstyle \fH^{2, 4} - \frac{2 (1 - \fq) (1 + 3 \fq) P_2 P_4}{1 + 2 \fq} = 2 \fv (1 - \fq) - \frac{2 (1 - \fq) (1 + 3 \fq)}{1 + 2 \fq} \cdot \frac{\fv (1 + 2 \fq)}{1 + 3 \fq} = 2 \fv (1 - \fq) - 2 \fv (1 - \fq) = 0
\end{equation}
demonstrates that
\begin{align}
& \textstyle A - \frac{2 (1 - \fq) (1 + 3 \fq) P P^{\top}}{1 + 2 \fq} = \bigl(\fH^{i, j} - \frac{2 (1 - \fq) (1 + 3 \fq) P_i P_j}{1 + 2 \fq}\bigr)_{(i, j) \in \{1, 2, 3, 4\}^2} \notag\\
& \textstyle =\!
\begin{psmallmatrix}
\frac{2 [\scra^2 (1 - \fq)^2 + \scrb^2 (1 + 2 \fq)^2 + \scrb \scra (1 + 4 \fq - 5 \fq^2)]}{3 \fv^{-2} (1 - \fq)^{-1} (1 + 3 \fq)} & \frac{\scra (1 - \fq^2) + \scrb (1 + 4 \fq + \fq^2)}{\fv^{-2} (1 - \fq)^{-1} (1 + 3 \fq)} & \frac{\scra (1 - \fq) + \scrb (2 + \fq)}{6 (1 - \fq)^{- 1/2} (1 + 3 \fq)^{1/2}} & \frac{\scra (1 - \fq) + \scrb (1 + \fq)}{\fv^{-1} (1 - \fq)^{-1}} \\[1ex]
\frac{\scra (1 - \fq^2) + \scrb (1 + 4 \fq + \fq^2)}{\fv^{-2} (1 - \fq)^{-1} (1 + 3 \fq)} & \frac{2 (1 - \fq) (1 + 2 \fq)}{\fv^{-2} (1 + 3 \fq)} & \frac{(1 - \fq)^{1/2}}{2 (1 + 3 \fq)^{1/2}} & 2 \fv (1 - \fq) \\[1ex]
\frac{\scra (1 - \fq) + \scrb (2 + \fq)}{6 (1 - \fq)^{- 1/2} (1 + 3 \fq)^{1/2}} & \frac{(1 - \fq)^{1/2}}{2 (1 + 3 \fq)^{1/2}} & \frac{1}{6 \fv^2 (1 + 3 \fq)} & \frac{(1 - \fq)^{1/2}}{2 \fv (1 + 3 \fq)^{1/2}} \\[1ex]
\frac{\scra (1 - \fq) + \scrb (1 + \fq)}{\fv^{-1} (1 - \fq)^{-1}} & 2 \fv (1 - \fq) & \frac{(1 - \fq)^{1/2}}{2 \fv (1 + 3 \fq)^{1/2}} & 2
\end{psmallmatrix} \notag\\
& \textstyle \quad - \frac{2 (1 - \fq) (1 + 3 \fq)}{1 + 2 \fq}
\begin{psmallmatrix}
\frac{[\scra (1 - \fq^2) + \scrb (1 + 4 \fq + \fq^2)]^2}{4 \fv^{-2} (1 + 3 \fq)^2} & \frac{\scra (1 - \fq^2) + \scrb (1 + 4 \fq + \fq^2)}{2 \fv^{-2} (1 + 2 \fq)^{-1} (1 + 3 \fq)^2} & \frac{\scra (1 - \fq^2) + \scrb (1 + 4 \fq + \fq^2)}{8 (1 - \fq)^{1/2} (1 + 3 \fq)^{3/2}} & \frac{\scra (1 - \fq^2) + \scrb (1 + 4 \fq + \fq^2)}{2 \fv^{-1} (1 + 3 \fq)} \\[1ex]
\frac{\scra (1 - \fq^2) + \scrb (1 + 4 \fq + \fq^2)}{2 \fv^{-2} (1 + 2 \fq)^{-1} (1 + 3 \fq)^2} & \frac{\fv^2 (1 + 2 \fq)^2}{(1 + 3 \fq)^2} & \frac{(1 + 2 \fq) (1 + 3 \fq)^{- 3/2}}{4 (1 - \fq)^{1/2}} & \frac{\fv (1 + 2 \fq)}{1 + 3 \fq} \\[1ex]
\frac{\scra (1 - \fq^2) + \scrb (1 + 4 \fq + \fq^2)}{8 (1 - \fq)^{1/2} (1 + 3 \fq)^{3/2}} & \frac{(1 + 2 \fq) (1 + 3 \fq)^{- 3/2}}{4 (1 - \fq)^{1/2}} & \frac{(1 + 3 \fq)^{- 1}}{16 \fv^2 (1 - \fq)} & \frac{(1 + 3 \fq)^{- 1/2}}{4 \fv (1 - \fq)^{1/2}} \\[1ex]
\frac{\scra (1 - \fq^2) + \scrb (1 + 4 \fq + \fq^2)}{2 \fv^{-1} (1 + 3 \fq)} & \frac{\fv (1 + 2 \fq)}{1 + 3 \fq} & \frac{(1 + 3 \fq)^{- 1/2}}{4 \fv (1 - \fq)^{1/2}} & 1
\end{psmallmatrix} \notag\\
& \textstyle =\!
\begin{psmallmatrix}
\frac{\fv^2 (\scrb - \scra)^2 (1 - \fq)^4}{6 (1 + 2 \fq)} & 0 & \frac{(\scrb - \scra) (1 - \fq)^{5/2}}{12 (1 + 2 \fq) (1 + 3 \fq)^{1/2}} & - \frac{(\scrb - \scra) (1 - \fq)^2 \fq \fv}{1 + 2 \fq} \\[1ex]
0 & 0 & 0 & 0 \\[1ex]
\frac{(\scrb - \scra) (1 - \fq)^{5/2}}{12 (1 + 2 \fq) (1 + 3 \fq)^{1/2}} & 0 & \frac{1 - \fq}{24 \fv^2 (1 + 2 \fq) (1 + 3 \fq)} & \frac{- \fq (1 - \fq)^{1/2}}{2 \fv (1 + 3 \fq)^{1/2} (1 + 2 \fq)} \\[1ex]
- \frac{(\scrb - \scra) (1 - \fq)^2 \fq \fv}{1 + 2 \fq} & 0 & \frac{- \fq (1 - \fq)^{1/2}}{2 \fv (1 + 3 \fq)^{1/2} (1 + 2 \fq)} & \frac{6 \fq^2}{1 + 2 \fq}
\end{psmallmatrix} \notag\\
& \textstyle = \frac{6 \fq^2}{1 + 2 \fq}
\begin{psmallmatrix}
\frac{(\scrb - \scra)^2 (1 - \fq)^4}{36 \fv^{-2} \fq^2} & 0 & \frac{(\scrb - \scra) (1 - \fq)^{5/2}}{72 \fq^2 (1 + 3 \fq)^{1/2}} & - \frac{(\scrb - \scra) (1 - \fq)^2}{6 \fv^{-1} \fq} \\[1ex]
0 & 0 & 0 & 0 \\[1ex]
\frac{(\scrb - \scra) (1 - \fq)^{5/2}}{72 \fq^2 (1 + 3 \fq)^{1/2}} & 0 & \frac{1 - \fq}{144 \fv^2 \fq^2 (1 + 3 \fq)} & \frac{- (1 - \fq)^{1/2}}{12 \fv \fq (1 + 3 \fq)^{1/2}} \\[1ex]
- \frac{(\scrb - \scra) (1 - \fq)^2}{6 \fv^{-1} \fq} & 0 & \frac{- (1 - \fq)^{1/2}}{12 \fv \fq (1 + 3 \fq)^{1/2}} & 1
\end{psmallmatrix}\! = \frac{6 \fq^2}{1 + 2 \fq}
\begin{psmallmatrix}
- \frac{(\scrb - \scra) (1 - \fq)^2}{6 \fv^{-1} \fq} \\[1ex]
0 \\[1ex]
\frac{- (1 - \fq)^{1/2}}{12 \fv \fq (1 + 3 \fq)^{1/2}} \\[1ex]
1
\end{psmallmatrix}
\begin{psmallmatrix}
- \frac{(\scrb - \scra) (1 - \fq)^2}{6 \fv^{-1} \fq} \\[1ex]
0 \\[1ex]
\frac{- (1 - \fq)^{1/2}}{12 \fv \fq (1 + 3 \fq)^{1/2}} \\[1ex]
1
\end{psmallmatrix}^{\top}\! .
\end{align}
The fact that $\scrb - \scra > 0 < \fv$, the fact that $0 < \alpha < \fq < \beta < 1$, the fact that $d \ge 1$, the fact that for all $n \in \N$, $\gamma, \delta \in [0, \infty)$, $x, y \in \R^{n \times 1}$, $\fA, \fB \in \R^{n \times n}$ it holds that
\begin{equation}
\textstyle \R^{n \times n} \ni (\gamma x x^{\top} + \delta y y^{\top})^{\top} = \allowbreak \gamma x x^{\top} \allowbreak + \delta y y^{\top} \allowbreak \succeq 0, \qquad \rk(x x^{\top}) = 1 - \indicator{\{0\}}(x),
\end{equation}
and $\rk(\fA + \fB) \le \rk(\fA) + \rk(\fB)$, \cref{eqn:lemma:rank&Hessian:P_Q}, and \cref{eqn:lemma:rank&Hessian:B_C_BT} \hence imply that
\begin{multline}
\textstyle \R^{4 \times 4} \ni A - B C^{-1} B^{\top} = A - \frac{2 (1 - \fq) (1 + 3 \fq) P P^{\top}}{1 + 2 \fq} + \frac{2 (\scrb - \scra)^2 (1 - \fq) (1 + 3 \fq) P P^{\top}}{(1 + 2 \fq) (3 (d - 1) (\scrb + \scra)^2 + (\scrb - \scra)^2)} \\
\textstyle = \frac{6 \fq^2}{1 + 2 \fq}
\begin{psmallmatrix}
- \frac{(\scrb - \scra) (1 - \fq)^2}{6 \fv^{-1} \fq} \\[1ex]
0 \\[1ex]
\frac{- (1 - \fq)^{1/2}}{12 \fv \fq (1 + 3 \fq)^{1/2}} \\[1ex]
1
\end{psmallmatrix}
\begin{psmallmatrix}
- \frac{(\scrb - \scra) (1 - \fq)^2}{6 \fv^{-1} \fq} \\[1ex]
0 \\[1ex]
\frac{- (1 - \fq)^{1/2}}{12 \fv \fq (1 + 3 \fq)^{1/2}} \\[1ex]
1
\end{psmallmatrix}^{\top}\! + \frac{2 (\scrb - \scra)^2 (1 - \fq) (1 + 3 \fq) P P^{\top}}{(1 + 2 \fq) (3 (d - 1) (\scrb + \scra)^2 + (\scrb - \scra)^2)} \succeq 0
\end{multline}
and
\begin{equation}
\textstyle \rk(A - B C^{-1} B^{\top}) \le \rk(P P^{\top}) + \rk
\begin{psmallmatrix}\!
\begin{psmallmatrix}
- \frac{(\scrb - \scra) (1 - \fq)^2}{6 \fv^{-1} \fq} \\[1ex]
0 \\[1ex]
\frac{- (1 - \fq)^{1/2}}{12 \fv \fq (1 + 3 \fq)^{1/2}} \\[1ex]
1
\end{psmallmatrix}
\begin{psmallmatrix}
- \frac{(\scrb - \scra) (1 - \fq)^2}{6 \fv^{-1} \fq} \\[1ex]
0 \\[1ex]
\frac{- (1 - \fq)^{1/2}}{12 \fv \fq (1 + 3 \fq)^{1/2}} \\[1ex]
1
\end{psmallmatrix}^{\top}\!
\end{psmallmatrix} \! = 2.
\end{equation}
This, \cref{eqn:lemma:rank&Hessian:A_B_C}, \cref{eqn:lemma:rank&Hessian:fH_A_B_C}, \cref{eqn:lemma:rank&Hessian:C_rank}, and \cref{lemma:Schur_complement} ensure that
\begin{equation}\label{eqn:lemma:rank&Hessian:fH_spectrum}
\textstyle \R^{(d + 3) \times (d + 3)} \ni \fH = (\fH^{i, j})_{(i, j) \in \{1, 2, \dots, d + 3\}^2} \succeq 0, \qquad \spectrum{\fH} \subseteq [0, \infty),
\end{equation}
\begin{multline}\label{eqn:lemma:rank&Hessian:rank_fH_upper}
\textstyle \andq \rk(\fH) = \indicator{(1, \infty)}(d) \rk(C) + \rk(A - B C^{-1} B^{\top}) \\
\textstyle = d - 1 + \rk(A - B C^{-1} B^{\top}) \le d + 1.
\end{multline}
This establishes \cref{item5:lemma:rank&Hessian}. Furthermore, \nobs that \cref{eqn:2nd_derivative_the_rest}, \cref{eqn:lemma:rank&Hessian:rank_fH_lower}, \cref{eqn:lemma:rank&Hessian:fH_spectrum}, and \cref{eqn:lemma:rank&Hessian:rank_fH_upper} prove that
\begin{equation}
\textstyle \R^{\dimension \times \dimension} \ni (\Hs \riskR)(\theta) \succeq 0, \qquad \spectrum{(\Hs \riskR)(\theta)} \subseteq [0, \infty),
\end{equation}
\begin{equation}
\textstyle \andqq \rk((\Hs \riskR) (\theta)) = \rk(\fH) = d + 1.
\end{equation}
This establishes \cref{item4:lemma:rank&Hessian,item6:lemma:rank&Hessian,item7:lemma:rank&Hessian}.
\end{cproof}

\subsection{On a submanifold of local minimum points of the ANN parameter space}
\label{subsec:submanifold_of_local_minimum_points_of_ANN_parameter_space}

\cfclear
\begin{lemma}\label{lemma:uncountably_many_realizations_of_critical_points}
Assume \cref{setting0} \cfload. Then
\begin{enumerate}[label=(\roman *)]
\item
\label{item1:lemma:uncountably_many_realizations_of_critical_points} it holds that $\bfM$ is an uncountable set,

\item
\label{item2:lemma:uncountably_many_realizations_of_critical_points} it holds for all $\theta \in \cM$ that $\riskR$ is differentiable at $\theta$,

\item
\label{item3:lemma:uncountably_many_realizations_of_critical_points} it holds for all $\theta \in \cM$ that $(\nabla \riskR) (\theta) = 0$,

\item
\label{item4:lemma:uncountably_many_realizations_of_critical_points} it holds for all $\theta, \vartheta \in \cM$ that $\riskR(\theta) = \riskR(\vartheta)$, and

\item
\label{item5:lemma:uncountably_many_realizations_of_critical_points} it holds that
\begin{multline}
\textstyle \bfM = \bigl\{v \in C([\scra, \scrb]^d, \R) \colon \\
\textstyle \bigl[\Exists \theta \in \{\vartheta \in \cM \colon v = \functionSNN^{\vartheta}\}, \eps \in (0, \infty) \colon \riskR(\theta) = \inf_{\vartheta \in [-\eps, \eps]^{\dimension}} \riskR(\theta + \vartheta) \bigr]\bigr\}\ifnocf.
\end{multline}
\end{enumerate}
\cfout[.]
\end{lemma}
\begin{cproof}{lemma:uncountably_many_realizations_of_critical_points}
\Nobs that \cref{item1:lemma:rank&Hessian} in \cref{lemma:rank&Hessian} and the fact that 
\begin{equation}
\textstyle \{q \in \R \colon [\Exists \theta = (\theta_1, \dots, \theta_{\dimension}) \in \cM \colon q = \fq_{\theta}]\} = (\alpha, \beta)
\end{equation}
ensure that $\bfM$ is an uncountable set. This establishes \cref{item1:lemma:uncountably_many_realizations_of_critical_points}. Next \nobs that \cref{lemma:Hessian_of_Risk_cM} and \cref{item3:lemma:rank&Hessian} in \cref{lemma:rank&Hessian} show for all $\theta \in \cM$ that $\riskR$ is differentiable at $\theta$ and $(\nabla \riskR)(\theta) = 0$. This establishes \cref{item2:lemma:uncountably_many_realizations_of_critical_points,item3:lemma:uncountably_many_realizations_of_critical_points}. In addition, \nobs that the integral transformation theorem, \cref{0_eq:f}, and \cref{item1:lemma:rank&Hessian} in \cref{lemma:rank&Hessian} prove that for all $\theta \in \cM$ it holds that
\begin{equation}
\begin{split}
\textstyle \riskR(\theta) & \textstyle = \int_{[\scra, \scrb]^d} (\functionSNN^{\theta}(x) - f(x))^2 \, \d x \\
& \textstyle = (\scrb - \scra)^{d - 1} \int_{\scra}^{\scrb} (\functionSNN^{\theta}(t, \scra, \dots, \scra) - f(t, \scra, \dots, \scra))^2 \, \d t \\
& \textstyle = (\scrb - \scra)^d \int_0^1 (\functionSNN^{\theta} ((\scrb - \scra) t + \scra, \scra, \dots, \scra) - f((\scrb - \scra) t + \scra, \scra, \dots, \scra))^2 \, \d t \\
& \textstyle = (\scrb - \scra)^d \int_0^1 (-\frac{(1 - \fq_{\theta})^{1/2}}{4 (1 + 3 \fq_{\theta})^{1/2}} + \frac{\max\{t - \fq_{\theta}, 0\}}{2 (1 - \fq_{\theta})^{3/2} (1 + 3 \fq_{\theta})^{1/2}} - \ff(t, 0, \dots, 0))^2 \, \d t.
\end{split}
\end{equation}
Combining this and \cref{item2:lemma:zero_derivative} in \cref{lemma:zero_derivative} with the fact that for all $\theta \in \cM$ it holds that $\alpha < \fq_{\theta} < \beta$ demonstrates for all $\theta \in \cM$ that
\begin{equation}
\textstyle \riskR(\theta) = (\scrb - \scra)^d (\int_0^1 [\ff(x, 0, \dots, 0)]^2 \, \d x - \frac{1}{48}).
\end{equation}
\Hence for all $\theta, \vartheta \in \cM$ that $\riskR(\theta) = \riskR(\vartheta)$. This establishes \cref{item4:lemma:uncountably_many_realizations_of_critical_points}. Furthermore, \nobs that \cref{prop:class:crit:points}, \cref{lemma:set_of_critical_points_submanifold}, \cref{lemma:Hessian_of_Risk_cM}, and \cref{item3:lemma:rank&Hessian,item6:lemma:rank&Hessian,item7:lemma:rank&Hessian} in \cref{lemma:rank&Hessian} demonstrate for all $\theta \in \cM$ that $\theta$ is a \locmin point of $\R^{\dimension} \ni \vartheta \mapsto \riskR(\vartheta) \in \R$ \cfload. This establishes \cref{item5:lemma:uncountably_many_realizations_of_critical_points}.
\end{cproof}

\subsection{On infinitely many realization functions of non-global local minimum points}
\label{subsec:existence_of_inf_many_realizations}

\cfclear
\begin{corollary}\label{cor:uncountably_many_realizations_of_non_global_local_minimas_with_fixed_delta}
Assume \cref{setting0} and assume $\width > 1$. Then there exists $\delta \in (0, \infty)$ such that
\begin{multline}
\textstyle \bigl\{v \in C([\scra, \scrb]^d, \R) \colon \bigl[\Exists \theta \in \{\vartheta \in \R^{\dimension} \colon v = \functionSNN^{\vartheta}\}, \eps \in (0, \infty) \colon \\
\textstyle \riskR(\theta) = \inf_{\vartheta \in [-\eps, \eps]^{\dimension}} \riskR(\theta + \vartheta) > \delta + \inf_{\vartheta \in \R^{\dimension}} \riskR(\vartheta)\bigr]\bigr\}
\end{multline}
is an uncountable set \cfout.
\end{corollary}
\begin{cproof}{cor:uncountably_many_realizations_of_non_global_local_minimas_with_fixed_delta}
Throughout this proof let $p \in (\alpha, \beta)$ and for every $\theta \in \R^{\dimension}$ let $\fN^{\theta}, \fg \allowbreak \in C([0, \allowbreak 1], \allowbreak \R)$ satisfy for all $x \in [0, 1]$ that $\fN^{\theta}(x) = \functionSNN^{\theta}((\scrb - \scra)x + \scra, \scra, \dots, \scra)$ and $\fg(x) \allowbreak = \ff(x, \allowbreak 0, \allowbreak \dots, \allowbreak 0)$ \cfload. \Nobs that \cref{lemma:uncountably_many_realizations_of_critical_points} ensures that
\begin{enumerate}[label=(\roman*)]
\item it holds that $\bfM$ is an uncountable set and
\item it holds that
\begin{multline}\label{eqn:cor:uncountably_many_realizations_of_non_global_local_minimas:bfM}
\textstyle \bfM = \bigl\{v \in C([\scra, \scrb]^d, \R) \colon \\
\textstyle \bigl[\Exists \theta \in \{\vartheta \in \cM \colon v = \functionSNN^{\vartheta}\}, \eps \in (0, \infty) \colon \riskR(\theta) = \inf_{\vartheta \in [-\eps, \eps]^{\dimension}} \riskR(\theta + \vartheta)\bigr]\bigr\}.
\end{multline}
\end{enumerate}
Next \nobs that \cref{item1:lemma:rank&Hessian} in \cref{lemma:rank&Hessian} shows for all $\theta \allowbreak \in \cM$, $x \allowbreak \in [0, \allowbreak 1]$ that
\begin{equation}
\textstyle \fN^{\theta}(x) = \functionSNN^{\theta}((\scrb - \scra)x + \scra, \scra, \dots, \scra) = -\frac{(1 - \fq_{\theta})^{1/2}}{4 (1 + 3 \fq_{\theta})^{1/2}} + \frac{\max\{x - \fq_{\theta}, 0\}}{2 (1 - \fq_{\theta})^{3/2} (1 + 3 \fq_{\theta})^{1/2}}.
\end{equation}
Combining this with the integral transformation theorem assures for all $\theta \allowbreak \in \cM$ that
\begin{equation}
\begin{split}
\textstyle \riskR(\theta) & \textstyle = \int_{[\scra, \scrb]^d} (\functionSNN^{\theta}(x) - f(x))^2 \, \d x \\
& \textstyle = (\scrb - \scra)^{d - 1} \int_{\scra}^{\scrb} (\functionSNN^{\theta}(t, \scra, \dots, \scra) - f(t, \scra, \dots, \scra))^2 \, \d t \\
& \textstyle = (\scrb - \scra)^d \int_0^1 (\functionSNN^{\theta} ((\scrb - \scra) t + \scra, \scra, \dots, \scra) - f((\scrb - \scra) t + \scra, \scra, \dots, \scra))^2 \, \d t \\
& \textstyle = (\scrb - \scra)^d \int_0^1 (-\frac{(1 - \fq_{\theta})^{1/2}}{4 (1 + 3 \fq_{\theta})^{1/2}} + \frac{\max\{t - \fq_{\theta}, 0\}}{2 (1 - \fq_{\theta})^{3/2} (1 + 3 \fq_{\theta})^{1/2}} - \fg(t))^2 \, \d t.
\end{split}
\end{equation}
The fact that for all $\theta \in \cM$ it holds that $\alpha < \fq_{\theta} < \beta$, \cref{item4:lemma:uncountably_many_realizations_of_critical_points} in \cref{lemma:uncountably_many_realizations_of_critical_points}, and \cref{item2:lemma:zero_derivative} in \cref{lemma:zero_derivative}
hence demonstrate that for all $\theta \in \cM$ it holds that
\begin{equation}\label{eqn:cor:uncountably_many_realizations_of_non_global_local_minimas:realization_simple}
\textstyle \riskR(\theta) = (\scrb - \scra)^d \int_0^1 (-\frac{(1 - p)^{1/2}}{4 (1 + 3 p)^{1/2}} + \frac{\max\{x - p, 0\}}{2 (1 - p)^{3/2} (1 + 3 p)^{1/2}} - \fg(x))^2 \, \d x.
\end{equation}
Next \nobs that \cref{0_eq:f} and the assumption that $0 < \alpha < p < \beta < 1$ ensure that
\begin{equation}
\textstyle \fg(p) = \frac{3p^2 - 1}{4(1 - p)^{1/2} (1 + 3 p)^{3/2}} > - \frac{(1 - p)^{1/2}}{4 (1 + 3 p)^{1/2}}.
\end{equation}
The fact that $[0, 1] \ni x \mapsto \fg(x) = \ff(x, 0, \dots, 0) \in \R$ is continuous \hence implies that there exists $\eps \in (0, \infty)$ which satisfies for all $x \in (p - \eps, p + \eps)$ that $(p - \eps, p + \eps) \subseteq (\alpha, \beta)$ and
\begin{equation}\label{eqn:cor:uncountably_many_realizations_of_non_global_local_minimas:double_ineq}
\textstyle \fg(x) > -\frac{(1 - p)^{1/2}}{4 (1 + 3 p)^{1/2}} + \frac{x - p + \eps}{4 (1 - p)^{3/2} (1 + 3 p)^{1/2}} > -\frac{(1 - p)^{1/2}}{4 (1 + 3 p)^{1/2}} + \frac{\max\{x - p, 0\}}{2 (1 - p)^{3/2} (1 + 3 p)^{1/2}}.
\end{equation}
In the following let $\vartheta = (\vartheta_1, \dots, \vartheta_{\dimension}) \in \R^{\dimension}$ satisfy
\begin{gather}\label{eqn:cor:uncountably_many_realizations_of_non_global_local_minimas:v1}
\textstyle\vartheta_1 = \vartheta_{d + 1} = \frac{1}{\scrb - \scra}, \qquad \vartheta_{d \width + 1} = -\frac{\scra}{\scrb - \scra} - p + \eps, \qquad \vartheta_{d \width + 2} = -\frac{\scra}{\scrb - \scra} - p - \eps, \\\label{eqn:cor:uncountably_many_realizations_of_non_global_local_minimas:v2}
\textstyle  \vartheta_{d \width + \width + 1} = \vartheta_{d \width + \width + 2} = \frac{1}{4 (1 - p)^{3/2} (1 + 3 p)^{1/2}}, \qquad \vartheta_{\dimension} = - \frac{(1 - p)^{1/2}}{4 (1 + 3 p)^{1/2}},
\end{gather}
and $\Forall j \in \{1, 2, \dots, \dimension\} \backslash \{1, \allowbreak d + 1, \allowbreak d \width + 1, \allowbreak d \width + 2, \allowbreak d \width + \width + 1, \allowbreak d \width + \width + 2, \allowbreak \dimension\} \colon \vartheta_j = 0$. \Nobs that \cref{eqn:cor:uncountably_many_realizations_of_non_global_local_minimas:v1}, \cref{eqn:cor:uncountably_many_realizations_of_non_global_local_minimas:v2}, and the integral transformation theorem ensure for all $x \in \allowbreak [0, \allowbreak 1]$, $y = \allowbreak (y_1, \allowbreak \dots, \allowbreak y_d) \allowbreak \in [\scra, \allowbreak \scrb]^d$ that
\begin{equation}
\textstyle \functionSNN^{\vartheta} (y) = - \frac{(1 - p)^{1/2}}{4 (1 + 3 p)^{1/2}} + \frac{\max\{\frac{y_1 - \scra}{\scrb - \scra} - p + \eps, 0\}}{4 (1 - p)^{3/2} (1 + 3p)^{1/2}} + \frac{\max\{\frac{y_1 - \scra}{\scrb - \scra} - p - \eps, 0\}}{4 (1 - p)^{3/2} (1 + 3p)^{1/2}},
\end{equation}
\begin{equation}
\textstyle \fN^{\vartheta}(x) = \functionSNN^{\vartheta} ((\scrb - \scra) x + \scra, \scra, \dots, \scra) = \frac{- (1 - p)^{1/2}}{4 (1 + 3 p)^{1/2}} + \frac{\max\{x - p + \eps, 0\}}{4 (1 - p)^{3/2} (1 + 3p)^{1/2}} + \frac{\max\{x - p - \eps, 0\}}{4 (1 - p)^{3/2} (1 + 3p)^{1/2}},
\end{equation}
\begin{equation}
\textstyle \andqq \riskR(\vartheta) = \int_{[\scra, \scrb]^d} (\functionSNN^{\vartheta}(x) - f(x))^2 \, \d x = (\scrb - \scra)^d \int_0^1 (\fN^{\vartheta}(x) - \fg(x))^2 \, \d x.
\end{equation}
Combining this with \cref{eqn:cor:uncountably_many_realizations_of_non_global_local_minimas:realization_simple}, \cref{eqn:cor:uncountably_many_realizations_of_non_global_local_minimas:double_ineq}, and the fact that for all $\fa, \fb, \fc \in \R$ with $\fa > \fb > \fc$ it holds that $(\fc - \fa)^2 > (\fb - \fa)^2$ proves for all $\theta \in \cM$ that
\begin{align}
\textstyle \frac{1}{(\scrb - \scra)^d} \riskR(\theta) & \textstyle = \int_0^1 (-\frac{(1 - p)^{1/2}}{4 (1 + 3 p)^{1/2}} + \frac{\max\{x - p, 0\}}{2 (1 - p)^{3/2} (1 + 3 p)^{1/2}} - \fg(x))^2 \, \d x \notag \\
& \textstyle = \int_0^{p - \eps} (-\frac{(1 - p)^{1/2}}{4 (1 + 3 p)^{1/2}} - \fg(x))^2 \, \d x \notag \\
& \quad \textstyle + \int_{p - \eps}^{p + \eps} (-\frac{(1 - p)^{1/2}}{4 (1 + 3 p)^{1/2}} + \frac{\max\{x - p, 0\}}{2 (1 - p)^{3/2} (1 + 3 p)^{1/2}} - \fg(x))^2 \, \d x \notag \\
& \quad \textstyle + \int_{p + \eps}^1 (-\frac{(1 - p)^{1/2}}{4 (1 + 3 p)^{1/2}} + \frac{\max\{x - p, 0\}}{2 (1 - p)^{3/2} (1 + 3 p)^{1/2}} - \fg(x))^2 \, \d x \notag \\
& \textstyle > \int_0^{p - \eps} (-\frac{(1 - p)^{1/2}}{4 (1 + 3 p)^{1/2}} - \fg(x))^2 \, \d x \notag \\
& \quad \textstyle + \int_{p - \eps}^{p + \eps} (-\frac{(1 - p)^{1/2}}{4 (1 + 3 p)^{1/2}} + \frac{x - p + \eps}{4 (1 - p)^{3/2} (1 + 3 p)^{1/2}} - \fg(x))^2 \, \d x \notag \\
& \quad \textstyle + \int_{p + \eps}^1 (-\frac{(1 - p)^{1/2}}{4 (1 + 3 p)^{1/2}} + \frac{\max\{x - p, 0\}}{2 (1 - p)^{3/2} (1 + 3 p)^{1/2}} - \fg(x))^2 \, \d x \notag \\
& \textstyle = \int_0^1 (- \frac{(1 - p)^{1/2}}{4 (1 + 3 p)^{1/2}} + \frac{\max\{x - p + \eps, 0\}}{4 (1 - p)^{3/2} (1 + 3p)^{1/2}} + \frac{\max\{x - p - \eps, 0\}}{4 (1 - p)^{3/2} (1 + 3p)^{1/2}} - \fg(x))^2 \, \d x \notag \\
& \textstyle = \int_0^1 (\fN^{\vartheta}(x) - \fg(x))^2 \, \d x = \frac{1}{(\scrb - \scra)^d} \riskR(\vartheta) \ge \frac{1}{(\scrb - \scra)^d}  \inf_{\Theta \in \R^{\dimension}} \riskR(\Theta).
\end{align}
This, \cref{eqn:cor:uncountably_many_realizations_of_non_global_local_minimas:bfM}, and \cref{item4:lemma:uncountably_many_realizations_of_critical_points} in \cref{lemma:uncountably_many_realizations_of_critical_points} show that there exists $\delta \in (0, \infty)$ which satisfies
\begin{equation}\label{eqn:cor:uncountably_many_realizations_of_non_global_local_minimas:bfM_new}
\begin{split}
\textstyle \bfM & = \bigl\{v \in C([\scra, \scrb]^d, \R) \colon \bigl[ \Exists \theta \in \{\Theta \in \cM \colon v = \functionSNN^{\Theta}\}, \epsilon \in (0, \infty) \colon \\
& \quad \quad \quad \textstyle \riskR(\theta) = \inf_{\Theta \in [-\epsilon, \epsilon]^{\dimension}} \riskR(\theta + \Theta) > \delta + \inf_{\Theta \in \R^{\dimension}} \riskR(\Theta) \bigr] \bigr\} \\
& \subseteq \bigl\{v \in C([\scra, \scrb]^d, \R) \colon \bigl[ \Exists \theta \in \{\Theta \in \R^{\dimension} \colon v = \functionSNN^{\Theta}\}, \epsilon \in (0, \infty) \colon \\
& \quad \quad \quad \textstyle \riskR(\theta) = \inf_{\Theta \in [-\epsilon, \epsilon]^{\dimension}} \riskR(\theta + \Theta) > \delta + \inf_{\Theta \in \R^{\dimension}} \riskR(\Theta) \bigr] \bigr\}.
\end{split}
\end{equation}
The fact that $\bfM$ is an uncountable set \hence demonstrates that
\begin{multline}
\textstyle \bigl\{v \in C([\scra, \scrb]^d, \R) \colon \bigl[\Exists \theta \in \{\Theta \in \R^{\dimension} \colon v = \functionSNN^{\Theta}\}, \epsilon \in (0, \infty) \colon \\
\textstyle \riskR(\theta) = \inf_{\Theta \in [-\epsilon, \epsilon]^{\dimension}} \riskR(\theta + \Theta) > \delta + \inf_{\Theta \in \R^{\dimension}} \riskR(\Theta)\bigr]\bigr\}
\end{multline}
is an uncountable set.
\end{cproof}

\cfclear
\begin{corollary}\label{cor:uncountably_many_realizations_of_non_global_local_minimas_with_varying_delta}
Let $d \in \N$, $\delta, \scra \in \R$, $\scrb \in (\scra, \infty)$, for every $\width \in \N$ let $\dimension_{\width} \in \N$ satisfy $\dimension_{\width} = \allowbreak d \width + \allowbreak 2 \width \allowbreak + 1$, for every $\width \in \N$, $\theta = \allowbreak (\theta_1, \allowbreak \dots, \allowbreak \theta_{\dimension_{\width}}) \in \R^{\dimension_{\width}}$ let $\cN_{\width}^{\theta} \in C([\scra, \scrb]^d, \R)$ satisfy for all $x = \allowbreak (x_1, \allowbreak \dots, \allowbreak x_d) \in \allowbreak [\scra, \allowbreak \scrb]^d$ that
\begin{equation}\label{eqn:cor:uncountably_many_realizations_delta:realization}
\textstyle \cN_{\width}^{\theta}(x) = \theta_{\dimension_{\width}} + \sum_{i = 1}^{\width} \theta_{d \width + \width + i} \max\bigl\{\theta_{d \width + i} + \sum_{j = 1}^d \theta_{(i - 1)d + j} x_j, 0\bigr\},
\end{equation}
and for every $\width \in \N$, $f \in C([\scra, \scrb]^d, \R)$ let $\cR_{f, \width} \colon \R^{\dimension_{\width}} \to \R$ satisfy for all $\theta \in \R^{\dimension_{\width}}$ that
\begin{equation}\label{eqn:cor:uncountably_many_realizations_delta:risk}
\textstyle \cR_{f, \width}(\theta) = \allowbreak \int_{[\scra, \scrb]^d} \abs{\cN_{\width}^{\theta}(y) - f(y)}^2 \, \d y.
\end{equation}
Then there exists an uncountable set $\cF \subseteq \{g \in C([\scra, \scrb]^d, \R) \colon g \text{ is Lipschitz continuous}\}$ such that for all $\scrf \in \cF$, $\width \in \N \cap (1, \infty)$ it holds that
\begin{multline}\label{eqn:cor:uncountably_many_realizations_of_non_global_local_minimas_with_varying_delta:res}
\textstyle \bigl\{v \in C([\scra, \scrb]^d, \R) \colon \bigl[\Exists \theta \in \{\vartheta \in \R^{\dimension_{\width}} \colon v = \cN_{\width}^{\vartheta}\}, \eps \in (0, \infty) \colon \\
\textstyle \cR_{\scrf, \width}(\theta) = \inf_{\vartheta \in [-\eps, \eps]^{\dimension_{\width}}} \cR_{\scrf, \width}(\theta + \vartheta) > \delta + \inf_{\vartheta \in \R^{\dimension_{\width}}} \cR_{\scrf, \width}(\vartheta)\bigr]\bigr\}
\end{multline}
is an uncountable set.
\end{corollary}
\begin{cproof}{cor:uncountably_many_realizations_of_non_global_local_minimas_with_varying_delta}
Throughout this proof let $\alpha \in (0, 1)$, $\beta \in (\alpha, 1)$, $\ff \allowbreak \in C([0, \allowbreak 1]^d, \allowbreak \R)$, $f \allowbreak \in C([\scra, \allowbreak \scrb]^d, \allowbreak \R)$ satisfy for all $x = (x_1, \dots, x_d) \in [0, 1]^d$ that
\begin{equation}\label{eqn:cor:uncountably_many_realizations_of_non_global_local_minimas_with_varying_delta:target_fn}
\ff(x) = 
\begin{cases}
\frac{4(x_1 - \alpha) + 3 \alpha^2 - 1}{4 (1 - \alpha)^{1/2} (1 + 3 \alpha)^{3/2}} & \colon x_1 \in [0, \alpha] \\[1ex]
\frac{3 [x_1]^2 - 1}{4 (1 - x_1)^{1/2} (1 + 3 x_1)^{3/2}} & \colon x_1 \in (\alpha, \beta] \\[1ex]
\frac{12 \beta [x_1]^2 - (18 \beta^2 + 8 \beta - 2) x_1 + 3 \beta^4 + 10 \beta^2 - 1}{4 (1 - \beta)^{5/2} (1 + 3 \beta)^{3/2}} & \colon x_1 \in (\beta, 1]
\end{cases}
\end{equation}
and $f((\scrb - \scra) x_1 + \scra, \dots, (\scrb - \scra) x_d + \scra) = \ff(x)$, let $\scrr \in (0, \infty)$ satisfy that
\begin{multline}\label{eqn:cor:uncountably_many_realizations_of_non_global_local_minimas_with_varying_delta:scrr}
\textstyle \bigl\{v \in C([\scra, \scrb]^d, \R) \colon \bigl[\Exists \theta \in \{\vartheta \in \R^{\dimension_2} \colon v = \cN_{2}^{\vartheta}\}, \eps \in (0, \infty) \colon \\
\textstyle \cR_{f, 2}(\theta) = \inf_{\vartheta \in [-\eps, \eps]^{\dimension_2}} \cR_{f, 2}(\theta + \vartheta) > \scrr + \inf_{\vartheta \in \R^{\dimension_2}} \cR_{f, 2}(\vartheta)\bigr]\bigr\}
\end{multline}
is an uncountable set (cf.\ \cref{cor:uncountably_many_realizations_of_non_global_local_minimas_with_fixed_delta}), let $\scrf = (\scrf_{\kappa})_{\kappa \in \R} \colon \R \to C([\scra, \scrb]^d, \R)$ satisfy
\begin{equation}\label{eqn:cor:uncountably_many_realizations_of_non_global_local_minimas_with_varying_delta:scrf_kappa}
\Forall (\kappa, x) \in \R \times [\scra, \scrb]^d \colon \scrf_{\kappa}(x) =
\begin{cases}
[\frac{\kappa}{\scrr}]^{1/2} f(x) & \colon \kappa \ge 0 \\
f(x) & \colon \kappa < 0,
\end{cases}
\end{equation}
and let $\cF \subseteq C([\scra, \scrb]^d, \R)$ satisfy $\cF = \cup_{\kappa \in [\max\{1,\delta\}, \infty)} \{\scrf_{\kappa}\}$. \Nobs that the chain rule ensures that for all $x \in (0, 1)$ it holds that
\begin{equation}
\textstyle \bigl[ \frac{4(x - \alpha) + 3 \alpha^2 - 1}{4(1 - \alpha)^{1/2} (1 + 3\alpha)^{3/2}} \bigr]' = \frac{1}{(1 - \alpha)^{1/2} (1 + 3\alpha)^{3/2}},
\end{equation}
\begin{equation}
\begin{split}
\textstyle \bigl[ \frac{3x^2 - 1}{4 (1 - x)^{1/2} (1 + 3x)^{3/2}} \bigr]' & \textstyle = \frac{([3x^2 - 1]') [(1 - x)^{1/2} (1 + 3x)^{3/2}] - [3x^2 - 1]([(1 - x)^{1/2} (1 + 3x)^{3/2}]')}{4[(1 - x)^{1/2} (1 + 3x)^{3/2}]^2} \\
& \textstyle = \frac{6x [(1 - x)^{1/2} (1 + 3x)^{3/2}] - [3x^2 - 1][ - \frac{1}{2} (1 - x)^{-1/2} (1 + 3x)^{3/2} + (1 - x)^{1/2} \frac{9}{2} (1 + 3x)^{1/2}]}{4(1 - x) (1 + 3x)^{3}} \\
& \textstyle = \frac{6x [(1 - x) (1 + 3x)] - [3x^2 - 1][ - \frac{1}{2} (1 + 3x) + \frac{9}{2} (1 - x)]}{4(1 - x)^{3/2} (1 + 3x)^{5/2}} = \frac{6x [1 + 2x - 3x^2] - [3x^2 - 1][4 - 6x]}{4(1 - x)^{3/2} (1 + 3x)^{5/2}} \\
& \textstyle = \frac{[6x + 12x^2 - 18x^3] - [12x^2 - 18 x^3 - 4 + 6x]}{4(1 - x)^{3/2} (1 + 3x)^{5/2}} = \frac{1}{(1 - x)^{3/2} (1 + 3x)^{5/2}},
\end{split}
\end{equation}
and
\begin{equation}
\textstyle \bigl[ \frac{12 \beta x^2 - (18 \beta^2 + 8 \beta - 2)x + 3 \beta^4 + 10 \beta^2 - 1}{4(1 - \beta)^{5/2} (1 + 3 \beta)^{3/2}} \bigr]' = \frac{24 \beta x - 18 \beta^2 - 8 \beta + 2}{4(1 - \beta)^{5/2} (1 + 3 \beta)^{3/2}} = \frac{12 \beta x - 9 \beta^2 - 4 \beta + 1}{2 (1 - \beta)^{5/2} (1 + 3 \beta)^{3/2}}.
\end{equation}
The fact that $[0,1]\ni x\mapsto\ff(x,0,\dots,0)\in\R$ is continuous, the fact that for all $x \in (\alpha, \beta)$ it holds that
\begin{equation}
\textstyle 0 < \frac{1}{(1 - x)^{3/2} (1 + 3x)^{5/2}} < \frac{1}{(1 - \beta)^{3/2} (1 + 3 \alpha)^{5/2}},
\end{equation}
and the fact that for all $x \in (\beta, 1)$ it holds that 
\begin{equation}
\textstyle \frac{3 \beta^2 - 4 \beta + 1}{2 (1 - \beta)^{5/2} (1 + 3 \beta)^{3/2}} < \frac{12 \beta x - 9 \beta^2 - 4 \beta + 1}{2 (1 - \beta)^{5/2} (1 + 3 \beta)^{3/2}} < \frac{-9 \beta^2 + 8 \beta + 1}{2 (1 - \beta)^{5/2} (1 + 3 \beta)^{3/2}}
\end{equation}
\hence assure that $[0, 1]^d \ni x \mapsto \ff(x) \in \R$ is Lipschitz continuous. This implies for all $\kappa \in \R$ that $[\scra, \scrb]^d \ni x \mapsto \scrf_{\kappa}(x) \in \R$ is Lipschitz continuous. Next \nobs that \cref{eqn:cor:uncountably_many_realizations_delta:realization} and \cref{eqn:cor:uncountably_many_realizations_delta:risk} show for all $g \in C([\scra, \scrb]^d, \R)$, $c \in \R$, $\width \in \N$, $\theta = (\theta_1, \dots, \theta_{\dimension_{\width}}), \vartheta = (\vartheta_1, \dots, \vartheta_{\dimension_{\width}}) \in \R^{\dimension_{\width}}$ with
\begin{equation}
\textstyle \bigwedge_{i = 1}^{\width} \bigl(\bigwedge_{j = (i - 1)d + 1}^{i d} ([\vartheta_{d \width + \width + i} = c \theta_{d \width + \width + i}] \wedge [\vartheta_j = \theta_j] \wedge [\vartheta_{d \width + i} = \theta_{d \width + i}] \wedge [\vartheta_{\dimension_{\width}} = c \theta_{\dimension_{\width}}])\bigr)
\end{equation}
that
$\cR_{c g, \width}(\vartheta) = \abs{c}^2 \cR_{g, \width}(\theta)$.
Moreover, \nobs that for every $c\in\R\backslash\{0\}$ the parameter transformation in the preceding display is an invertible linear map and its inverse is obtained by replacing $c$ with $c^{-1}$. This shows that this transformation preserves local minimum points, multiplies the infimum of the risk by $\abs{c}^2$, and maps distinct realization functions to distinct realization functions. Combining this with \cref{eqn:cor:uncountably_many_realizations_of_non_global_local_minimas_with_varying_delta:scrr} establishes for all $\kappa \in (0, \infty)$ that
\begin{multline}
\textstyle \bigl\{v \in C([\scra, \scrb]^d, \R) \colon \bigl[\Exists \theta \in \{\vartheta \in \R^{\dimension_2} \colon v = \cN_{2}^{\vartheta}\}, \eps \in (0, \infty) \colon \\
\textstyle \cR_{\scrf_{\kappa}, 2}(\theta) = \inf_{\vartheta \in [-\eps, \eps]^{\dimension_2}} \cR_{\scrf_{\kappa}, 2}(\theta + \vartheta) > \kappa + \inf_{\vartheta \in \R^{\dimension_2}} \cR_{\scrf_{\kappa}, 2}(\vartheta)\bigr]\bigr\}
\end{multline}
is an uncountable set. Next \nobs that every two-neuron parameter vector can be embedded into the parameter space with $\width>2$ hidden neurons by retaining the first two neurons and the outer bias and by setting the inner weights and outer weights of every added neuron to zero and its inner bias to $-1$. The compactness of $[\scra,\scrb]^d$ ensures that these added neurons remain inactive throughout a sufficiently small parameter neighborhood. This implies that the embedding preserves the realization function and that every local minimum point of the two-neuron risk yields a local minimum point of the wider-network risk with the same risk value. Together with the identity embedding in the case $\width=2$, this proves for all $g \in C([\scra, \scrb]^d, \R)$, $\width \in \N \cap (1, \infty)$ that
\begin{equation}
\textstyle \bigcup_{\vartheta \in \R^{\dimension_2}} \{\cR_{g, 2}(\vartheta)\} \subseteq \bigcup_{\vartheta \in \R^{\dimension_{\width}}} \{\cR_{g, \width}(\vartheta)\}, \quad \inf_{\vartheta \in \R^{\dimension_2}} \cR_{g, 2}(\vartheta) \ge \inf_{\vartheta \in \R^{\dimension_{\width}}} \cR_{g, \width}(\vartheta), 
\end{equation}
\begin{equation}
\begin{split}
\textstyle \andqq \bigl\{v \in C([\scra, \scrb]^d, \R) \colon \bigl[\Exists \theta \in \{\vartheta \in \R^{\dimension_2} \colon & \textstyle v = \cN_{2}^{\vartheta}\}, \eps \in (0, \infty) \colon \\
& \textstyle \cR_{g, 2}(\theta) = \inf_{\vartheta \in [-\eps, \eps]^{\dimension_2}} \cR_{g, 2}(\theta + \vartheta)\bigr]\bigr\} \\
\textstyle \subseteq \bigl\{v \in C([\scra, \scrb]^d, \R) \colon \bigl[\Exists \theta \in \{\vartheta \in \R^{\dimension_{\width}} \colon & \textstyle v = \cN_{\width}^{\vartheta}\}, \eps \in (0, \infty) \colon \\
& \textstyle \cR_{g, \width}(\theta) = \inf_{\vartheta \in [-\eps, \eps]^{\dimension_{\width}}} \cR_{g, \width}(\theta + \vartheta)\bigr]\bigr\}.
\end{split}
\end{equation}
Combining this with the scaled risk gap above and the fact that every $\kappa\in[\max\{1,\delta\},\infty)$ satisfies $\kappa\geq\delta$ and $\kappa>0$ demonstrates for all $g \in \cF$, $\width \in \N \cap (1, \infty)$ that
\begin{multline}
\textstyle \bigl\{v \in C([\scra, \scrb]^d, \R) \colon \bigl[\Exists \theta \in \{\vartheta \in \R^{\dimension_{\width}} \colon v = \cN_{\width}^{\vartheta}\}, \eps \in (0, \infty) \colon \\
\textstyle \cR_{g, \width}(\theta) = \inf_{\vartheta \in [-\eps, \eps]^{\dimension_{\width}}} \cR_{g, \width}(\theta + \vartheta) > \delta + \inf_{\vartheta \in \R^{\dimension_{\width}}} \cR_{g, \width}(\vartheta)\bigr]\bigr\}
\end{multline}
is an uncountable set. Moreover, \nobs that $f$ is not the zero function, since its normalized restriction to $(\alpha,\beta)$ has strictly positive derivative. The definition of $\scrf_{\kappa}$ and the fact that $[\max\{1,\delta\},\infty)$ is uncountable \hence show that $\cF$ is uncountable. Combining this with the fact that every $\scrf_{\kappa}$ is Lipschitz continuous establishes \cref{eqn:cor:uncountably_many_realizations_of_non_global_local_minimas_with_varying_delta:res}.
\end{cproof}

\begin{lstlisting}[language=Python, caption={{\sc Python} source code used to create \cref{fig:local:minimum:points}.}, label=list:python]
import numpy as np
import matplotlib.pyplot as plt
import scipy.integrate as integrate
import math

# hidden neurons, initializations, learning rate, and seed
𝔥 = 4; M = 50; eta = 1/20; np.random.seed(0)

# target function f : [0,1] -> \R
def f(x):
    if (x <= 1/3): return (math.sqrt(3)*(2*x - 1)/8)
    elif (x >= 2/3): return (6*x**2-17/2*x+109/36)
    else: return ((3*x**2-1)/(4*math.sqrt(1-x)*math.sqrt((1+3*x)**3)))
    
# realization funct. g : (1/3,2/3) x [0,1] -> \R of non-global loc. min. pt.
def g(q, x):
    return (-math.sqrt(1-q)/(4*math.sqrt(1+3*q))
            +(np.maximum(x-q,0))/(2*math.sqrt((1-q)**3)*math.sqrt(1+3*q)))
    
# realization function associated to an ANN
def realization(theta, x):
    return (np.dot(theta[2*𝔥:3*𝔥],np.maximum(theta[0:𝔥]*x+theta[𝔥:2*𝔥],0))
            +theta[3*𝔥])

# generalized gradient function gen_grad : \R^{3*𝔥+1} -> \R^{3*𝔥+1}
def gen_grad(theta):
    g = np.ones(3*𝔥+1)
    for j in range(𝔥):
        g[j] = 2*theta[2*𝔥+j]*integrate.quad(lambda x:
                x*(realization(theta,x)-f(x))*(lambda y: 1 if y>0 else 0)
                (theta[j]*x+theta[𝔥+j]), 0, 1, limit=100)[0]
        g[𝔥+j] = 2*theta[2*𝔥+j]*integrate.quad(lambda x:
                (realization(theta,x)-f(x))*(lambda y: 1 if y>0 else 0)
                (theta[j]*x+theta[𝔥+j]), 0, 1, limit=100)[0]
        g[2*𝔥+j] = 2*integrate.quad(lambda x:
                np.maximum(theta[j]*x+theta[𝔥+j],0)
                *(realization(theta,x)-f(x)), 0, 1, limit=100)[0]
    g[3*𝔥] = 2*integrate.quad(lambda x:
            realization(theta,x)-f(x), 0, 1, limit=100)[0]
    return g

# vectorized versions of f, g, and the realization function
vec_f = np.vectorize(f)
vec_g = np.vectorize(g)
vec_realization = np.vectorize(realization, excluded=['theta'])

# initialize each ANN with the Xavier initialization
parameters = np.random.normal(scale=math.sqrt(2/(𝔥+1)), size=(3*𝔥+1)*M)
for i in range(M):
    parameters[i*(3*𝔥+1)+𝔥:i*(3*𝔥+1)+2*𝔥] = parameters[i*(3*𝔥+1)+3*𝔥] = 0
theta = np.reshape(parameters, (M,3*𝔥+1))

# train each ANN with the GD optimization method
for i in range(M):
    while(True):
        grad = gen_grad(theta[i])
        if np.linalg.norm(grad,np.inf) < 1e-4:
            break
        theta[i] = theta[i] - eta * grad

# divide the trained ANNs into groups with a similar realization function
groups = []
for i in range(M):
    for group in groups:
        if(integrate.quad(lambda x: (realization(theta[i],x)
                -realization(group[1],x))**2, 0, 1)[0] < 1e-4):
            group.append(theta[i])
            break
    else:
        groups.append([integrate.quad(lambda x:
            (realization(theta[i],x)-f(x))**2, 0, 1)[0], theta[i]])
groups = sorted(groups, key=lambda g: g[0])

# plot the simulation with one representative realization function per group
plt.rcParams["axes.prop_cycle"] = plt.cycler('color', ['tab:blue','tab:red',
    'tab:orange','tab:green','tab:purple'])
plt.figure()
plt.margins(x=0)
plt.xticks(np.arange(0,1.1,0.1))
plt.ylim((-0.25,0.55))
x = np.linspace(0,1,10000)
plt.plot(x, vec_f(x), label="risk"+r'$=$'+"0 (target function)", zorder=M+1)
for group in groups:
    plt.plot(x, vec_realization(theta=group[1],x=x),
             label="risk"+r'$\approx$'+format(group[0], 'f'))
plt.legend()
plt.savefig("simulation.pdf", bbox_inches="tight")

# plot ten realization funct. of the considered non-global loc. min. points
plt.figure()
plt.margins(x=0)
plt.xticks(np.arange(0,1.1,0.1))
plt.ylim((-0.25,0.55))
plt.plot(x, vec_f(x), label="risk"+r'$=$'+"0 (target function)", zorder=11)
for q in np.linspace(0.3333333334, 0.6666666666, 10):
    plt.plot(x, vec_g(q,x), color="tab:orange")
plt.plot(x, vec_g(2/3,x), label="risk"+r'$\approx$'+format(integrate.quad(
    lambda x: (g(2/3,x)-f(x))**2, 0, 1)[0], 'f'), color="tab:orange")
plt.legend()
plt.savefig("non-global_local_minimum_points.pdf", bbox_inches="tight")
\end{lstlisting}

\section{On finitely many realization functions of critical points}
\label{sec:on_finitely_many_realization_functions}

In this section we prove in \cref{cor:finite_number_of_realizations} in \cref{subsec:main_positive_results} below in the special situation where the target function $f \colon [\scra, \scrb] \to \R$ is continuous and piecewise polynomial and where the considered shallow \ReLU\ \ANNs\ have one neuron on each layer that there exist only finitely many different realization functions $\functionSNN^{\theta} \colon [\scra, \allowbreak \scrb] \allowbreak \to \R$, $\theta \in \R^4$, of all critical points of the risk function $\riskR \colon \R^4 \to \R$ in the sense that there exist only finitely many different realization functions $\functionSNN^{\theta} \colon [\scra, \allowbreak \scrb] \allowbreak \to \R$, $\theta \in \R^4$, of zeros of the generalized gradient function $ \gradientG \colon \R^4 \to \R^4 $. In \cref{cor:convergence_to_global_minima} in \cref{subsec:main_positive_results} we extend \cref{cor:finite_number_of_realizations} by using \cite[Item~$($v$)$ in Theorem~1.1]{Adrian2021GradientFlows} and \cite[Theorem~1.2]{eberle2021existence} to establish that in the training of such \ANNs\ we have that the risk of every non-divergent \GF\ trajectory converges to the risk of a global minimum point provided that the initial risk is sufficiently small.

The remainder of this section is organized in the following way. In \cref{setting2} in \cref{subsec:ANNs_one-dim_input_and_hidden_layer} below we present our mathematical setup of shallow \ANNs\ with one neuron on each layer, in the elementary result in \cref{lemma:positiv:case_1} in \cref{subsec:critical_points_with_constant_realization_functions} below we analyze critical points with constant realization functions, in the elementary result in \cref{lemma:positiv:case_2} in \cref{subsec:critical_points_with_affine_linear_realization_functions} below we analyze critical points with affine linear realization functions, in \cref{lemma:positiv:case_3} in \cref{subsec:critical_points_with_non-decreasing_non-affine_linear_realization_functions} below we analyze critical points with right-active ($\theta_1 > 0$) non-affine linear realization functions, and in \cref{lemma:positiv:case_4} in \cref{subsec:critical_points_with_non-increasing_non-affine_linear_realization_functions} below we analyze critical points with left-active ($\theta_1 < 0$) non-affine linear realization functions.

In \cref{subsec:main_positive_results} we combine \cref{lemma:positiv:case_1}, \cref{lemma:positiv:case_2}, \cref{lemma:positiv:case_3}, and \cref{lemma:positiv:case_4} to establish \cref{cor:finite_number_of_realizations}.

\subsection{ANNs with one-dimensional input and hidden layer}
\label{subsec:ANNs_one-dim_input_and_hidden_layer}

\begin{setting}\label{setting2}
Let $n \in \N$, $\fx_0, \fx_1, \dots, \fx_n, \scra, \scrb \in \R$ satisfy $\scra \allowbreak = \fx_0 \allowbreak < \fx_1 \allowbreak < \dots \allowbreak < \fx_n \allowbreak = \scrb$, let $f \allowbreak \in C([\scra, \allowbreak \scrb], \allowbreak \R)$ satisfy  for all $j \in \{1, 2, \dots, n\}$ that $f|_{[\fx_{j - 1}, \fx_j]}$ is a polynomial, let $\Rect_r \allowbreak \colon \R \allowbreak \to \R$, $r \in \N \cup \{\infty\}$, satisfy for all $x \in \R$ that $(\cup_{r \in \N} \{\Rect_r \}) \subseteq C^1 (\R, \R)$, $\Rect_\infty(x) = \max \{x, 0\}$, $\sup_{r \in \N} \sup_{y \in [- \abs{x}, \abs{x}]} \abs{(\Rect_r)'(y)} < \infty$, and
\begin{equation}
\limsup\nolimits_{r \to \infty} (\abs{\Rect_r(x) - \Rect_\infty(x)} + \abs{(\Rect_r)'(x) - \indicator{(0, \infty)}(x)}) = 0,
\end{equation}
for every $r \in \allowbreak \N \allowbreak \cup \{\infty\}$ let $\realapprox{}{r} \allowbreak = (\realapprox{\theta}{r})_{\theta \in \R^4} \colon \R^4 \allowbreak \to C([\scra, \allowbreak \scrb], \allowbreak \R)$ and $\cR_r \colon \R^4 \to \R$ satisfy for all $\theta \allowbreak = (\theta_1, \allowbreak \dots, \allowbreak \theta_4) \allowbreak \in \R^4$, $x \allowbreak \in [\scra,\scrb]$ that $\realapprox{\theta}{r}(x) = \theta_4 + \theta_3 [\Rect_r (\theta_2 + \theta_1 x)]$ and $\cR_r(\theta) \allowbreak = \int_{\scra}^{\scrb} \abs{\realapprox{\theta}{r}(y) \allowbreak - f(y)}^2 \, \d y$, for every $\theta = (\theta_1, \dots, \theta_4) \in \R^4$ let $I^{\theta} \allowbreak \subseteq \R$ satisfy $I^{\theta} \allowbreak = \{x \allowbreak \in [\scra, \allowbreak \scrb] \colon \theta_2 \allowbreak + \theta_1 x \allowbreak > 0\}$, and let $\gradientG \colon \R^{4} \to \R^{4}$ satisfy for all $\theta \in \{\vartheta \in \R^{4} \colon ((\nabla \cR_r)(\vartheta))_{r\in \N} \text{ is convergent}\}$ that $\gradientG(\theta) \allowbreak = \lim_{r \to \infty} (\nabla \cR_r)(\theta)$.
\end{setting}

\subsection{Critical points with constant realization functions}
\label{subsec:critical_points_with_constant_realization_functions}

\cfclear
\begin{lemma}\label{lemma:positiv:case_1}
Assume \cref{setting:snn0} and assume $\mu([\scra, \allowbreak \scrb]^d) \allowbreak > 0$. Then
\begin{equation}\label{eqn:lemma:positiv:case_1}
\#(\{v \in C([\scra, \scrb]^d, \R) \colon ( \Exists \theta \in \gradientG^{-1}(\{0\}) \colon [\Forall x, y \in [\scra, \scrb]^d \colon v(x) = \realapprox{\theta}{\infty}(y)])\}) = 1.
\end{equation}
\end{lemma}
\begin{cproof}{lemma:positiv:case_1}
Throughout this proof let $V \subseteq \R^{\dimension}$ satisfy
\begin{equation}
\textstyle V = \{\theta \in \gradientG^{-1}(\{0\}) \colon [\Forall \allowbreak x, \allowbreak y \in [\scra, \scrb]^d \colon \allowbreak \realapprox{\theta}{\infty}(x) = \realapprox{\theta}{\infty}(y)]\}.
\end{equation}
\Nobs that \cref{prop:limit:lr} and \cref{prop:loss:diff:vc} ensure that
\begin{enumerate}[label = (\roman*)]
\item it holds for all $\theta \in V$ that $[\scra, \scrb]^d \ni x \mapsto \realapprox{\theta}{\infty}(x) \in \R$ is constant,

\item it holds for all $\theta = (\theta_1, \dots, \theta_{\dimension}) \in V$ that $\R \ni t \mapsto \riskRR(\theta_1, \theta_2, \dots, \theta_{\dimension - 1}, t) \in \R$ is differentiable, and

\item it holds for all $\theta \in V$ that
\begin{equation}
\textstyle (\frac{\partial}{\partial \theta_{\dimension}} \riskRR)(\theta) \textstyle = 2 \int_{[\scra, \scrb]^d} (\realapprox{\theta}{\infty}(x) - f(x)) \, \mu(\d x) = 0.
\end{equation}
\end{enumerate}
\Hence for all $\theta \in V$, $x \in [\scra, \scrb]^d$ that
\begin{equation}\label{eqn:lemma:positiv:case_1:pre_final}
\begin{split}
\textstyle \realapprox{\theta}{\infty}(x) & \textstyle = \frac{1}{\mu([\scra, \scrb]^d)} [\int_{[\scra, \scrb]^d} (\realapprox{\theta}{\infty}(y) - f(y)) \, \mu(\d y) + \int_{[\scra, \scrb]^d} f(y) \, \mu(\d y)] \\
& \textstyle = \frac{1}{\mu([\scra, \scrb]^d)} [\int_{[\scra, \scrb]^d} f(y) \, \mu(\d y)].
\end{split}
\end{equation}
This and \cref{prop:limit:lr} ensure that for all $\theta \allowbreak \in \R^{\dimension}$ with
\begin{equation}
\textstyle \Forall i \allowbreak \in \{1, \allowbreak 2, \allowbreak \dots, \allowbreak \width\} \colon (\norm{\fw_i^{\theta}} \allowbreak + \abs{\fv_i^{\theta}} \allowbreak = 0) \wedge (\fb_i^{\theta} \allowbreak = -1)
\end{equation}
and $\fc^{\theta} = [\mu([\scra,\scrb]^d)]^{-1} \allowbreak \int_{[\scra, \scrb]^d} f(x) \, \mu(\d x)$ it holds that $\theta \in V$. \Hence that $V \neq \emptyset$. Combining this with \cref{eqn:lemma:positiv:case_1:pre_final} establishes \cref{eqn:lemma:positiv:case_1}.
\end{cproof}

\subsection{Critical points with non-constant affine linear realization functions}
\label{subsec:critical_points_with_affine_linear_realization_functions}

\cfclear
\begin{lemma}\label{lemma:positiv:case_2}
Assume \cref{setting:snn0} and assume $d = 1 = \width$ and $\mu = \lambda$. Then
\begin{multline}\label{eqn:lemma:positiv:case_2}
\#(\{v \in C([\scra, \scrb], \R) \colon (\Exists \theta = (\theta_1, \dots, \theta_{4}) \in \gradientG^{-1}(\{0\}) \colon \\
\textstyle [(\theta_1 \theta_{3} \neq 0) \wedge ((\scra, \scrb) \subseteq I_1^{\theta}) \wedge (v = (\realapprox{\theta}{\infty})|_{[\scra,\scrb]})])\}) \leq 1.
\end{multline}
\end{lemma}
\begin{cproof}{lemma:positiv:case_2}
Throughout this proof let $V \subseteq \R^{4}$ satisfy
\begin{equation}
V = \{\theta = (\theta_1, \dots, \theta_{4}) \in \gradientG^{-1}(\{0\}) \colon [(\theta_1 \theta_{3} \neq 0) \wedge ((\scra, \scrb) \subseteq I_1^{\theta})]\}.
\end{equation}
\Nobs that \cref{prop:limit:lr}, \cref{prop:loss:diff:vc}, and \cref{prop:loss:diff:wb} show that
\begin{enumerate}[label = (\roman*)]
\item it holds for all $\theta \in V$ that $\riskRR$ is differentiable at $\theta$ and

\item it holds for all $\theta \in V$ that $(\nabla \riskRR)(\theta) = \gradientG(\theta) = 0$.
\end{enumerate}
\Hence that for all $\theta = (\theta_1, \dots, \theta_{4}) \in V$ it holds that
\begin{equation}
\begin{split}
\textstyle (\frac{\partial}{\partial \theta_{4}} \riskRR)(\theta) & \textstyle = 2 \int_{\scra}^{\scrb} (\realapprox{\theta}{\infty}(x) - f(x)) \, \d x = 0 \\
& \textstyle = 2 \int_{\scra}^{\scrb} (\theta_1 x + \theta_{2})(\realapprox{\theta}{\infty}(x) - f(x)) \, \d x = (\frac{\partial}{\partial \theta_{3}} \riskRR)(\theta).
\end{split}
\end{equation}
This implies that for all $\theta \in V$ it holds that
\begin{equation}
\textstyle \int_{\scra}^{\scrb} (\realapprox{\theta}{\infty}(x) - f(x)) \, \d x = 0 = \int_\scra^{\scrb} x (\realapprox{\theta}{\infty}(x) - f(x)) \, \d x.
\end{equation}
Combining this with the fact that for all $\theta = (\theta_1, \dots, \theta_{4}) \in V$, $x \in [\scra, \scrb]$ it holds that $\realapprox{\theta}{\infty}(x) = \theta_1 \theta_{3} x + \allowbreak (\theta_{4} + \theta_{3} \theta_{2})$ ensures that for all $\theta = (\theta_1, \dots, \theta_4) \in V$ it holds that
\begin{equation}
\textstyle \int_{\scra}^{\scrb} (\theta_1 \theta_3 x + (\theta_4 + \theta_3 \theta_2) - f(x)) \, \d x = \theta_1 \theta_3 \bigl[\frac{\scrb^2 - \scra^2}{2}\bigr] + (\theta_4 + \theta_3 \theta_2) (\scrb - \scra) - \int_{\scra}^{\scrb} f(x) \, \d x = 0
\end{equation}
and
\begin{equation}
\textstyle \int_\scra^{\scrb} x(\theta_1 \theta_3 x + (\theta_4 + \theta_3 \theta_2) - f(x)) \, \d x = \theta_1 \theta_3 \bigl[ \frac{\scrb^3 - \scra^3}{3} \bigr] + (\theta_4 + \theta_3 \theta_2) \bigl[ \frac{\scrb^2 - \scra^2}{2} \bigr] - \int_{\scra}^{\scrb} x f(x) \, \d x = 0.
\end{equation}
The fact that 
\begin{equation}
\textstyle \bigl[\frac{(\scrb^2 - \scra^2)^2}{4}\bigr] - [\scrb - \scra] \bigl[\frac{\scrb^3 - \scra^3}{3}\bigr] = - \frac{(\scrb - \scra)^2}{12} [4 (\scrb^2 + \scrb \scra + \scra^2) - 3 (\scrb + \scra)^2] = - \frac{(\scrb - \scra)^4}{12} \neq 0
\end{equation}
\hence shows that there exist $c_1, c_2 \in \R$ which satisfy for all $\theta = (\theta_1, \dots, \theta_4) \in V$ that 
\begin{equation}\label{eqn:lemma:positiv:case_2:c_1andc_2}
\textstyle \theta_1 \theta_3 = c_1 \qqandqq \theta_4 + \theta_3 \theta_2 = c_2.
\end{equation}
\Nobs that \cref{eqn:lemma:positiv:case_2:c_1andc_2} implies that for all $\theta \in V$, $x \in [\scra, \scrb]$ it holds that $\realapprox{\theta}{\infty}(x) = c_1 x + c_2$. This establishes \cref{eqn:lemma:positiv:case_2}.
\end{cproof}

\subsection[Critical points with right-active non-affine linear realization functions]{Critical points with right-active non-affine linear realization\\ functions}
\label{subsec:critical_points_with_non-decreasing_non-affine_linear_realization_functions}

\cfclear
\begin{lemma}\label{lemma:int_eq_finiteness}
Let $\fn \in \N$, $\scrx_0, \scrx_1, \dots, \scrx_{\fn} \in \R$ satisfy $0 = \scrx_0 \allowbreak < \scrx_1 \allowbreak < \dots \allowbreak < \scrx_{\fn} \allowbreak = 1$, let $\scrf \allowbreak \in C([0, \allowbreak 1], \allowbreak \R)$ satisfy for all $j \in \{1, 2, \dots, \fn\}$ that $\scrf|_{[\scrx_{j - 1}, \scrx_j]}$ is a polynomial, and let $j \allowbreak \in \{1, \allowbreak 2, \allowbreak \dots, \allowbreak \fn\}$. Then
\begin{multline}\label{eqn:lemma:int_eq_finiteness:set}
\bigl\{q \in [\scrx_{j - 1}, \scrx_j] \backslash \{0, 1\} \colon \\
\textstyle \bigl[\bigl(\frac{(1 - q)^2}{6 q} \int_0^{q} \scrf(x) \, \d x = \int_{q}^1 (\frac{q + 2}{3} - x) \scrf(x) \, \d x\bigr) \wedge \bigl(\int_0^1 \scrf(x) \, \d x \neq \frac{1}{q} \int_0^{q} \scrf(x) \, \d x\bigr)\bigr]\bigr\}
\end{multline}
is a finite set.
\end{lemma}
\begin{cproof}{lemma:int_eq_finiteness}
Throughout this proof we use the convention $\operatorname{deg}(0)=-\infty$, assume w.l.o.g.\ that the set in \cref{eqn:lemma:int_eq_finiteness:set} is non-empty, and let $d \in \N_0$, $q \in [\scrx_{j-1}, \allowbreak \scrx_j] \backslash \{0,\allowbreak 1\}$ satisfy
\begin{equation}\label{eqn:lemma:int_eq_finiteness:conditions}
\textstyle \frac{(1 - q)^2}{6 q} \int_0^{q} \scrf(x) \, \d x = \int_{q}^1 (\frac{q + 2}{3} - x) \scrf(x) \, \d x, \qquad\int_0^1 \scrf(x) \, \d x \neq \frac{1}{q} \int_0^{q} \scrf(x) \, \d x,
\end{equation}
and $d=\max\{0,\operatorname{deg}(\scrf|_{[\scrx_{j-1},\scrx_j]})\}$. In the following we distinguish between the case $d = 0$ and the case $d > 0$. First we prove \cref{eqn:lemma:int_eq_finiteness:set} in the case $d = 0$. \Nobs that the assumption that $d = 0$ implies that there exists $\fc \in \R$ which satisfies for all $x \in [\scrx_{j-1}, \scrx_{j}]$ that $\scrf(x) = \fc$. \Hence that $(1-q)^2 [\int_0^q \scrf(x) \, \d x] = (1-q)^2 [\int_0^{\scrx_{j-1}} \scrf(x) \, \d x \allowbreak + \int_{\scrx_{j-1}}^q \fc \, \d x]$ and $2q \int_q^1 (q+2 - 3x) \scrf(x) \, \d x = 2q [\int_q^{\scrx_j} (q+2 - 3x) \fc \, \d x + \int_{\scrx_j}^1 (q+2 - 3x) \scrf(x) \, \d x]$. Combining this with 
\cref{eqn:lemma:int_eq_finiteness:conditions} assures that
\begin{equation}
\textstyle (1-q)^2 \bigl[\int_0^{\scrx_{j-1}} \scrf(x) \, \d x + \int_{\scrx_{j-1}}^q \fc \, \d x\bigr] = 2q \bigl[\int_q^{\scrx_j} (q+2 - 3x) \fc \, \d x + \int_{\scrx_j}^1 (q+2 - 3x) \scrf(x) \, \d x\bigr].
\end{equation}
The fact that
\begin{equation}
\begin{split}
& \textstyle (1-q)^2 \bigl[\int_0^{\scrx_{j-1}} \scrf(x) \, \d x + \int_{\scrx_{j-1}}^q \fc \, \d x\bigr] \\
& \textstyle = (q^2 - 2q + 1) \bigl[ \int_0^{\scrx_{j-1}} \scrf(x) \, \d x + \fc q - \fc \scrx_{j - 1} \bigr] = \fc q^3 + q^2 \bigl[ \int_0^{\scrx_{j-1}} \scrf(x) \, \d x - \fc \scrx_{j - 1} - 2 \fc \bigr] \\
& \textstyle \quad + q \bigl[-2 \int_{0}^{\scrx_{j - 1}} \scrf(x) \, \d x + 2 \fc \scrx_{j - 1} + \fc\bigr] + \bigl[ \int_0^{\scrx_{j-1}} \scrf(x) \, \d x - \fc \scrx_{j - 1}\bigr]
\end{split}
\end{equation}
and the fact that
\begin{equation}
\begin{split}
& \textstyle 2q \bigl[\int_q^{\scrx_j} (q+2 - 3x) \fc \, \d x + \int_{\scrx_j}^1 (q+2 - 3x) \scrf(x) \, \d x\bigr] \\
& \textstyle = 2q \bigl[\fc \scrx_j q - \fc q^2 + 2 \fc \scrx_j - 2 \fc q - \frac{3}{2} \fc [\scrx_j]^2 + \frac{3}{2} \fc q^2 + q \int_{\scrx_j}^1 \scrf(x) \, \d x + 2 \int_{\scrx_j}^1 \scrf(x) \, \d x \\
& \textstyle \quad - 3 \int_{\scrx_j}^1 x \scrf(x) \, \d x\bigr] = \fc q^3 + q^2 [2 \fc \scrx_j - 4 \fc + 2 \int_{\scrx_j}^1 \scrf(x) \, \d x] \\
& \textstyle \quad + q [4 \fc \scrx_j - 3 \fc [\scrx_j]^2 + 4 \int_{\scrx_j}^1 \scrf(x) \, \d x - 6 \int_{\scrx_j}^1 x \scrf(x) \, \d x]
\end{split}
\end{equation}
hence show that
\begin{multline}\label{eqn:cor:finite_number_of_realizations:quadratic_eq_on_q}
\textstyle q^2 \bigl[\int_0^{\scrx_{j-1}} \scrf(x) \, \d x - \fc \scrx_{j-1} + 2 \fc - 2 \fc \scrx_j - 2 \int_{\scrx_j}^1 \scrf(x) \, \d x\bigr] \\
\textstyle + q \bigl[\fc + 2 \fc \scrx_{j-1} - 4 \fc \scrx_j + 3 \fc [\scrx_j]^2 - 2 \int_0^{\scrx_{j-1}} \scrf(x) \, \d x - 4 \int_{\scrx_j}^1 \scrf(x) \, \d x \\
\textstyle + 6 \int_{\scrx_j}^1 x \scrf(x) \, \d x\bigr] + \bigl[\int_0^{\scrx_{j-1}} \scrf(x) \, \d x - \fc \scrx_{j-1}\bigr] = 0.
\end{multline}
Next \nobs that \cref{eqn:lemma:int_eq_finiteness:conditions} ensures that
\begin{equation}
\begin{split}
& \textstyle \int_0^1 \scrf(x) \, \d x - \frac{1}{q} \int_0^q \scrf(x) \, \d x \\
& \textstyle = \int_0^{\scrx_{j-1}} \scrf(x) \, \d x + \int_{\scrx_{j-1}}^{\scrx_j} \fc \, \d x + \int_{\scrx_j}^1 \scrf(x) \, \d x - \frac{1}{q} \bigl[\int_0^{\scrx_{j-1}} \scrf(x) \, \d x + \int_{\scrx_{j-1}}^q \fc \, \d x\bigr] \\
& \textstyle = \bigl(1 - \frac{1}{q}\bigr) \bigl[\int_0^{\scrx_{j-1}} \scrf(x) \, \d x - \fc \scrx_{j - 1}\bigr] + \bigl[\int_{\scrx_j}^1 \scrf(x) \, \d x - (\fc - \fc \scrx_j)\bigr] \neq 0.
\end{split}
\end{equation}
This shows that
\begin{equation}
\textstyle \bigl[\int_0^{\scrx_{j-1}} \scrf(x) \, \d x - \fc \scrx_{j-1} + 2 (\fc - \fc \scrx_j - \int_{\scrx_j}^1 \scrf(x) \, \d x)\bigr]^2 + \bigl[\int_0^{\scrx_{j-1}} \scrf(x) \, \d x - \fc \scrx_{j-1}\bigr]^2 \neq 0.
\end{equation}
Combining this with \cref{eqn:cor:finite_number_of_realizations:quadratic_eq_on_q} establishes \cref{eqn:lemma:int_eq_finiteness:set} in the case $d = 0$. In the next step we prove \cref{eqn:lemma:int_eq_finiteness:set} in the case $d > 0$. \Nobs that the assumption that $d > 0$ assures that there exist $a \in \R \backslash \{0\}$ and $Q \colon \R \to \R$ which satisfy for all $x \in [\scrx_{j-1}, \scrx_{j}]$ that $Q$ is a polynomial with $\deg(Q) \le d - 1$ and $\scrf(x) = a x^{d} + Q(x)$. \Hence that 
\begin{equation}
\textstyle (1-q)^2 \bigl[\int_0^q \scrf(x) \, \d x\bigr] = (1-q)^2 \bigl[\int_0^{\scrx_{j-1}} \scrf(x) \, \d x + \int_{\scrx_{j-1}}^q (ax^{d} + Q(x)) \, \d x\bigr]
\end{equation}
and
\begin{equation}
\textstyle 2q \int_q^1 (q+2 - 3x) \scrf(x) \, \d x = 2q \bigl[\int_q^{\scrx_j} (q+2 - 3x) (ax^{d} + Q(x)) \, \d x + \int_{\scrx_j}^1 (q+2 - 3x) \scrf(x) \, \d x\bigr].
\end{equation} 
Combining this with 
\cref{eqn:lemma:int_eq_finiteness:conditions} shows that
\begin{equation}\label{eqn:lemma:int_eq_finiteness:eq_poly}
\begin{split}
& \textstyle (1-q)^2 \bigl[\int_0^{\scrx_{j-1}} \scrf(x) \, \d x + \int_{\scrx_{j-1}}^q (ax^{d} + Q(x)) \, \d x\bigr] \\
& \textstyle = 2q \bigl[\int_q^{\scrx_j} (q+2 - 3x) (ax^{d} + Q(x)) \, \d x + \int_{\scrx_j}^1 (q+2 - 3x) \scrf(x) \, \d x\bigr].
\end{split}
\end{equation}
Next \nobs that the fact that $\operatorname{deg} (Q) \le d - 1$ demonstrates that there exist polynomials $\scrP \colon \R \to \R$ and $\fP \colon \R \to \R$ which do not depend on $q$ and satisfy $\max\{\operatorname{deg} (\scrP), \operatorname{deg} (\fP)\} \le d + 2$,
\begin{equation}
\begin{split}
& \textstyle (1-q)^2 \bigl[\int_0^{\scrx_{j-1}} \scrf(x) \, \d x + \int_{\scrx_{j-1}}^q (ax^{d} + Q(x)) \, \d x\bigr] \\
& \textstyle = (q^2 - 2q + 1) \bigl[ \frac{a}{d + 1} (q^{d + 1} - [\scrx_{j - 1}]^{d + 1}) + \int_0^{\scrx_{j - 1}} \scrf(x) \, \d x + \int_{\scrx_{j - 1}}^q Q(x) \, \d x \bigr] \\
& \textstyle = q^{d + 3} \bigl[ \frac{a}{d + 1}\bigr] + \scrP(q),
\end{split}
\end{equation}
and
\begin{equation}
\begin{split}
& \textstyle 2q \bigl[\int_q^{\scrx_j} (q+2 - 3x) (ax^{d} + Q(x)) \, \d x + \int_{\scrx_j}^1 (q+2 - 3x) \scrf(x) \, \d x\bigr] \\
& \textstyle = 2q \bigl[ \int_q^{\scrx_j} (q - 3x) a x^{d} \, \d x + \int_q^{\scrx_j} (q - 3x) Q(x) \, \d x + 2 \int_q^{\scrx_j} (a x^{d} + Q(x)) \, \d x \\
& \textstyle \quad  + \int_{\scrx_j}^1 (q+2 - 3x) \scrf(x) \, \d x \bigr] = 2q \bigl[ \frac{a q}{d + 1} ([\scrx_j]^{d + 1} - q^{d + 1}) - \frac{3a}{d + 2} ([\scrx_j]^{d + 2} - q^{d + 2}) \\
& \textstyle \quad + \int_q^{\scrx_j} (q - 3x) Q(x) \, \d x + 2 \int_q^{\scrx_j} (a x^{d} + Q(x)) \, \d x + \int_{\scrx_j}^1 (q+2 - 3x) \scrf(x) \, \d x\bigr] \\
& \textstyle = q^{d + 3} \bigl[ \frac{-2a}{d + 1} + \frac{6a}{d + 2} \bigr] + \fP(q) = q^{d + 3} \bigl[ \frac{2a (2 d + 1)}{(d + 1)(d + 2)} \bigr] + \fP(q).
\end{split}
\end{equation}
Combining this with \cref{eqn:lemma:int_eq_finiteness:eq_poly} ensures that
\begin{equation}\label{eqn:cor:finite_number_of_realizations:eq_on_q_d+3}
\textstyle q^{d + 3} \bigl[\frac{2a(2d + 1)}{(d + 1)(d + 2)} - \frac{a}{d + 1}\bigr] + \fP(q) - \scrP(q) = q^{d + 3} \bigl[\frac{3a d}{(d + 1)(d + 2)}\bigr] + \fP(q) - \scrP(q) = 0.
\end{equation}
The fact that $a \neq 0$, the fact that $d > 0$, and the fact that $\max\{\operatorname{deg}(\scrP), \allowbreak \operatorname{deg}(\fP)\} \le d + 2$ \hence establish \cref{eqn:lemma:int_eq_finiteness:set} in the case $d > 0$.
\end{cproof}

\begin{corollary}\label{cor:int_eq_finiteness}
Let $\fn \in \N$, $\scrx_0, \scrx_1, \dots, \scrx_{\fn} \in \R$ satisfy $0 = \scrx_0 \allowbreak < \scrx_1 \allowbreak < \dots \allowbreak < \scrx_{\fn} \allowbreak = 1$ and let $\scrf \allowbreak \in C([0, \allowbreak 1], \allowbreak \R)$ satisfy for all $j \in \{1, 2, \dots, \fn\}$ that $\scrf|_{[\scrx_{j - 1}, \scrx_j]}$ is a polynomial. Then
\begin{equation}\label{eqn:cor:int_eq_finiteness:set}
\bigl\{ \textstyle q \in (0, 1) \colon \bigl[\bigl(\frac{(1 - q)^2}{6 q} \int_0^q \scrf(x) \, \d x = \int_q^1 (\frac{q + 2}{3} - x) \scrf(x) \, \d x\bigr) \wedge \bigl(\int_0^1 \scrf(x) \, \d x \neq \frac{1}{q} \int_0^q \scrf(x) \, \d x\bigr) \bigr]\bigr\}
\end{equation}
is a finite set.
\end{corollary}
\begin{cproof}{cor:int_eq_finiteness}
Throughout this proof for every $j \in \{1, 2, \dots, \fn\}$ let $S_j \subseteq (0, 1)$ satisfy
\begin{multline}\label{eqn:cor:int_eq_finiteness:set_pieces}
S_j = \bigl\{q \in [\scrx_{j - 1}, \scrx_j] \backslash \{0, 1\} \colon \\
\textstyle \bigl[\bigl(\frac{(1 - q)^2}{6 q} \int_0^{q} \scrf(x) \, \d x = \int_{q}^1 (\frac{q + 2}{3} - x) \scrf(x) \, \d x\bigr) \wedge \bigl(\int_0^1 \scrf(x) \, \d x \neq \frac{1}{q} \int_0^{q} \scrf(x) \, \d x\bigr)\bigr]\bigr\}.
\end{multline}
\Nobs that \cref{lemma:int_eq_finiteness} and \cref{eqn:cor:int_eq_finiteness:set_pieces} ensure for all $j \in \{1, 2, \dots, \fn\}$ that $S_j$ is a finite set. The fact that
\begin{multline}
\textstyle \cup_{j \in \{1, 2, \dots, \fn\}} S_{j} = \bigl\{ \textstyle q \in (0, 1) \colon \\
\textstyle \bigl[\bigl(\frac{(1 - q)^2}{6 q} \int_0^{q} \scrf(x) \, \d x = \int_{q}^1 (\frac{q + 2}{3} - x) \scrf(x) \, \d x\bigr) \wedge \bigl(\int_0^1 \scrf(x) \, \d x \neq \frac{1}{q} \int_0^{q} \scrf(x) \, \d x\bigr) \bigr]\bigr\}
\end{multline}
\hence establishes \cref{eqn:cor:int_eq_finiteness:set}.
\end{cproof}

\cfclear
\begin{lemma}\label{lemma:positiv:case_3}
Assume \cref{setting2}. Then
\begin{multline}\label{eqn:lemma:positiv:case_3}
\textstyle \{v \in C([\scra, \scrb], \R) \colon (\Exists \theta = (\theta_1, \dots, \theta_{4}) \in \gradientG^{-1}(\{0\}) \colon \\
\textstyle [(\theta_1 > 0 \neq \theta_{3}) \wedge (\scra < -(\theta_{2} / \theta_1) < \scrb) \wedge (v = \functionSNN^{\theta})])\}
\end{multline} 
is a finite set.
\end{lemma}
\begin{cproof}{lemma:positiv:case_3}
Throughout this proof let $\scrf \in C([0, 1], \R)$ satisfy for all $x \allowbreak \in [0, \allowbreak 1]$ that $\scrf(x) \allowbreak = f(x(\scrb - \scra) \allowbreak + \scra)$, let $\phi \allowbreak = (\phi_1, \allowbreak \dots, \allowbreak \phi_4) \allowbreak \in \gradientG^{-1}(\{0\})$ satisfy $\phi_1 \allowbreak > 0 \allowbreak \neq \phi_{3}$ and $\scra \allowbreak < - \phi_2 / \phi_1 \allowbreak < \scrb$, let $\theta \allowbreak = (w, \allowbreak b, \allowbreak v, \allowbreak c) \allowbreak \in \R^{4}$ satisfy $w \allowbreak = \phi_1(\scrb - \scra)$, $b \allowbreak= \phi_1 \scra \allowbreak + \phi_2$, $v \allowbreak = \phi_3$, and $c \allowbreak = \phi_4$, and let $q \allowbreak \in \R$ satisfy $q \allowbreak = - b / w$. \Nobs that the fact that $\abs{\phi_1} + \abs{\phi_2} > 0$, \cref{prop:loss:diff:vc:item1} in \cref{prop:loss:diff:vc}, and \cref{prop:loss:diff:wb:item1} in \cref{prop:loss:diff:wb} show that $\riskR$ is differentiable at $\phi$. The fact that $\gradientG(\phi) = 0$, \cref{prop:limit:lr}, \cref{prop:loss:diff:vc}, \cref{prop:loss:diff:wb}, and the integral transformation theorem \hence ensure that
\begin{equation}
\begin{aligned}[b]
& \textstyle (\frac{\partial}{\partial \phi_1} \riskR)(\phi) \\
& \textstyle = 2 \phi_3 \int_{I^{\phi}} x (\functionSNN^{\phi}(x) - f(x)) \, \d x = 2 \phi_3 \int_{I^{\phi}} x (\phi_4 + \phi_3 \max\{\phi_2 + \phi_1 x, 0\} - f(x)) \, \d x \\
& \textstyle = 2(\scrb - \scra) \phi_3 \int_{q}^1 ((\scrb - \scra) x + \scra) \bigl(\phi_4 + \phi_3 [\phi_2 + \phi_1 ((\scrb - \scra) x + \scra)] - f((\scrb - \scra) x + \scra)\bigr) \, \d x \\
& \textstyle = 2(\scrb - \scra) v \int_{q}^1 ((\scrb - \scra) x + \scra) (c + v w (x - q) - \scrf(x)) \, \d x = 0,
\end{aligned}
\end{equation}
\begin{equation}
\begin{split}
& \textstyle (\frac{\partial}{\partial \phi_2} \riskR)(\phi) \\
& \textstyle = 2 \phi_3 \int_{I^{\phi}} (\functionSNN^{\phi}(x) - f(x)) \, \d x = 2 \phi_3 \int_{I^{\phi}} (\phi_4 + \phi_3 \max\{\phi_2 + \phi_1 x, 0\} - f(x)) \, \d x \\
& \textstyle = 2 (\scrb - \scra) \phi_3 \int_{q}^1 (\phi_4 + \phi_3 [\phi_2 + \phi_1 ((\scrb - \scra) x + \scra)] - f((\scrb - \scra) x + \scra)) \, \d x \\
& \textstyle = 2 (\scrb - \scra) v \int_{q}^1 (c + v w (x - q) - \scrf(x)) \, \d x = 0,
\end{split}
\end{equation}
\begin{equation}
\begin{split}
& \textstyle (\frac{\partial}{\partial \phi_3} \riskR)(\phi) = 2 \int_{I^{\phi}} (\phi_2 + \phi_1 x)(\functionSNN^{\phi} (x) - f(x)) \, \d x \\
& \textstyle = 2 \int_{I^{\phi}} (\phi_2 + \phi_1 x)(\phi_4 + \phi_3 \max\{\phi_2 + \phi_1 x, 0\} - f(x)) \, \d x \\
& \textstyle = 2 (\scrb - \scra) \int_{q}^1 (\phi_1 ((\scrb - \scra) x + \scra) + \phi_2)(\phi_4 + \phi_3 [\phi_2 + \phi_1 ((\scrb - \scra) x + \scra)] \\
& \textstyle \quad - f((\scrb - \scra) x + \scra)) \, \d x = 2 (\scrb - \scra) \int_{q}^1 (w x + b) (c + v w (x - q) - \scrf(x)) \, \d x = 0,
\end{split}
\end{equation}
and
\begin{equation}
\begin{split}
& \textstyle (\frac{\partial}{\partial \phi_{4}} \riskR)(\phi) \\ 
& \textstyle = 2 \int_{\scra}^{\scrb} (\functionSNN^{\phi} (x) - f(x)) \, \d x = 2 \int_{\scra}^{\scrb} (\phi_4 + \phi_3 \max\{\phi_2 + \phi_1 x, 0\} - f(x)) \, \d x \\
& \textstyle = 2 (\scrb - \scra) \int_0^1 (\phi_4 + \phi_3 \max\{\phi_2 + \phi_1 ((\scrb - \scra) x + \scra), 0\} - f((\scrb - \scra) x + \scra)) \, \d x \\ 
& \textstyle = 2 (\scrb - \scra) \int_0^1 (c + v w \max\{x - q, 0\} - \scrf(x)) \, \d x  = 0.
\end{split}
\end{equation}
\Hence that
\begin{equation}\label{eqn:cor:finite_number_of_realizations:int_equations_w>0}
\textstyle \int_0^q (c - \scrf(x)) \, \d x = \int_q^1 (c + vw (x-q) - \scrf(x)) \, \d x = \int_q^1 x(c + vw (x-q) - \scrf(x)) \, \d x = 0.
\end{equation}
This implies that
\begin{equation}
\textstyle cq - \int_0^q \scrf(x) \, \d x = 0, \qquad c (1 - q) + v w \bigl[\frac{1 - q^2}{2} - q (1 - q) \bigr] - \int_q^1 \scrf(x) \, \d x = 0,
\end{equation}
\begin{equation}
\textstyle \text{and} \qquad c \bigl[ \frac{1 - q^2}{2} \bigr] + v w \bigl[ \frac{1 - q^3}{3} - \frac{q(1 - q^2)}{2} \bigr] - \int_q^1 x \scrf(x) \, \d x = 0.
\end{equation}
\Hence that
\begin{equation}\label{eqn:cor:finite_number_of_realizations:int_eq_c_vw}
\textstyle c = \frac{1}{q} \!\pb{\int_0^q \scrf(x) \, \d x}\!, \qquad vw = \frac{2}{(1-q)^2} \! \pb{\int_0^1 \scrf(x) \, \d x - \frac{1}{q} \int_0^q \scrf(x) \, \d x}\!,\\
\end{equation}
\begin{equation}\label{eqn:cor:finite_number_of_realizations:int_eq_q}
\textstyle \text{and} \qquad \frac{(1-q)^2}{6q} \! \pb{\int_0^q \scrf(x) \, \d x} \! = \int_q^1 \! \pa{\frac{q+2}{3} - x} \! \scrf(x) \, \d x.
\end{equation}
Moreover, \nobs that the fact that $\phi_1 \phi_3 \neq 0$ assures that $v w \neq 0$. Combining this with \cref{eqn:cor:finite_number_of_realizations:int_eq_c_vw} shows that $\int_0^1 \scrf(x) \, \d x \neq \frac{1}{q} \int_0^q \scrf(x) \, \d x$. \cref{cor:int_eq_finiteness}, \cref{eqn:cor:finite_number_of_realizations:int_eq_c_vw}, and \cref{eqn:cor:finite_number_of_realizations:int_eq_q} \hence establish \cref{eqn:lemma:positiv:case_3}.
\end{cproof}

\subsection{Critical points with left-active non-affine linear realization functions}
\label{subsec:critical_points_with_non-increasing_non-affine_linear_realization_functions}

\cfclear
\begin{lemma}\label{lemma:positiv:case_4}
Assume \cref{setting2}. Then
\begin{multline}\label{eqn:lemma:positiv:case_4:finite_set}
\textstyle \{v \in C([\scra, \scrb], \R) \colon (\Exists \theta = (\theta_1, \dots, \theta_{4}) \in \gradientG^{-1}(\{0\}) \colon \\
\textstyle [(\theta_1 < 0 \neq \theta_{3}) \wedge (\scra < -(\theta_{2} / \theta_1) < \scrb) \wedge (v = \functionSNN^{\theta})])\}
\end{multline} 
is a finite set.
\end{lemma}
\begin{cproof}{lemma:positiv:case_4}
Throughout this proof let $\scrf, \ff \in C([0, 1], \R)$ satisfy for all $x \allowbreak \in [0, \allowbreak 1]$ that $\scrf(x) \allowbreak = f(x(\scrb - \scra) \allowbreak + \scra)$ and $\ff(x) \allowbreak = \scrf(1 - x)$, let $\phi \allowbreak = (\phi_1, \allowbreak \dots, \allowbreak \phi_4) \allowbreak \in \gradientG^{-1}(\{0\})$ satisfy $\phi_1 \allowbreak < 0 \allowbreak \neq \phi_3$ and $\scra \allowbreak < - \phi_2 / \phi_1 \allowbreak < \scrb$, let $\theta = (w, \allowbreak b, \allowbreak v, \allowbreak c) \allowbreak \in \R^{4}$ satisfy $w \allowbreak = \phi_1(\scrb - \scra)$, $b \allowbreak = \phi_1 \scra \allowbreak + \phi_2$, $v \allowbreak = \phi_3$, and $c \allowbreak = \phi_4$, and let $q, \allowbreak \scrq \allowbreak \in \R$ satisfy $q \allowbreak = -\frac{b}{w}$ and $\scrq \allowbreak = 1 - q$. \Nobs that the fact that $\abs{\phi_1} \allowbreak + \abs{\phi_2} \allowbreak > 0$, \cref{prop:loss:diff:vc:item1} in \cref{prop:loss:diff:vc}, and \cref{prop:loss:diff:wb:item1} in \cref{prop:loss:diff:wb} show that $\riskR$ is differentiable at $\phi$. The fact that $\gradientG(\phi) = 0$, \cref{prop:limit:lr}, \cref{prop:loss:diff:vc}, \cref{prop:loss:diff:wb}, and the integral transformation theorem \hence ensure that
\begin{equation}
\begin{aligned}[b]
& \textstyle (\frac{\partial}{\partial \phi_1} \riskR)(\phi) \\
& \textstyle = 2 \phi_3 \int_{I^{\phi}} x (\functionSNN^{\phi} (x) - f(x)) \, \d x = 2 \phi_3 \int_{I^{\phi}} x (\phi_4 + \phi_3 \max\{\phi_2 + \phi_1 x, 0\} - f(x)) \, \d x \\
& \textstyle = 2(\scrb - \scra) \phi_3 \int_{0}^{q} ((\scrb - \scra) x + \scra) \bigl(\phi_4 + \phi_3 [\phi_2 + \phi_1 ((\scrb - \scra) x + \scra)] - f((\scrb - \scra) x + \scra)\bigr) \, \d x \\
& \textstyle = 2(\scrb - \scra) v \int_{0}^{q} ((\scrb - \scra) x + \scra) (c + v w (x - q) - \scrf(x)) \, \d x = 0,
\end{aligned}
\end{equation}
\begin{equation}
\begin{split}
& \textstyle (\frac{\partial}{\partial \phi_2} \riskR)(\phi) \\
& \textstyle = 2 \phi_3 \int_{I^{\phi}} (\functionSNN^{\phi} (x) - f(x)) \, \d x = 2 \phi_3 \int_{I^{\phi}} (\phi_4 + \phi_3 \max\{\phi_2 + \phi_1 x, 0\} - f(x)) \, \d x \\
& \textstyle = 2 (\scrb - \scra) \phi_3 \int_{0}^{q} (\phi_4 + \phi_3 [\phi_2 + \phi_1 ((\scrb - \scra) x + \scra)] - f((\scrb - \scra) x + \scra)) \, \d x \\
& \textstyle = 2 (\scrb - \scra) v \int_{0}^{q} (c + v w (x - q) - \scrf(x)) \, \d x = 0,
\end{split}
\end{equation}
\begin{equation}
\begin{split}
& \textstyle (\frac{\partial}{\partial \phi_3} \riskR)(\phi) = 2 \int_{I^{\phi}} (\phi_2 + \phi_1 x)(\functionSNN^{\phi} (x) - f(x)) \, \d x \\
& \textstyle = 2 \int_{I^{\phi}} (\phi_2 + \phi_1 x)(\phi_4 + \phi_3 \max\{\phi_2 + \phi_1 x, 0\} - f(x)) \, \d x \\
& \textstyle = 2 (\scrb - \scra) \int_{0}^{q} (\phi_2 + \phi_1 ((\scrb - \scra) x + \scra)) (\phi_4 + \phi_3 [\phi_2 + \phi_1 ((\scrb - \scra) x + \scra)] \\
& \textstyle \quad - f((\scrb - \scra) x + \scra)) \, \d x = 2 (\scrb - \scra) \int_{0}^{q} (w x + b) (c + v w (x - q) - \scrf(x)) \, \d x = 0,
\end{split}
\end{equation}
and
\begin{equation}
\begin{split}
& \textstyle (\frac{\partial}{\partial \phi_{4}} \riskR)(\phi) \\
& \textstyle = 2 \int_{\scra}^{\scrb} (\functionSNN^{\phi} (x) - f(x)) \, \d x = 2 \int_{\scra}^{\scrb} (\phi_4 + \phi_3 \max\{\phi_2 + \phi_1 x, 0\} - f(x)) \, \d x \\ 
& \textstyle = 2 (\scrb - \scra) \int_0^1 (\phi_4 + \phi_3 \max\{\phi_2 + \phi_1 ((\scrb - \scra) x + \scra), 0\} - f((\scrb - \scra) x + \scra)) \, \d x \\
& \textstyle = 2 (\scrb - \scra) \int_0^1 (c + v w \min\{x - q, 0\} - \scrf(x)) \, \d x = 0.
\end{split}
\end{equation}
\Hence that
\begin{equation}
\textstyle \int_q^1 (c - \scrf(x)) \, \d x = \int_0^q (c + vw (x-q) - \scrf(x)) \, \d x = \int_0^q x(c + vw (x-q) - \scrf(x)) \, \d x = 0.
\end{equation}
This implies that
\begin{equation}
\textstyle c(1 - q) - \int_q^1 \scrf(x) \, \d x = 0, \qquad  c q + vw \bigl[ \frac{q^2}{2} - q^2 \bigr] - \int_0^q \scrf(x) \, \d x = 0,
\end{equation}
\begin{equation}
\textstyle \text{and} \qquad c \bigl[\frac{q^2}{2}\bigr] + vw \bigl[ \frac{q^3}{3} - \frac{q^3}{2} \bigr] - \int_0^q x \scrf(x) \, \d x = 0.
\end{equation}
\Hence that
\begin{equation}\label{eqn:lemma:positiv:case_4_c_vw}
\textstyle c = \frac{1}{1 - q} \bigl[\int_q^1 \scrf(x) \, \d x\bigr], \qquad vw = \frac{2}{q^2} \bigl[\frac{1}{1 - q} \int_q^1 \scrf(x) \, \d x - \int_0^1 \scrf(x) \, \d x\bigr],
\end{equation}
\begin{equation}\label{eqn:lemma:positiv:case_4:int_eq}
\textstyle \text{and} \qquad \frac{q^2}{6(1 - q)} \bigl[\int_q^1 \scrf(x) \, \d x\bigr] = \int_0^q \! \pa{x - \frac{q}{3}} \! \scrf(x) \, \d x.
\textstyle 
\end{equation}
Combining this with the integral transformation theorem ensures that
\begin{equation}\label{eqn:lemma:positiv:case_4_vw_trans}
\textstyle v w = \frac{2}{(1 - \scrq)^2}\!\pb{\frac{1}{\scrq} \int_0^{\scrq} \ff(y) \, \d y - \int_0^1 \ff(y) \, \d y} \qandq \frac{(1 - \scrq)^2}{6 \scrq} \int_0^{\scrq} \ff(y) \, \d y = \int_{\scrq}^1 (\frac{\scrq + 2}{3} - y) \ff(y) \, \d y.
\end{equation}
The fact that $\phi_1 \phi_3 \neq 0$, \cref{cor:int_eq_finiteness}, \cref{eqn:lemma:positiv:case_4_c_vw}, and \cref{eqn:lemma:positiv:case_4:int_eq} \hence establish \cref{eqn:lemma:positiv:case_4:finite_set}.
\end{cproof}

\subsection{On finitely many realization functions of critical points}
\label{subsec:main_positive_results}
\cfclear
\begin{corollary}\label{cor:finite_number_of_realizations}
Assume \cref{setting2}. Then 
\begin{equation}\label{eqn:cor:finite_number_of_realizations}
\{v \in C([\scra, \scrb], \R) \colon ( \Exists \theta \in \gradientG^{-1} (\{0\}) \colon v = \functionSNN^{\theta} )\}
\end{equation}
is a finite set \cfout.
\end{corollary}
\begin{cproof}{cor:finite_number_of_realizations}
Throughout this proof let $V\subseteq C([\scra,\scrb],\R)$ denote the set in \cref{eqn:cor:finite_number_of_realizations} and let $V_0,V_1,V_2\subseteq V$ satisfy
\begin{gather}\label{eqn:cor:finite_number_of_realizations:constant}
\textstyle V_0 = \{v \in V\colon v\text{ is constant}\}, \\ \label{eqn:cor:finite_number_of_realizations:affine}
\textstyle V_1 = \{v \in V\colon v\text{ is affine linear and not constant}\}, \\ \label{eqn:cor:finite_number_of_realizations:non-affine}
\textstyle \andqq V_2 = \{v\in V\colon v\text{ is not affine linear}\}.
\end{gather}
\Nobs that the fact that the Lebesgue measure of $[\scra, \allowbreak \scrb]$ is $\scrb \allowbreak - \scra \allowbreak > 0$ and \cref{lemma:positiv:case_1} (applied with $d = 1 = \width$, $\mu = \lambda$ in the notation of \cref{lemma:positiv:case_1}) show that $V_0$ is finite. Moreover, \nobs that for all $\theta \allowbreak = (\theta_1, \allowbreak \dots, \allowbreak \theta_4) \allowbreak \in \gradientG^{-1}(\{0\})$ with $\functionSNN^{\theta} \allowbreak \in V_1$ it holds that $\theta_1 \theta_3 \allowbreak \neq 0$ and $(\scra, \allowbreak \scrb) \allowbreak \subseteq I^{\theta}$. \cref{lemma:positiv:case_2} (applied with $d = 1 = \width$, $\mu = \lambda$ in the notation of \cref{lemma:positiv:case_2}) \hence ensures that $V_1$ has at most one element. Next \nobs that for all $\theta \allowbreak = (\theta_1, \allowbreak \dots, \allowbreak \theta_4) \allowbreak \in \gradientG^{-1}(\{0\})$ with $\functionSNN^{\theta} \allowbreak \in V_2$ it holds that $\theta_1 \theta_3 \allowbreak \neq 0$ and $\scra \allowbreak < - \theta_2/\theta_1 \allowbreak < \scrb$. \cref{lemma:positiv:case_3,lemma:positiv:case_4} therefore show that $V_2$ is finite. Combining this and the fact that $V_0$ and $V_1$ are finite with the fact that $V = V_0 \allowbreak \cup \allowbreak V_1 \allowbreak \cup V_2$ establishes that $V$ is finite.
\end{cproof}

\begin{corollary}\label{cor:convergence_to_global_minima}
Assume \cref{setting2}. Then 
\begin{enumerate} [label = (\roman*)]
\item
\label{item1:cor:convergence_to_global_minima} it holds that $\{v \in C([\scra, \scrb], \R) \colon ( \Exists \theta \in \gradientG^{-1} (\{0\}) \colon v = \functionSNN^{\theta} )\}$ is a finite set,

\item 
\label{item2:cor:convergence_to_global_minima} it holds that
\begin{multline}
\bigl\{ v \in C([\scra, \scrb], \R) \colon \bigl[ \Exists \theta \in \{ \vartheta \in \R^4 \colon v = \functionSNN^{\vartheta} \}, \varepsilon \in (0,\infty) \colon \\
\riskR (\theta) = \inf\nolimits_{ \vartheta \in [-\varepsilon,\varepsilon]^4 } \riskR(\theta + \vartheta) \bigr] \bigr\}
\end{multline}
is a finite set, and

\item
\label{item3:cor:convergence_to_global_minima} there exists $\eps \in (0, \infty)$ which satisfies for all $\Theta = (\Theta_t)_{t \in [0, \infty)} \colon [0, \infty) \to \R^4$ with $\riskR(\Theta_0) \allowbreak \le \eps \allowbreak + \inf_{\vartheta \in \R^4} \riskR(\vartheta)$, $\Forall t \allowbreak \in [0, \allowbreak \infty) \allowbreak \colon \Theta_t \allowbreak = \Theta_0 \allowbreak - \int_0^t \gradientG(\Theta_s) \, \d s$, and $\liminf_{t \to \infty} \allowbreak \norm{\Theta_t} < \infty$ that
\begin{equation}
\textstyle \limsup_{t \to \infty} \riskR(\Theta_t) = \inf_{\vartheta \in \R^4} \riskR(\vartheta).
\end{equation}
\end{enumerate}
\cfout[.]
\end{corollary}
\begin{cproof}{cor:convergence_to_global_minima}
\Nobs that \cref{cor:finite_number_of_realizations} establishes \cref{item1:cor:convergence_to_global_minima}. Moreover, \nobs that \allowbreak \cref{prop:local:minima:gradient} implies that 
\begin{multline}
\textstyle \bigl\{ v \in C([\scra, \scrb], \R) \colon \bigl[ \Exists \theta \in \{ \vartheta \in \R^4 \colon v = \functionSNN^{\vartheta} \}, \eps \in (0,\infty) \colon \\
\textstyle \riskR (\theta) = \inf\nolimits_{ \vartheta \in [-\varepsilon,\varepsilon]^4 } \riskR (\theta + \vartheta)\bigr]\bigr\} \subseteq \{v \in C([\scra, \scrb], \R) \colon ( \Exists \theta \in \gradientG^{-1} (\{0\}) \colon v = \functionSNN^{\theta})\}.
\end{multline}
Combining this with \cref{item1:cor:convergence_to_global_minima} establishes \cref{item2:cor:convergence_to_global_minima}. Next \nobs that \cref{item1:cor:convergence_to_global_minima} implies that
\begin{equation}
\textstyle \{v \in C([\scra, \scrb], \R) \colon (\Exists \theta \in \gradientG^{-1}(\{0\}) \cap (\riskR)^{-1}\bigl((\inf_{\vartheta \in \R^{4}} \riskR(\vartheta), \infty)\bigr) \colon v = \functionSNN^{\theta})\}
\end{equation}
is a finite set. \Hence that there exists a finite set $S \subseteq \R^{4}$ which satisfies 
\begin{equation}
\textstyle \{v \in \R \colon (\Exists \theta \in \gradientG^{-1}(\{0\}) \cap (\riskR)^{-1}\bigl((\inf_{\vartheta \in \R^{4}} \riskR(\vartheta), \infty)\bigr) \colon v = \riskR(\theta))\} = \allowbreak \cup_{\theta \in S} \allowbreak \{\riskR(\theta)\}. 
\end{equation}
Let $\eps\in(0,\infty)$ satisfy $\eps=1$ in the case $S=\emptyset$ and $\eps< [\min_{\theta\in S}\riskR(\theta)]-\inf_{\vartheta\in\R^4}\riskR(\vartheta)$ in the case $S\neq\emptyset$ and let $\Theta \allowbreak = (\Theta_t)_{t \in [0, \infty)} \allowbreak \colon [0, \allowbreak \infty) \allowbreak \to \R^4$ satisfy
\begin{gather}
\textstyle \liminf_{t \to \infty} \norm{\Theta_t} < \infty, \qquad \Forall t \in [0, \infty) \colon \Theta_t = \Theta_0 - \int_0^t \gradientG(\Theta_s) \, \d s, \\
\textstyle \andqq \riskR(\Theta_0) \le \eps + \inf_{\vartheta \in \R^{4}} \riskR(\vartheta)\ifnocf.
\end{gather}
\cfload[.]\Nobs that the piecewise polynomial assumption on $f$, the fact that the density of the input data is constant, and \cite[Theorem~1.2]{eberle2021existence} ensure that $\sup_{t \in [0, \infty)}\norm{\Theta_t} < \infty$ \cfload. The fact that for all $\theta \in \gradientG^{-1}(\{0\}) \allowbreak \cap (\riskR)^{-1}((\inf_{\vartheta \in \R^{4}} \riskR(\vartheta), \allowbreak \infty))$ it holds that
\begin{equation}
\textstyle \riskR(\Theta_0) \le \eps + \inf_{\vartheta \in \R^{4}} \riskR(\vartheta) < \riskR(\theta)
\end{equation}
and \cite[Item $($v$)$ in Theorem~1.1]{Adrian2021GradientFlows} \hence show that 
\begin{equation}
\textstyle \limsup_{t \to \infty} \riskR(\Theta_t) = \inf_{\vartheta \in \R^{4}} \riskR(\vartheta).
\end{equation}
This establishes \cref{item3:cor:convergence_to_global_minima}.
\end{cproof}

\begin{remark}[\resname{Intuition for the construction of the target function in \cref{theorem:negative} resulting in a risk function with infinitely many realization functions of non-global local minimum points}]\label{remark:intuition}
Consider the situation of \cref{setting2} with $\scra=0=1-\scrb$ and with an arbitrary continuous target function $f$ instead of a piecewise polynomial target function. \Nobs that \cref{prop:local:minima:gradient} shows that every local minimum point of $\riskR\colon\R^4\to\R$ is a zero of the generalized gradient $\gradientG\colon\R^4\to\R^4$. For a fixed $q\in(0,1)$ and a fixed sign of $\theta_1$, the identities in \cref{eqn:cor:finite_number_of_realizations:int_eq_c_vw,eqn:lemma:positiv:case_4_c_vw} show that there exists at most one realization function of a zero $\theta=(\theta_1,-q\theta_1,\theta_3,\theta_4)\in\gradientG^{-1}(\{0\})$ with $\theta_1\theta_3\neq0$. In particular, fixing the breakpoint alone leaves two possible activation orientations.

The possible breakpoints satisfy the integral equation in \cref{eqn:cor:finite_number_of_realizations:int_eq_q} in the right-active case ($\theta_1 \allowbreak > 0$) and the integral equation in \cref{eqn:lemma:positiv:case_4:int_eq} in the left-active case ($\theta_1 \allowbreak < 0$), together with the corresponding nonzero-slope condition in \cref{eqn:cor:finite_number_of_realizations:int_eq_c_vw,eqn:lemma:positiv:case_4_c_vw}. For a continuous piecewise polynomial target function, \cref{cor:int_eq_finiteness} and the transformation $x \allowbreak \mapsto 1 - x$ imply that each of these equations has only finitely many solutions satisfying that nonzero-slope condition. The nonzero-slope condition is essential, since a constant target function satisfies the integral equations for every $q\in(0,1)$ but gives only a constant realization function through the displayed moment identities.

When the piecewise polynomial assumption is relaxed, the target function in \cref{0_eq:f} is Lipschitz continuous on $[0,1]$ and smooth on $(\alpha,\beta)$, and the right-active integral equation in \cref{eqn:cor:finite_number_of_realizations:int_eq_q} has a solution with nonzero slope for every $q \allowbreak \in (\alpha, \allowbreak \beta)$. The residual identities in \cref{lemma:zero_derivative} express precisely this property. The manifold and Hessian arguments in \cref{lemma:set_of_critical_points_submanifold}, \cref{lemma:rank&Hessian}, and \cref{lemma:uncountably_many_realizations_of_critical_points} then establish local minimality, and the two-neuron comparison in \cref{cor:uncountably_many_realizations_of_non_global_local_minimas_with_fixed_delta} establishes nonglobality when the hidden layer has more than one neuron. This explains the construction underlying \cref{theorem:negative} and \cref{cor:uncountably_many_realizations_of_non_global_local_minimas_with_varying_delta}.
\end{remark}

\subsubsection*{Acknowledgements}
This work has been partially funded by the National Natural Science Foundation of China (NSFC) under grant number W\-2\-5\-3\-1\-0\-1\-0. This work has also been partially funded by the European Union (ERC, MONTECARLO, 1\-0\-1\-0\-4\-5\-8\-1\-1). The views and the opinions expressed in this work are however those of the authors only and do not necessarily reflect those of the European Union or the European Research Council (ERC). Neither the European Union nor the granting authority can be held responsible for them. We also gratefully acknowledge the Cluster of Excellence EXC 2044/2-390685587, Mathematics Münster: Dynamics-Geometry-Structure funded by the Deutsche Forschungsgemeinschaft (DFG, German Research Foundation). The first author dedicates this paper to the memory of his grandmother, Khosiyat Ibragimova, who passed away during the preparation of the initial version of this paper.


\subsubsection*{Use of large language models}
{\sc GPT 6 Astra Pro} has helped us to correct several mathematical inaccuracies and minor gaps in this work and {\sc GPT 6 Astra Pro} has also helped us to improve the formulation of some of the sentences in this work. The entire article (including all statements and all arguments in this article) has been carefully written by the authors and the authors take full responsibility for each of the sentences/statements made in the paper.

\bibliographystyle{acm}
\bibliography{bibfile}

\end{document}